%% file: nlshape.tex
\documentclass[a4paper]{article}

\include{preamble}

\author{ Volker Schulz\thanks{Universitaet Trier, D-54286 Trier, Germany; Email: volker.schulz@uni-trier.de, schusterm@uni-trier.de, vollmann@uni-trier.de} \and Matthias Schuster\footnotemark[1] \and Christian Vollmann\footnotemark[1] }

\title{\textbf{Shape optimization for interface identification in nonlocal models}}
\date{}

\begin{document}

\maketitle

\small
\textbf{Abstract.} 
Shape optimization methods have been proven useful for identifying interfaces in models governed by partial differential equations. Here we consider a class of shape optimization problems constrained by nonlocal equations which involve interface--dependent kernels. We derive a novel shape derivative associated to the nonlocal system model and solve the problem by established numerical techniques.
\normalsize

\section{Introduction}
%
Many physical relations and data-based coherences cannot satisfactorily be described by classical differential equations. Often they inherently possess some features, which are not purely local. In this regard, mathematical models which are governed by nonlocal operators enrich our modeling spectrum and present useful alternative as well as supplemental approaches. That is why they appear in a large variety of applications including among others, anomalous or fractional diffusion \cite{fractional_anomaDiff1,fractional_anomaDiff2,fractlapl}, peridynamics \cite{peridym1,Silling_birthperidym,peridym2, coupling_nl_to_l_Marta}, image processing \cite{imageproc1,imageproc3,imageproc2}, cardiology \cite{fractional_cardiology}, machine learning \cite{machinelearning}, as well as finance and jump processes \cite{finance_levy,jumpprocess,nonlocaldirichletform,finance_jumpproc, nonsym_ker_1}. Nonlocal operators are integral operators allowing for interactions between two distinct points in space. The nonlocal models investigated in this paper involve kernels that are not necessarily symmetric and which are assumed to have a finite range of nonlocal interactions; see, e.g, \cite{Du_book,tian_diss, wellposedness, nonsym_ker_1} and the references therein.

Not only the problem itself but also various optimization problems involving nonlocal models of this type are treated in literature. For example matching-type problems are treated in \cite{optcontrol_forcingterm,optcontrol_fractional1,optcontrol_parameter_kernel} to identify system parameters such as the forcing term or a scalar diffusion parameter. The control variable is typically modeled to be an element of a suitable function space. 
Moreover, nonlocal interface problems have become popular in recent years\cite{D'EliaGlusaInterface,RodriguezGimperleinInterface,EnergyCoupling}.
However, shape optimization techniques applied to nonlocal models can hardly be found in literature. For instance, the articles \cite{shape_nonloc_fractional_eigvals,shape_nonloc_fractional_eigvals_ball_is_minimizer,shape_nonlocal_fractional_min_functionalofeigvals} deal with minimizing (functions of) eigenvalues of the fractional Laplacian with respect to the domain of interest. Also, in \cite{shape_nonloc_energy_functional_fractional1,shape_nonloc_energy_functional_fractional3} the energy functional related to fractional equations is minimized. In \cite{shape_nonloc_energy_functional_general} a functional involving a more general kernel is considered. All of the aforementioned papers are of theoretical nature only. To the best of our knowledge, shape optimization problems involving nonlocal constraint equations with truncated kernels and numerical methods for solving such problems cannot yet be found in literature.

Instead, shape optimization problems which are constrained by partial differential equations appear in many fields of application \cite{shape_aero,shape_image_seg,shape_siebenborn1,shape_cellular} and particularly for inverse problems where the parameter to be estimated, e.g., the diffusivity in a heat equation model, is assumed to be defined piecewise on certain subdomains. Given a rough picture of the configuration, shape optimization techniques can be successfully applied to identify the detailed shape of these subdomains \cite{shape_volker1, shape_volker2, shape_volker3,Welker_diss}. 

In this paper we transfer the problem of parameter identification into a nonlocal regime. Here, the parameter of interest is given by the kernel which describes the nonlocal model. We assume that this kernel is defined piecewise with respect to a given partition $\{\Omega_i\}_{i}$ of the domain of interest $\Omega$. 
Thereby, the state of such a nonlocal model depends on the interfaces between the respective subdomains $\Omega_i$. Under the assumption that we know the rough setting but are lacking in details, we can apply the techniques developed in the aforementioned shape optimization papers to identify these interfaces from a given measured state.

For this purpose we formulate a shape optimization problem which is constrained by an interface--dependent nonlocal convection--diffusion model. Here, we do not aim at investigating conceptual improvements of existing shape optimization algorithms. On the contrary, we want to study the applicability of established methods for problems of this type. Thus this paper can be regarded as a feasibility study where we set a focus on the numerical implementation. 

The realization of this plan basically requires two ingredients both of which are worked out here. First, we define a reasonable interface--dependent nonlocal model and provide a finite element code which discretizes a variational formulation thereof. Second, we need to derive the shape derivative of the corresponding nonlocal bilinear form which is then implemented into an overall shape optimization algorithm. 

This leads to the following organization of the present paper. In Section \ref{sec:problemformulation} we formulate the shape optimization problem including an interface--dependent nonlocal model. Once established, we briefly recall basic concepts from the shape optimization regime in Section \ref{sec:basic_concepts_in_shape_opt}. Then Section \ref{sec:nonlocal_shape_opt} is devoted to the task of computing the shape derivative of the nonlocal bilinear form and the reduced objective functional. Finally we present numerical illustrations in Section \ref{sec:shape_numerics} which corroborate theoretical findings.


\section{Problem formulation} \label{sec:problemformulation}
\begin{figure}[h!] 
	\centering
	\def\svgwidth{0.5\textwidth}
	{\Large 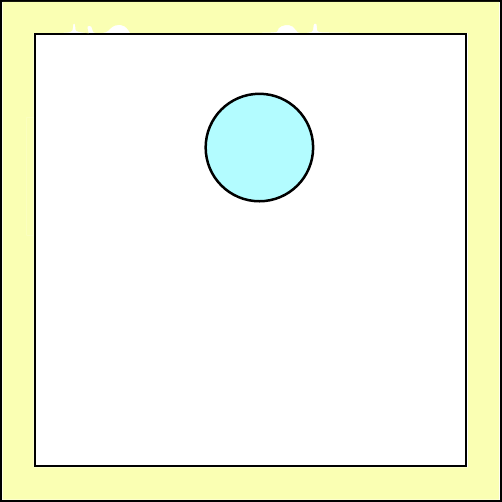}
	\caption[Interface configuration in the nonlocal setting.]{
		An example configuration. 
	}
	\label{fig:shape_problem}
\end{figure}

The system model to be considered is the \emph{homogeneous steady-state nonlocal Dirichlet problem with volume constraints}, given by 
\begin{equation}\label{def:problem}
\left\{
\begin{aligned}
-\mcL_\shape u &= \bfff_\shape \quad \text{on } \omgs \\
u  &= 0 \quad \text{~on } \omgc,
\end{aligned}
\right.
\end{equation} posed on a bounded domain $\Omega \subset \R^d$; see, e.g, \cite{Du_book,RossiBuch, wellposedness, nonsym_ker_1,tian_diss} and the references therein. Here, we assume that this domain is partitioned into a simply connected interior subdomain $\Omega_1 \subset \Omega$ with boundary $\shape \defas \partial \Omega_1$ and a domain $\Omega_2 \defas \Omega \backslash \overline{\Omega}_1$. Thus we have $\Omega = \Omega(\shape) = \Omega_1 \dot{\cup} \shape \dot{\cup} \Omega_2$, where $\dot{\cup}$ denotes the disjoint union. In the following, the boundary $\shape$ of the interior domain $\Omega_1$ is called the \emph{interface} and is assumed to be an element of an appropriate shape space; see also Section \ref{sec:basic_concepts_in_shape_opt} for a related discussion. The governing operator $\mcL_\shape$ is an \emph{interface--dependent, nonlocal convection-diffusion operator} of the form
\begin{align}
-\Lu_\shape u(\xb) 	&:= \int_{\R^d} \left(u(\xb)\g_\shape(\xb,\yb) - u(\yb)\g_\shape(\yb,\xb)\right) d\yb, \label{def:nonloc_op}
\end{align}
which is determined by a nonnegative, interface--dependent \emph{(interaction) kernel} $\g_\shape \colon \R^d \times \R^d \to \R$. The second equation in \eqref{def:problem} is called \textit{Dirichlet volume constraint}. It specifies the values of $u$ on the \emph{interaction domain} 
$$\Omega_I := \lk \yb \in \R^d\backslash \Omega\colon~\exists \xb\in \Omega: \g_\shape(\xb,\yb) \neq 0 \rk ,$$ 
which consists of all points in the complement of $\Omega$ that interact with points in $\Omega$.
For ease of exposition, we set $u=0$ on $\Omega_I$, but generally we can use the constrained $u = g$ on $\Omega_I$, if $g$ satisfies some appropriate assumptions.\\ 
Furthermore, we assume that the kernel depends on the interface in the following way
\begin{equation} \label{def:shape_special_mixed_kernel}
\g_\shape(\xb, \yb) = \sum_{i,j = 1,2} \gij(\xb,\yb) \chi_{\Omega_i \times \Omega_j}(\xb,\yb)
    + \sum_{i = 1,2} \giI(\xb,\yb)\chi_{\Omega_i \times \Omega_I}(\xb,\yb),
\end{equation}
where $\chi_{\Omega_i \times \Omega_j}$ denotes the indicator of the set $\Omega_i \times \Omega_j$.
For instance, in \cite{SelesonGunzburgerParksInterface2013} the authors refer to $\gij$ and $\giI$ as inter-- and intra--material coefficients. Notice that we do not need kernels $\gamma_{Ii}$, since $u=0$ on $\Omega_I$. Furthermore, throughout this work we consider truncated interaction kernels, which can be written as
\begin{align}
\gij(\xb,\yb) = \ggg_{ij}(\xb,\yb) \chi_{{\SSS_{i}(\xb)}}(\yb) \text{ and } \giI(\xb,\yb) = \ggg_{iI}(\xb,\yb) \chi_{{\SSS_{i}(\xb)}}(\yb) \text{ for } i,j=1,2
\label{def:intro_kerneldecompositon}
\end{align}
for appropriate positive functions $\ggg_{ij} \colon \R^d \times \R^d \to \R$ and $\ggg_{iI} \colon \R^d \times \R^d \to \R$, which we refer to as \emph{kernel functions}. 
We assume for $i \in \{1,2\}$ that there exist two radii $0<\varepsilon_{i}^1,\varepsilon_{i}^2<\infty$ such that $B_{\varepsilon_{i}^1}(\xb) \subset\SSS_{i}(\xb)\subset B_{\varepsilon_{i}^2}(\xb)$ for all $\xb \in \Omega$, where $B_{\varepsilon_{i}^k}(\xb)$ denotes the Euclidean ball of radius $\varepsilon_{i}^k$. 
In this paper we differentiate between \emph{square integrable kernels} and \emph{singular symmetric kernels}. For square integrable kernels we require $\gij \in L^2(\Omega \times \Omega)$ and $\giI \in L^2(\Omega \times \Omega_I)$, which also implies $\gamma_\shape \in L^2(\OuO \times \OuO)$.
We do not assume that \eqref{def:shape_special_mixed_kernel} is symmetric for this type of kernels.
In  the case of singular symmetric kernels we require the existence of constants $0 <\gamma_* \leq \gamma^* < \infty$, such that
\begin{align*}
    \gamma_* \leq \gamma(\xb,\yb)||\xb - \yb||_2^{d + 2s} \leq \gamma^*
\end{align*}
\text{ for } $\xb \in \Omega \text{ and } \yb \in S_{1}(\xb) \cup S_{2}(\xb)$. Also, since the singular kernel is required to be symmetric, the condition $\gamma(\xb,\yb)=\gamma(\yb,\xb)$, and, respectively, $\phi_{12}(\xb,\yb)=\phi_{21}(\yb,\xb)$, $\phi_{ii}(\xb,\yb)=\phi_{ii}(\yb,\xb)$ has to hold. Because we do not need to define $\gamma_{Ii}$, as described above, there is no further symmetry condition for $\giI$ required.
\begin{ex}
\label{ex:singular_kernel}
One example of such a singular symmetric kernel is given by
\begin{align*}
    \gij(\xb,\yb) \defas \frac{\sigma_{ij}(\xb,\yb)}{||\xb - \yb||^{d + 2s}_2}\ballxy,\quad \giI(\xb,\yb) \defas \frac{\sigma_{iI}(\xb,\yb)}{||\xb - \yb||^{d + 2s}_2}\ballxy,\quad \text{for } i,j=1,2,
\end{align*}
where $0<\varepsilon<\infty$ and the functions $\sigma_{ij},\sigma_{iI}:\R^d \times \R^d \rightarrow \R$ are bounded from below and above by some positive constants $\gamma_*$ and $\gamma^*$. Additionally, the $\sigma_{ii}$ are assumed to be symmetric on $\Omega\times\Omega$ and $\sigma_{12}(\xb,\yb) = \sigma_{21}(\yb,\xb)$ holds for ${\xb,\yb \in \Omega\cup\Omega_I}$.
\end{ex}
For the forcing term $f_\shape$ in $\eqref{def:problem}$ we assume a dependency on the interface in the following way
\begin{equation}\label{def:shape_rhs}
\begin{aligned}
f_\shape (\xb) \defas 
\begin{cases}
f_1(\xb) : &\xb \in \Omega_1\\
f_2(\xb): &\xb \in \Omega_2,
\end{cases}
\end{aligned}
\end{equation}
where we assume that $f_i \in H^1(\Omega)$, $i=1,2$, because we need that $f$ is weakly differentiable in Section \ref{sec:nonlocal_shape_opt}. Figure \ref{fig:shape_problem} illustrates our setting.
%
%
%
%
%

Next, we introduce a variational formulation of problem \eqref{def:problem}. For this purpose we define the corresponding forms
\begin{align} \label{def:shape_nonlocal_forms}
\AAA_{\shape}(u,  \advar) 	\defas \left(-\Lu_\shape u, \advar\right)_{L^2(\Omega)} ~~~~\text{and}~~~~
F_{\shape}(  \advar) 	\defas (f_\shape,  \advar)_{L^2(\Omega)}
\end{align}
for some functions $u,v\colon\OuO \to \R$. By inserting the definitions of the nonlocal operator \eqref{def:nonloc_op} with the kernel given in \eqref{def:intro_kerneldecompositon} and the definition of the forcing term \eqref{def:shape_rhs}, we obtain the \emph{nonlocal bilinear form}
\begin{align}
\AAA_{\shape}(u, \advar)
&=\int_{\Omega} \advar(\xb)\int_{\R^d}   (u(\xb)\g_\shape(\xb,\yb) - u(\yb)\g_\shape(\yb,\xb))  d\yb d\xb   \nonumber\\
\begin{split}
\label{bilinear_form_representation_1}
    &=\sum_{i,j= 1,2} \int_{\Omega_i} \advar(\xb) \int_{\Omega_j} \left(u(\xb)\gij(\xb,\yb) - u(\yb)\gji(\yb,\xb)\right) d\yb d\xb \\
    &~~~~+ \sum_{i=1,2} \int_{\Omega_i} \advar(\xb)u(\xb) \int_{\Omega_I} \giI(\xb,\yb) d\yb d\xb
\end{split}\\
\begin{split}
\label{bilinear_form_representation_2}
    &=\sum_{i,j= 1,2} \frac{1}{2} \int_{\Omega_i} \int_{\Omega_j}\left(\advar(\xb) - \advar(\yb) \right) \left(u(\xb)\gij(\xb,\yb) - u(\yb)\gji(\yb,\xb)\right) d\yb d\xb\\
&~~~~+ \sum_{i=1,2} \int_{\Omega_i} \advar(\xb)u(\xb) \int_{\Omega_I} \giI(\xb,\yb) d\yb d\xb
\end{split}
\end{align}
and the linear functional
\begin{align} \label{def:rhs_shapeprob}
F_\shape( \advar) = \int_{\Omega} f_\shape \advar ~d\xb = \int_{\Omega_1} f_1 \advar ~d\xb + \int_{\Omega_2} f_2  \advar~ d\xb.
\end{align}
In order to derive the second bilinear form \eqref{bilinear_form_representation_2} we used Fubini’s theorem. We make use of both representations \eqref{bilinear_form_representation_1} and \eqref{bilinear_form_representation_2} of the nonlocal bilinear form in the proofs of Section \ref{sec:nonlocal_shape_opt} . For singular symmetric kernels we also use another equivalent representation of the nonlocal bilinear form given by
\begin{align*}
    \AAA_{\shape}(u, \advar) = \frac{1}{2}\iint\limits_{(\OuO)^2} (\advar(\xb) - \advar(\yb))(u(\xb) - u(\yb))\g_\shape(\xb,\yb) ~d\yb d\xb,
\end{align*}
where we again used Fubini's theorem and applied that $u,\advar = 0$ on $\Omega_I$. Next, we use the nonlocal bilinear form to define a seminorm
\begin{align*}
    |||u||| \defas \sqrt{\AAA_\Gamma(u,u)}.
\end{align*}
With this seminorm, we further define the energy spaces
\begin{align*}
    V(\OuO) &\defas \{u \in L^2(\OuO): ||u||_{V(\OuO)} \defas |||u||| + ||u||_{L^2(\OuO)} < \infty\} \text{ and}\\
    V_c(\OuO) &\defas \{u \in V(\OuO): u = 0 \text{ on } \Omega_I \}.
\end{align*}
We now formulate the variational formulation corresponding to problem \eqref{def:problem} as follows
\begin{equation}\label{weakform}
\begin{tabular}{c}
{\em given $f_\shape \in H^1(\Omega)$ find $u \in V_c(\OuO) $ such that } 
\\[1ex]
$\AAA_\shape(u,v) = F_\shape(v)~~\text{\em for all}~~ v\in V_c(\OuO).$\\[1ex]
\end{tabular}
\end{equation}
For square integrable kernels one can show under appropriate assumptions on the kernel the equivalence between $V(\OuO)$ and $L^2(\OuO)$ and, respectively, the equivalence of $V_c(\OuO)$ and $L_c^2(\OuO)$, see related results in \cite{wellposedness, FelsingerKassmannVoigtDirichlet, Vollmann_diss}. Moreover, the well-posedness of problem \eqref{weakform} for symmetric (square integrable) kernels is proven in \cite{wellposedness} and in \cite{FelsingerKassmannVoigtDirichlet} the well-posedness for some nonsymmetric cases is also covered (again under certain conditions on the kernel and the forcing term $f$). For the singular symmetric kernels the well-posedness of problem \eqref{weakform}, the equivalence between $V(\OuO)$ and the fractional Sobolev space $H^s(\OuO)$ and between $V_c(\OuO)$ and $H_c^s(\OuO)$ is shown in \cite{wellposedness}.\\
%
Finally, let us suppose we are given certain measurements $\bar{u}\colon \Omega \to \R$ on the domain $\Omega$, which we assume to follow the nonlocal model \eqref{weakform} with the interface--dependent kernel $\g_\shape$ and the forcing term $f_\shape$ defined in \eqref{def:shape_special_mixed_kernel} and \eqref{def:shape_rhs}, respectively. In oder to formulate the shape derivative in Chapter \ref{sec:nonlocal_shape_opt} we need $\bar{u} \in H^1(\Omega)$. Then, given the data $\bar{u}$ we aim at identifying the interface $\shape$ for which the corresponding nonlocal solution $u(\shape)$ is the ``best approximation'' to these measurements. Mathematically spoken, we formulate an optimal control problem with a tracking-type objective functional where the interface $\shape$, modeled as a shape, represents the control variable. We now assume
$\Omega := (0,1)^2$ and introduce the following nonlocally constrained shape optimization problem
\begin{framed}
	\begin{equation}\label{def:shape_opt_problem}
	\begin{aligned}
	\min\limits_{\shape}\; &J(u,\shape) \\
	\mbox{s.t. } &\AAA_\shape(u, \advar ) = F_{\shape}( \advar )~~~\text{for all}~ \advar \in V_c(\OuO).
	\end{aligned}
	\end{equation}
\end{framed}
\noindent The objective functional is given by
\begin{align*}
J(u,\shape)			&\defas j(u,\shape)+j_{reg}(\shape)\defas \frac{1}{2} \int_{\Omega} (u - \bar{u})^2 ~ d \xb + \nu \int_{\shape}1 ~ds .
\end{align*}
The first term $j(u,\shape)$ is a standard $L^2$ tracking-type functional ``projecting'' the data on the set of reachable solutions, whereas the second term $j_{reg}(\shape)$ is known as the perimeter regularization, which is commonly used in the related literature to overcome possible ill-posedness of optimization problems \cite{perimeter_regular}.
\section{Basic concepts in shape optimization} \label{sec:basic_concepts_in_shape_opt}
For solving the constrained shape optimization problem \eqref{def:shape_opt_problem} we want to use the same shape optimization algorithms as they are developed in \cite{shape_volker1, shape_volker2, shape_volker4} for problem classes that are comparable in structure. Thus, in this section we briefly introduce the basic concepts and ideas of the therein applied shape formalism. For a rigorous introduction to shape spaces, shape derivatives and shape calculus in general, we refer to the monographs \cite{Delfour-Zolesio-2011,SokoZol,Welker_diss}.

\subsection{Notations and definitions} \label{subsec:shape_notations}

Based on our perception of the interface, we now refer to the image of a simple closed and smooth curve as a \emph{shape}, i.e., the spaces of interest are subsets of 
\begin{align}
\A \defas \lk \shape \defas \varphi(S^1)\colon \varphi \in C^\infty(S^1, \Omega) ~\text{injective};~\varphi' \neq 0 \rk. \label{def:A_shapes}
\end{align} 
By the Jordan curve theorem \cite{jordan_curvetheorem} such a shape $\shape \in \A$ divides the plane into two (simply) connected components with common boundary $\shape$. One of them is the bounded interior, which in our situation can then be identified with $\Omega_1$.\\ 
Functionals $J \colon \A \to \R$ which assign a real number to a shape are called \emph{shape functionals}. Since this paper deals with minimizing such shape functionals, i.e., with so-called shape optimization problems, we need to introduce the notion of an appropriate \emph{shape derivative}. To this end we consider a family of mappings $F_t\colon \overline{\Omega} \to \R^d$ with $F_0 = id$, where $t\in [0,T]$ and $T > 0$, which transform a shape $\shape$ into a family of \emph{perturbed shapes} $\lk \shapet \rk_{t \in [0,T]}$, where $ \shapet \defas F_t(\shape) $ with $\shape^0 = \shape$. Here the family of mappings $\lk F_t \rk_{t\in[0,T]}$ is described by the \emph{perturbation of identity}, which for a smooth vector field $\Vb \in C_0^k(\Omega, \R^d)$, $k \in \mathbb{N}$, is defined by 
$$F_t(\xb ) \defas \xb + t \Vb(\xb).$$ 
We note that for sufficiently small $t\in [0,T]$ the function $F_t$ is injective, and thus $\shapet \in \A$. Then the \emph{Eulerian or directional derivative} of a shape functional $J$ at a shape $\shape$ in direction of a vector field $\Vb \in C_0^k(\Omega, \R^d)$, $k \in \mathbb{N}$, is defined by 
\begin{align}\label{def:shape_derivative}
D_\shape J(\shape)[\Vb]  \defas \left. \frac{d}{dt}\right|_{t = 0^+} J(F_t(\shape))= \lim_{t \searrow 0} \frac{\left( J(F_t(\shape)) - J(\shape)\right)}{t}.
\end{align}
If $D_\shape J(\shape)[\Vb] $ exists for all $\Vb \in C_0^k(\Omega, \R^d)$, $\Vb \mapsto DJ(\shape)[\Vb] $ is continuous and in the dual space $\left(C_0^k(\Omega, \R^d)\right)^*$, then $DJ(\shape)[\Vb] $ is called the \emph{shape derivative} of $J$ \cite[Definition 4.6]{Welker_diss}.\\
At this point, let us also define the \emph{material derivative} of a family of functions ${ \{\advar^t\colon \Omega \to \R: t \in [ 0, T ] \} }$ in direction $\Vb$ by
\begin{align*}
D_m \advar(\xb) \defas \left.\frac{d}{dt}\right|_{t = 0^+} \advar^t(F_t(\xb)).
\end{align*}
For functions $\advar$, which do not explicitly depend on the shape, i.e., $\advar^t = \advar~\text{for all } t\in [0,T]$, we derive 
\begin{align*}
    D_m \advar = \nabla \advar^\top \Vb.
\end{align*}
For more details on shape optimization we refer to the literature, e.g., \cite{SokoNov}.
\subsection{Optimization approach: Averaged adjoint method} \label{subsec:saddlepointsystem}
Let us assume, that for each admissible shape $\shape$, there exists a unique solution $u(\shape)$ of the constraint equation, i.e., $u(\shape)$ satisfies $\AAA_\shape(u(\shape), \advar) = F_{\shape}( \advar)$ for all $\advar \in V_c(\OuO)$. Then we can consider the \emph{reduced problem}
\begin{align}
\label{eq_reduced_functional}
\min\limits_{\shape}\; & J^{red}(\shape) := J(u(\shape), \shape).
\end{align}
In order to employ derivative based minimization algorithms, we need to derive the shape derivative of the reduced objective functional $J^{red}$. By formally applying the chain rule, we obtain
\begin{align*}
D_\shape J^{red}(\shape)[\Vb]  = D_uJ(u(\shape), \shape) D_\shape u(\shape)[\Vb]  + D_\shape J(u(\shape), \shape)[\Vb] ,
\end{align*}
where $D_uJ$ and $D_\shape J$ denote the partial derivatives of the objective $J$ with respect to the state variable $u$ and the control $\shape$, respectively. In applications we typically do not have an explicit formula for the control-to-state mapping $u(\shape)$, so that we cannot analytically quantify the sensitivity of the unique solution $u(\shape)$ with respect to the interface $\shape$. Thus, a formula for the shape derivative $D_\shape u(\shape)[\Vb] $ is unattainable. One possible approach to access this derivative is the \emph{averaged adjoint method (AAM)} of \cite{AAM,Sturm_diss}, which is a Lagrangian method, where the so-called \emph{Langrangian functional} is defined as
\begin{align*}
    L(u,\shape,\advar) \defas J(u,\shape) + \AAA_\shape(u,\advar) - F_\shape(\advar).
\end{align*}
The basic idea behind Lagrangian methods is the aspect, that we can express the reduced functional as
\begin{align*}
    J^{red}(\shape)= L(u(\shape),\shape,\advar), \quad \forall v \in V_c(\OuO).
\end{align*}
Now let $\shape$ be fixed and denote by $\shapet \defas F_t(\shape)$ and $\Omega_i^t \defas F_t(\Omega_i)$ the deformed interior boundary and respectively the deformed domains. Furthermore let $\Omega^t\defas F_t(\Omega)$ be the domain which corresponds to the interior boundary $\shapet$. Then we consider the reduced objective functional regarding $\shapet$, i.e.,
\begin{align}
    \label{deformed_problem}
    J^{red}(\shapet)= L(u(\shapet),\shapet,\advar), \quad \forall v \in V_c(\OtuO),
\end{align}
where $u(\shapet) \in V_c(\OtuO)$. If we now try to differentiate $L$ with respect to $t$ in order to derive the shape derivative, we would have to compute the derivative for $u(\shapet) \circ F_t$ and $\advar \circ F_t$, where $u(\shapet),\advar \in V_c(\OtuO)$ may not be differentiable. Additionally the norm $||\cdot||_{V(\OtuO)}$, and therefore the space $V_c(\OtuO)$, is also dependent on $t$. Instead, since $F_t$ is a homeomorphism, we can use that for $u,\advar \in V_c(\OtuO)$, there exist functions $\tilde{u},\tilde{\advar} \in V_c(\OuO)$, such that
\begin{align*}
    u = \tilde{u} \circ F_t^{-1} ~\text{and}~ \advar = \tilde{\advar} \circ F_t^{-1}.
\end{align*}
Moreover we define
\begin{align}
    J:[0,T]\times &V_c(\OuO)\rightarrow \R,\nonumber\\
    &J(t,u) \defas J(u \circ F_t^{-1},\shapet), \nonumber\\
    A:[0,T]\times &V_c(\OuO)\times V_c(\OuO) \rightarrow \R, \nonumber\\
    &A(t,u,\advar)\defas A_{\shapet}(u\circ F_t^{-1}, \advar \circ F_t^{-1}), \nonumber\\
    F:[0,T]\times &V_c(\OuO)\rightarrow \R,\nonumber\\
    &F(t,\advar)\defas F_{\shapet}(\advar \circ F_t^{-1}), \nonumber\\
    \begin{split}
    G:[0,T] \times &V_c(\OuO) \times V_c(\OuO) \rightarrow \R,\\
    &G(t,u,\advar) \defas L(u \circ F_t^{-1},\shapet,\advar \circ F_t^{-1})
    = J(t, u)  + \AAA(t,u,\advar) - F(t, \advar). \label{AAM_function}
    \end{split}
\end{align}
Then we can reformulate (\ref{deformed_problem}) as
\begin{align*}
    J^{red}(\shapet)= G(t,u^t,\advar),\quad \forall v \in V_c(\OuO),
\end{align*}
where $u^t \in V_c(\OuO)$ is the unique solution of the nonlocal equation corresponding to $\shape^t$
\begin{align*}
    \AAA(t,u,\advar) - F(t,\advar) = 0,\quad \forall \advar \in V_c(\OuO).
\end{align*}
Furthermore $\AAA(t,u,\advar) - F(t,\advar)$ is obviously linear in $\advar$ for all $(t,u) \in [0,T] \times V_c(\OuO)$, which is one prerequisite of the AAM.
Then, in order to use the AAM to compute the shape derivative, the following additional assumptions have to be met.
\begin{itemize}
\item \textbf{Assumption (H0)}:
For every $(t,\advar) \in [0,T] \times V_c(\OuO)$
\begin{enumerate}
\item $[0,1] \ni s \rightarrow G(t,su^t+(1-s)u^0,\advar)$ is absolutely continuous and
\item $[0,1] \ni s \rightarrow d_uG(t,su^t + (1-s)u^0,\advar)[\tilde{u}] \in L^1(0,1)$ for all $\tilde{u} \in V_c(\OuO)$.
\end{enumerate}
\item For every $t \in [0,T]$ there exists a unique solution $\advar^t \in L^2(\Omega)$, such that $\advar^t$ solves the average adjoint equation
\begin{equation}
\label{AAE}
\int_0^1 d_u G(t,su^t + (1-s)u^0,\advar^t)[\tilde{u}]ds = 0 \quad \text{for all } \tilde{u} \in V_c(\OuO).
\end{equation}
\item \textbf{Assumption (H1)}:\\
Assume that the following equation holds
\begin{equation*}
\lim_{t \searrow 0} \frac{G(t,u^0,\advar^t)- G(0,u^0,\advar^t)}{t}= \partial_t G(0,u^0,\advar^0).
\end{equation*}
\end{itemize}
Then the next theorem yields a practical formula for deriving the shape derivative:
\begin{theo}[{\cite[Theorem 3.1]{AAM}}]
Let the assumptions (H0) and (H1) be satisfied and suppose there exists a unique solution $\advar^t$ to the average adjoint equation (\ref{AAE}). Then for $\advar \in V_c(\OuO)$ we obtain
\begin{align}
\label{shape_derivative_reduced_functional}
    D_\shape J^{red}(\Gamma)[\Vb]=\left.\frac{d}{dt}\right|_{t=0^+}J^{red}(\shape^t)=\left.\frac{d}{dt}\right|_{t=0^+}G(t,u^t,\advar) = \partial_t G(0,u^0,\advar^0).
\end{align}
\end{theo}
\begin{proof}
See proof of \cite[Theorem 3.1]{AAM}.
\end{proof}
For $t=0$ the average adjoint equation \eqref{AAE} can be written as
\begin{align}
\label{eq:adjoint}
    \AAA(t,\tilde{u},\advar^0) &= - \int_\Omega (u^0 - \bar{u})\tilde{u} ~d\xb \quad \forall \tilde{u} \in V_c(\OuO).
\end{align}
Here we also call \eqref{eq:adjoint} adjoint equation and the solution $\advar^0$ is referred to as the adjoint solution. Moreover the nonlocal problem \eqref{weakform} is also called state equation and the solution $u^0$ is named state solution.  

\subsection{Optimization algorithm} \label{subsec:optimization_algorithm_shape}
Let us assume for a moment that we have an explicit formula for the shape derivative of the reduced objective functional. We now briefly recall the techniques developed in \cite{shape_volker1} and describe how to exploit this derivative for implementing gradient based optimization methods or even Quasi-Newton methods, such as L-BFGS, to solve the constrained shape optimization problem \eqref{def:shape_opt_problem}.  

In order to identify gradients we need to require the notion of an inner product, or more generally a Riemannian metric. Unfortunately, shape spaces typically do not admit the structure of a linear space. However, in particular situations it is possible to define appropriate quotient spaces, which can be equipped with a Riemannian structure. For instance consider the set $\A$ introduced in \eqref{def:A_shapes}. Since we are only interested in the image of the defining embedding, a re-parametrization thereof does not lead to a different shape. Consequently, two curves that are equal modulo (diffeomorphic) re-parametrizations define the same shape. This conception naturally leads to the quotient space $\Emb(S^1, \R^d) / \Diff(S^1, S^1) $, which can be considered an infinite-dimensional Riemannian manifold \cite{michor}. This example already intimates the difficulty of translating abstract shape derivatives into discrete optimization methods; see, e.g., the thesis \cite{Welker_diffeological} on this topic. A detailed discussion of these issues is not the intention of this work and we now outline Algorithm \ref{alg:shape_optimization}.

The basic idea can be intuitively explained in the following way. Starting with an initial guess $\shape_0$, we aim to iterate in a steepest-descent fashion over interfaces $\shape_k$ until we reach a ``stationary point'' of the reduced objective functional $J^{red}$. The interface $\shape_k$ is encoded in the finite element mesh and transformations thereof are realized by adding vector fields $\Ub\colon \Omega \to \R^d$ (which can be interpreted as tangent vectors at a fixed interface) to the finite element nodes which we denote by $\Omega_k$. 
Thus, the essential part is to update the finite element mesh after each iteration by adding an appropriate transformation vector field. For this purpose, we use the solution $\Ub(\shape)\colon\Omega(\shape) \to \R^d$ of the so-called \emph{deformation equation}
\begin{equation}\label{eq:deformation_equation}
a_{\shape}(\Ub(\shape),\Vb)= D_\shape J^{red}(\shape)[\Vb]  ~~~\text{for all}~\Vb\in H^1_0(\Omega(\shape),\mathbb{R}^2).
\end{equation}
The right-hand side of this equation is given by the shape derivative of the reduced objective functional \eqref{shape_derivative_reduced_functional} and the left-hand side denotes an inner product on the vector field space $H^1_0(\Omega,\mathbb{R}^2)$. In the view of the manifold interpretation, we can consider $a_{\shape}$ as inner product on the tangent space at $\shape$, so that $\Ub(\shape)$ is interpretable as the gradient of the shape functional $J^{red}$ at $\shape$. The solution $\Ub(\shape)\colon \Omega \to \mathbb{R}^2$ of (\ref{eq:deformation_equation}) is then added to the coordinates $\Omega_k$ of the finite element nodes.\\
A common choice for $a_\shape$ is the bilinear form associated to the linear elasticity equation given by
\begin{equation*}
a_\shape(\Ub,\Vb)= \int\limits_{\Omega(\shape)}\sigma(\Ub):\epsilon(\Vb)\,dx,
\end{equation*}
for $\Ub ,\Vb \in H^1_0(\Omega,\mathbb{R}^2)$, where 
\begin{align} \label{def:strain_tensor}
\sigma(\Ub):= \lambda \text{tr}(\epsilon(\Ub)) \Id + 2 \mu \epsilon(\Ub)
\end{align}  
and 
$$\epsilon(\Ub):= \frac{1}{2}(\nabla \Ub + \nabla \Ub^T)$$ 
are the strain and stress tensors, respectively. 
Deformation vector fields $\Vb$ which do not change the interface do not have an impact on the reduced objective functional, so that 
$$D_\shape J^{red}(\shape)[\Vb] =0 ~~~\text{for all}~  \Vb \text{ with }
\text{supp}(\Vb)\cap\Gamma=\emptyset.$$ 
Therefore, the right-hand side $D_\shape J^{red}(\shape)[\Vb] $ is only assembled for test vector fields whose support intersects with the interface $\shape$ and set to zero for all other basis vector fields. This prevents wrong mesh deformations resulting from discretization errors as outlined and illustrated in \cite{shape_volker2}. 
Furthermore, $\lambda$ and $\mu$ in \eqref{def:strain_tensor} denote the Lam\'{e} parameters which do not need to have a physical meaning here. It is more important to understand their effect on the mesh deformation. They enable us to control the stiffness of the material and thus can be interpreted as some sort of step size. In \cite{shape_volker4}, it is observed that locally varying Lam\'{e} parameters have a stabilizing effect on the mesh. A good strategy is to choose $\lambda=0$ and $\mu$ as solution of the following Laplace equation
\begin{equation} \label{eq:locally_varying_lame}
\begin{aligned}
-\Delta \mu &= 0 ~~\quad\quad\text{in } \Omega\\
\mu &= \mu_{\text{max}} \quad \text{on }\shape\\
\mu &= \mu_{\text{min}} ~\quad \text{on }\partial \Omega.
\end{aligned}
\end{equation}
Therefore $\mu_\text{min},\mu_\text{max}\in \mathbb{R}$ influence the step size of the optimization algorithm. A small step is achieved by the choice of a large $\mu_\text{max}$. Note that $a_\shape$ then depends on the interface $\shape$ through the parameter $\mu = \mu(\shape) \colon \Omega(\shape) \to \R$.\\
\begin{algorithm}[H] \label{alg:shape_optimization}
	Initialize: $\g_\shape, f_\shape, \shape^0, \bar{u}$, $k=1$,  $\texttt{maxiter} \in \Nb$  \\
	\While{$k \leq$ \texttt{maxiter} (alternatively $\|DJ^{red}(\shape_k)\|$ > tol)}{
		Interpolate $\bar{u}$ onto the current finite element mesh $\Omega_k$ \label{shape_algo:interpolation}\\
		Assemble $\AAA_\shape$ and solve state \eqref{weakform} and adjoint equation \eqref{eq:adjoint} \label{shape_algo:assembly_state_adjoint}\\
		~~~~~~~ $\rightarrow$ $u(\shape_k)$, $\advar(\shape_k)$\\
		Compute the mesh deformation\\
		~~~~~~~Assemble shape derivative \\
		~~~~~~~~~~~~~~~~~~$D_\shape J^{red}(\shape_k)[\Vb]  = D_\shape L(u(\shape_k), \shape_k, \advar(\shape_k))[\Vb] $ \eqref{shape_derivative_reduced_functional} \label{shape_algo:shape_der}\\
		~~~~~~~Set $D_\shape J^{red}(\shape_k)[\Vb]  = 0$ for all $\Vb$ with $\text{supp}(\Vb) \cap \shape_k = \emptyset$\\
		~~~~~~~Compute locally varying Lam\'e parameter by solving \eqref{eq:locally_varying_lame}\\
		~~~~~~~Assemble linear elasticity $a_{\shape_k}$ and solve the deformation equation \eqref{eq:deformation_equation}\\
		~~~~~~~$\rightarrow$ $\Ub_k$ \\
		\eIf{curvature condition is satisfied}
		{$\tilde{\Ub}_k = \text{ L-BFGS-Update} \label{shape_algo:lbfgs}$}
        {$\tilde{\Ub}_k = - U_k$}		
		Backtracking line search (with parameters $\alpha=1, \tau ,c\in (0,1)$)\\	
		\SetInd{30pt}{6pt}
		~~~~~~~\testwhile{$J^{red}((id + \alpha \tilde{\Ub}_k)(\shape_k)) \geq c J^{red}(\shape_k)$ }{ \label{shape_algo:line search}
			$\alpha = \tau \alpha$} 
		~~~~~~~$\rightarrow$ $\alpha_k$ \\
		Update mesh\\
		~~~~~~~$\Omega_{k+1} = (id + \alpha_k \tilde{\Ub}_k)(\Omega_k)$\\
		$k = k+1$
	}
	\caption{Shape optimization algorithm}
\end{algorithm}

How to perform the limited memory L-BFGS update in Line \ref{shape_algo:lbfgs} of Algorithm \ref{alg:shape_optimization} within the shape formalism is investigated in \cite[Section 4]{shape_volker3}. Here, we only mention that the therein examined vector transport is approximated with the identity operator, so that we finally treat the gradients $\Ub_k\colon\Omega_k \to \R^d$ as vectors in $\R^{d|\Omega_k|}$ and implement the standard L-BFGS update \cite[Section 5]{shape_volker4}.

\section{Shape derivative of the reduced objective functional}
\label{sec:nonlocal_shape_opt}
In Section \ref{sec:basic_concepts_in_shape_opt} we have depicted the optimization methodology, that we follow in this work to numerically solve the constrained shape optimization problem \eqref{def:shape_opt_problem}. In order to proof the requirements of the AAM, we need some additional assumptions.\\ \\
\textbf{Assumption (P0):}
\begin{itemize}
    \item For every $t \in [0,T]$, there exist unique solutions $u^t,\advar^t \in V_c(\OuO)$, such that
\begin{align}
    A(t,u^t,\advar) &= F(t,\advar) \text{ for all } \advar \in V_c(\OuO) \text{ and } \nonumber \\
    A(t,u,\advar^t) &= \left( -(\frac{1}{2}(u^t + u^0) - \bar{u})\xt,u \right)_{L^2(\OuO)} \text{ for all } u \in V_c(\OuO), \label{eq:nonlocal_AAE}
\end{align}
where $A(t,u,\advar)$ and $F(t,\advar)$ are defined as in \eqref{AAM_function} and $\xt(\xb) \defas \det DF_t(\xb)$.
\item Additionally assume that there exists a constant $0 < C_0 < \infty$, such that
\begin{align*}
    A(t,u,u) \geq C_0 ||u||^2_{L^2(\Omega)} \text{ for all } t \in [0,T] \text{ and } u \in V_c(\OuO),
\end{align*}
where $A$ is defined as in \eqref{AAM_function}.
\end{itemize}
\textbf{Assumption (P1):}
\begin{itemize}
    \item For singular kernels:\\
    $|\nabla_\xb\gamma(\xb,\yb)\Vb(\xb) + \nabla_\yb\gamma(\xb,\yb)\Vb(\yb)|||\xb - \yb||_2^{d+2s} \in L^{\infty}((\OuO)^2)$.
    \item For square integrable kernels:\\
    Let the kernel functions satisfy the requirements ${\ggg_{ij} \in H^1(\Omega \times \Omega)}$, 
    ${\ggg_{ij}, \nabla \ggg_{ij} \in L^\infty(\Omega \times \Omega)}$, \\
    ${\ggg_{iI} \in H^1(\Omega \times \Omega_I)}$ and ${\ggg_{iI}, \nabla \ggg_{iI} \in L^\infty(\Omega \times \Omega_I)}$.
\end{itemize}
Singular kernels already satisfy assumption (P0) since the first condition is fulfilled by the theory of \cite{wellposedness, Vollmann_diss} and the second requirement is shown in the following Lemma:
\begin{lem}
In the case of a singular kernel, there exists a constant $0 < C_0 < \infty$, so that
\begin{align*}
A(t,u,u) \geq C_0||u||_{L^2(\Omega)},\quad \text{for every } t \in [0,T],\ u \in H^s(\Omega).
\end{align*}
\end{lem}
\begin{proof}
Let $\varepsilon:=\min\{\varepsilon^1_1,\varepsilon^1_2\}$. Applying \cite[Lemma 4.3]{wellposedness} there exists a constant $C_* > 0$ for the kernel $\frac{\gamma_*}{||\xb - \yb||^{2+2s}}\ballxy$, s.t.
\begin{align*}
C_*||u||_{L^2(\Omega)} &\leq \iint\limits_{(\OuO)^2} \frac{1}{2} (u(\xb) - u(\yb))^2\frac{\gamma_*}{||\xb - \yb||_2^{d+2s}}\ballxy ~ d\yb d\xb \\
&= \iint\limits_{(F_t(\Omega) \cup \Omega_I)^2} \frac{1}{2} (u(\xb) - u(\yb))^2\frac{\gamma_*}{||\xb - \yb||_2^{d+2s}}\ballxy ~ d\yb d\xb \\
&\leq \iint\limits_{(F_t(\Omega) \cup \Omega_I)^2}\frac{1}{2} (u(\xb) - u(\yb))^2\gamma(\xb,\yb)\ballxy ~ d\yb d\xb = A_{\shapet}(u, u)
\end{align*}
So we conclude
\begin{align*}
C_* ||u \circ F_t^{-1}||_{L^2(\Omega)}^2 \leq A_{\shapet}(u \circ F_t^{-1},u \circ F_t^{-1}) = A(t,u,u).
\end{align*}
Since $T$ is chosen small enough, $[0,T] \times \bar{\Omega}$ is a compact set and $\xt$ is continuous on $[0,T] \times \bar{\Omega}$, there exists $\xi_* > 0$, s.t. $\xt(\xb) \geq \xi_*$ for every $t \in [0,T]$ and $\xb \in \bar{\Omega}$. Therefore, by using that $F_t(\Omega)=\Omega$, we derive
\begin{align*}
||u \circ F_t^{-1}||_{L^2(\Omega)}^2 &= \int_\Omega (u \circ F_t^{-1})^2 ~d\xb = \int_{F_t(\Omega)} (u \circ F_t^{-1})^2 ~d\xb = \int_\Omega u^2 \xt ~d\xb \geq \xi_* \int_\Omega u^2 ~d\xb \\ 
&= \xi_* ||u||_{L^2(\Omega)}^2.
\end{align*}
\end{proof}
In the following we proof that assumption (P1) also holds for a standard example of a singular symmetric kernel.
\begin{ex}
For $\gamma(\xb,\yb)= \frac{\sigma(\xb,\yb)}{||\xb-\yb||^{d+2s}}\ballxy$ of Example \ref{ex:singular_kernel}, where additionally there exists a constant $\sigma^* \in (0,\infty)$ with ${|\nabla_\xb \sigma|,|\nabla_\yb \sigma| \leq \sigma^*}$, the assumption (P1) holds, since $\OuO$ is a bounded domain and
\begin{align*}
    &|\nabla_\xb \gamma(\xb,\yb)\Vb(\xb) + \nabla_\yb \gamma(\xb,\yb)\Vb(\yb)|||\xb - \yb||_2^{d+2s}\\ 
    &\leq |\sigma(\xb,\yb)\frac{(\xb-\yb)^\top(\Vb(\xb) - \Vb(\yb))}{||\xb-\yb||_2^2}| + |\nabla_\xb \sigma(\xb,\yb)\Vb(\xb) + \nabla_\yb \sigma(\xb,\yb)\Vb(\yb)| \\
    &\leq L \gamma^* + 2\sigma^*\Vb^* < \infty,
\end{align*}
where we used that $\Vb \in C_0^k(\Omega)$ is Lipschitz continuous for some Lipschitz constant $L > 0$ and that there exists a $\Vb^*>0$ with $|\Vb(\xb)| \leq\Vb^* $ for $\xb \in \OuO$.
\end{ex}
We will see in the proof of the following Lemma \ref{lem:AAM_nonlocal_shape_problem}, that in our case the average adjoint equation \eqref{AAE} is equivalent to equation \eqref{eq:nonlocal_AAE}.
Now we can show, that the additional requirements of AAM are satisfied by problem (\ref{def:shape_opt_problem}):
\begin{lem}
\label{lem:AAM_nonlocal_shape_problem}
Let $G$ be defined as in (\ref{AAM_function}) and let the assumptions (P0) and (P1) be fulfilled. Then the assumptions (H0) and (H1) are satisfied and for every $t \in [0,T]$ there exists a solution ${\advar^t \in V_c(\OuO)}$ that solves the average adjoint equation (\ref{AAE}).
\end{lem}
\begin{proof}
Because of the length of the proof, we move it to Appendix \ref{chap:Proof_AAM_nl_shape_problem}.
\end{proof}
The missing piece to implement the respective algorithmic realization presented in \\Subsection \ref{subsec:optimization_algorithm_shape} is the shape derivative of the reduced objective functional, which is used in \\Line \ref{shape_algo:shape_der} of Algorithm \ref{alg:shape_optimization} and given by
\begin{align}
\begin{split}
D_\shape J^{red}(\shape)[\Vb] = \partial_t G(0,u^0,\advar^0)
= \left.\frac{d}{dt}\right|_{t=0^+} J(t, u^0) + \left.\frac{d}{dt}\right|_{t=0^+} \AAA(t, u^0, \advar^0) - \left.\frac{d}{dt}\right|_{t=0^+} F(t, \advar^0).
\label{eq:Jred_shape_derivative}
\end{split}
\end{align}
As a first step, we formulate the shape derivative of the objective functional $J$ and the linear functional $F$, which can also be found in the standard literature.
\begin{theo}[\textbf{\textit{Shape derivative of the reduced objective functional}}]
\label{theo:shape_derivative_reduced_functional}
Let the assumptions (P0) and (P1) be satisfied. Further let $\shape$ be a shape with corresponding state variable $u^0$ and adjoint variable $\advar^0$. Then, for a vector field $\Vb \in \vecfields$ we find
\begin{equation} \label{eq:shapeder_reduced_final_general}
\begin{aligned}
D_\shape J^{red}(\shape)[\Vb]  
= &\int_{\Omega} -(u^0 - \bar{u})\nabla\bar{u}^\top\Vb + (u^0 - \bar{u})^2 \di \Vb ~d\xb + \nu \int_\Gamma \di \Vb - \nb ^\top \nabla \Vb^\top \nb ~ds\\ &- \int_\Omega D_m f_\shape \advar^0 + \di \Vb (f\advar^0)~d\xb + D_\shape\AAA_{\shape}(u^0, \advar^0) [\Vb].
\end{aligned}
\end{equation}	 
\end{theo}
\begin{proof}
In order to proof this theorem, we just have to compute the shape derivative of the objective function $J(u^0,\shape)$ and of the linear functional $F_{\shape}(\advar^0)$. Therefore, let $\xt(\xb) \defas \det DF_t(\xb)$. \\
Then, we have $\xi^0(x)= \det DF_0(\xb)= \det(I)=1$ and $\left. \frac{d}{dt}\right|_{t = 0^+} \xt = \di \Vb$(see e.g. \cite{Schmidt_diss}), such that the shape derivative of the right-hand side $F_\shape$ can be derived as follows
\begin{align*}
D_\shape F_\shape(\advar^0)[\Vb]  &= \left. \frac{d}{dt}\right|_{t = 0^+} F_{\shapet}(\advar^0 \circ F_t^{-1}) = \int_\Omega \left. \frac{d}{dt}\right|_{t = 0^+} (f_\shape \circ F_t) \advar^0 \xt ~d\xb \\ 
&= \int_{\Omega} D_m f_\shape \advar^0 ~d \xb +\int_{\Omega}f_\shape \advar^0 ~\di \Vb~d\xb.
\end{align*}
Moreover, the shape derivative of the objective functional can be written as
\begin{align*}
    D_{\shape}J(u^0, \shape)[\Vb] = D_{\shape}j(u^0, \shape)[\Vb] + D_{\shape}j_{reg}(\shape)[\Vb] = \left. \frac{d}{dt} \right|_{t = 0^+} j(u^0 \circ F_t^{-1}, \shape^t) + \left. \frac{d}{dt} \right|_{t = 0^+} j_{reg}(\shape^t).
\end{align*}
Here the shape derivative of the regularization term is an immediate consequence of \cite[Theorem 4.13]{Welker_diss} and is given by 
\begin{align*}
D_{\shape}j_{reg}(u^0,\shape)[\Vb]  = \nu\int_{{\Gamma}}\di_{\Gamma} \Vb ~ds 
= \nu\int_{{\Gamma}}\di  \Vb  - \nb^\top \nabla \Vb^\top  \nb~ds,
\end{align*}
where $\nb$ denotes the outer normal of $\Omega_1$. Additionally, we obtain for the shape derivative of the tracking-type functional
\begin{align*}
D_\shape j(u^0,\shape)[\Vb] &= \left. \frac{d}{dt}\right|_{t = 0^+} j(u^0 \circ F_t^{-1}, \shapet) = \frac{1}{2} \left. \frac{d}{dt}\right|_{t = 0^+} \int_{F_t(\Omega)} (u^0 \circ F_t^{-1} - \bar{u})^2 ~d\xb\\ 
&= \frac{1}{2} \int_{\Omega} \left. \frac{d}{dt}\right|_{t = 0^+} (u^0 - \bar{u} \circ F_t)^2 \xt ~d\xb = \int_{\Omega} -(u^0 - \bar{u})\nabla\bar{u}^\top\Vb + (u^0 - \bar{u})^2 \di \Vb ~d\xb.
\end{align*}
Putting the above terms into equation (\ref{eq:Jred_shape_derivative}) yields the formula of Theorem \ref{theo:shape_derivative_reduced_functional}.
\end{proof}
The last step to derive the shape derivative of the reduced objective functional \eqref{eq:Jred_shape_derivative} is to compute the shape derivative of the nonlocal bilinear form $\AAA_\shape$, which is shown in the next Lemma. 
\begin{lem}[\textbf{\textit{Shape derivative of the nonlocal bilinear form}}]
\label{lem:shape_derivativative_nonlocal_bilinear_form}
Let the assumptions (P0) and (P1) be satisfied. Further let $\shape$ be a shape with corresponding state variable $u^0$ and adjoint variable $\advar^0$. Then for a vector field $\Vb \in \vecfields$ we find for a square integrable kernel $\g$ that
 	\begin{equation} \label{eq:shape_derivative_nonlocbilin_general}
 	\begin{aligned}
 	\left.\frac{d}{dt}\right|_{t=0^+} \AAA(t, u^0, \advar^0) &= D_\shape\AAA_{\shape}(u^0, \advar^0) [\Vb]\\ 
 	= \sum_{i,j=1,2} \int_{\Omega_i}\int_{\Omega_j} &\left( \advar^0(\xb) - \advar^0(\yb) \right)\left(u^0(\xb)\nabla_x \gij(\xb,\yb) - u^0(\yb)\nabla_y \gji(\yb,\xb) \right)^\top \Vb(\xb) \\
 	&+ (\advar^0(\xb) - \advar^0(\yb))(u^0(\xb)\gij(\xb,\yb) - u^0(\yb)\gji(\yb,\xb))\di \Vb(\xb) ~d\yb d\xb\\
+ \sum_{i=1,2} \int_{\Omega_i}\int_{\Omega_I} &u^0(\xb)\advar^0(\xb) (\nabla_x \giI(\xb,\yb)^\top \Vb(\xb) + \nabla_y \giI(\xb,\yb)^\top \Vb(\yb)) \\ 
&+ u^0(\xb)\advar^0(\xb) \giI(\xb,\yb) (\di \Vb(\xb) + \di \Vb(\yb))~d\yb d\xb.
 	\end{aligned}
 	\end{equation}
 	and for a singular kernel $\g$ that
 	\begin{align*}
 	&D_\shape\AAA_{\shape}(u^0, \advar^0) [\Vb]\\*
 	&= \sum_{i,j=1,2}\frac{1}{2} \int_{\Omega_i}\int_{\Omega_j}(u^0(\xb) -u^0(\yb))(\advar^0(\xb) - \advar^0(\yb))\left( \nabla_{\xb}\gij(\xb,\yb)\Vb(\xb) + \nabla_{\yb}\gij(\xb,\yb)\Vb(\yb) \right) ~d\yb d\xb \\*
    &+ \sum_{i,j=1,2}~ \int_{\Omega_i}\int_{\Omega_j}(u^0(\xb) -u^0(\yb))(\advar^0(\xb) - \advar^0(\yb))\gij(\xb,\yb) \di \Vb(\xb) ~d\yb d\xb\\*
    &+\sum_{i=1,2} \int_{\Omega_i}\int_{\Omega_I}(u^0(\xb) -u^0(\yb))(\advar^0(\xb) - \advar^0(\yb))\left( \nabla_{\xb}\giI(\xb,\yb)\Vb(\xb) + \nabla_{\yb}\giI(\xb,\yb)\Vb(\yb) \right) ~d\yb d\xb \\*
    &+ \sum_{i=1,2} \int_{\Omega_i}\int_{\Omega_I}(u^0(\xb) -u^0(\yb))(\advar^0(\xb) - \advar^0(\yb))\giI(\xb,\yb)\left(\di \Vb(\xb) + \di \Vb(\yb)\right) ~d\yb d\xb.
 	\end{align*}
\end{lem}
\begin{proof}
Define $\xt(\xb)\defas \det DF_t(\xb)$ and $\gijt(\xb,\yb) \defas \gij(F_t(\xb), F_t(\yb))$.\\
\textbf{Case 1: Square integrable kernels}\\
Then, we can write by using representation \eqref{bilinear_form_representation_2} of the nonlocal bilinear form $\AAA$
\begin{align*}
&\AAA(t, u^0, \advar^0) = \AAA_{\shape^t}(u^0 \circ F_t^{-1},\advar^0 \circ F_t^{-1})\\
&=\frac{1}{2} \sum_{i,j=1,2} \int_{\Omega_i}\int_{\Omega_j} \left(\advar^0(\xb) - \advar^0(\yb) \right)\left(u^0(\xb)\gijt(\xb,\yb) - u^0(\yb)\gjit(\yb,\xb) \right)\xt(\xb)\xt(\yb)~d\yb d\xb \\
&+ \sum_{i=1,2} \int_{\Omega_i}\int_{\Omega_I} u^0(\xb)\advar^0(\xb)\giI(F_t(\xb),\yb)\xt(\xb)~d\yb d\xb.
\end{align*}
So we derive the shape derivative of the nonlocal bilinear form
\begin{align*}
\left.\frac{d}{dt}\right|_{t=0^+} \AAA(t, u^0, \advar^0)\\ 
= \frac{1}{2} \sum_{i,j=1,2} \int_{\Omega_i}\int_{\Omega_j} &\left( \advar^0(\xb) - \advar^0(\yb) \right)\left(u^0(\xb)\nabla_x \gij(\xb,\yb) - u^0(\yb)\nabla_y \gji(\yb,\xb) \right)^\top\Vb(\xb) \\ 
&+ \left( \advar^0(\xb) - \advar^0(\yb) \right)\left(u^0(\xb)\nabla_y \gij(\xb,\yb) - u^0(\yb)\nabla_x \gji(\yb,\xb) \right)^\top\Vb(\yb) \\
&+(\advar^0(\xb) - \advar^0(\yb))(u^0(\xb)\gij(\xb,\yb) - u^0(\yb)\gji(\yb,\xb))(\di \Vb(\xb) + \di \Vb(\yb))~d\yb d\xb\\
+ \sum_{i=1,2} \int_{\Omega_i}\int_{\Omega_I} &u^0(\xb)\advar^0(\xb)(\nabla_x \giI(\xb,\yb)^\top \Vb(\xb) + \nabla_y \giI(\xb,\yb)^\top \Vb(\yb))\\ 
&+ u^0(\xb)\advar^0(\xb)\giI (\di \Vb + \di \Vb(\yb))~d\yb d\xb \\
= \sum_{i,j=1,2} \int_{\Omega_i}\int_{\Omega_j} &\left(\advar^0(\xb) - \advar^0(\yb) \right)\left(u^0(\xb)\nabla_x \gij(\xb,\yb) - u^0(\yb)\nabla_y \gji(\yb,\xb) \right)^\top \Vb(\xb) \\
&+ (\advar^0(\xb) - \advar^0(\yb))(u^0(\xb)\gij(\xb,\yb) - u^0(\yb)\gji(\yb,\xb))\di \Vb(\xb)~d\yb d\xb\\
+ \sum_{i=1,2} \int_{\Omega_i}\int_{\Omega_I} &u^0(\xb)\advar^0(\xb)(\nabla_x \giI(\xb,\yb)^\top \Vb(\xb) + \nabla_y \giI(\xb,\yb)^\top \Vb(\yb)) \\
&+ u^0(\xb)\advar^0(\xb)\giI(\xb,\yb) (\di \Vb(\xb) + \di \Vb(\yb))~d\yb d\xb.
\end{align*}
For the second equation, the following computations are used, which can be obtained by applying Fubini's theorem and by swapping $\xb$ and $\yb$ we obtain
\begin{align*}
    &\int_{\Omega_i} \int_{\Omega_j} (\advar^0(\xb)-\advar^0(\yb))(- u^0(\yb)\nabla_x\gji(\yb,\xb)^\top \Vb(\yb)) ~d\yb d\xb \\
    &= \int_{\Omega_j} \int_{\Omega_i} (\advar^0(\xb)-\advar^0(\yb))(u^0(\xb)\nabla_x\gji(\xb,\yb)^\top \Vb(\xb)) ~d\yb d\xb, \\ \\
    &\int_{\Omega_i} \int_{\Omega_j}(\advar^0(\xb)-\advar^0(\yb))u^0(\xb)\nabla_y\gij(\xb,\yb)^\top \Vb(\yb) ~d\yb d\xb \\
    &= -\int_{\Omega_j} \int_{\Omega_i}(\advar^0(\xb)-\advar^0(\yb))u^0(\yb)\nabla_y\gij(\yb,\xb)^\top \Vb(\xb) ~d\yb d\xb \text{ and }\\ \\
    &\int_{\Omega_i} \int_{\Omega_j}(\advar^0(\xb)-\advar^0(\yb))(u^0(\xb)\gij(\xb,\yb)-u^0(\yb)\gji(\yb,\xb))\di \Vb(\yb)~d\yb d\xb \\ 
    &= \int_{\Omega_j} \int_{\Omega_i}(\advar^0(\xb)-\advar^0(\yb))(u^0(\xb)\gji(\xb,\yb) - u^0(\yb)\gij(\yb,\xb))\di \Vb(\xb) ~d\yb d\xb.
\end{align*}
\textbf{Case 2: Singular kernels}\\
Analogously we get for the singular symmetric kernel
\begin{align*}
    &\left. \frac{d}{dt} \right|_{t=0} A(t,u^0,\advar^0) = \left. \frac{d}{dt} \right|_{t=0} \frac{1}{2}\iint\limits_{(\Omega \cup \Omega_I)^2} (u^0(\xb) - u^0(\yb))(v^0(\xb) - v^0(\yb))\gt(\xb,\yb)\xt(\xb)\xt(\yb) ~d\yb d\xb \\
    &= \sum_{i,j=1,2}\frac{1}{2} \int_{\Omega_i}\int_{\Omega_j}(u^0(\xb) -u^0(\yb))(\advar^0(\xb) - \advar^0(\yb))\left( \nabla_{\xb}\gij(\xb,\yb)\Vb(\xb) + \nabla_{\yb}\gij(\xb,\yb)\Vb(\yb) \right) ~d\yb d\xb \\
    &+ \sum_{i,j=1,2}~ \int_{\Omega_i}\int_{\Omega_j}(u^0(\xb) -u^0(\yb))(\advar^0(\xb) - \advar^0(\yb))\gij(\xb,\yb) \di \Vb(\xb) ~d\yb d\xb\\
    &+\sum_{i=1,2} \int_{\Omega_i}\int_{\Omega_I}(u^0(\xb) - u^0(\yb))(\advar^0(\xb) - \advar^0(\yb))\left( \nabla_{\xb}\giI(\xb,\yb)\Vb(\xb) + \nabla_{\yb}\giI(\xb,\yb)\Vb(\yb) \right) ~d\yb d\xb \\
    &+ \sum_{i=1,2} \int_{\Omega_i}\int_{\Omega_I}(u^0(\xb) - u^0(\yb))(\advar^0(\xb) - \advar^0(\yb))\giI(\xb,\yb)\left(\di \Vb(\xb) + \di \Vb(\yb)\right) ~d\yb d\xb.
\end{align*}
\end{proof}

We can now derive the shape derivative of the reduced objective functional.
If we formally set $u^0 = \bar{u}$ and $\nu = 0$, we can conclude from \eqref{eq:adjoint} that $\advar^0 = 0$ and therefore $D_\shape J^{red}(\shape)[\Vb] = 0$ for all $\Vb \in \vecfields$.
So, if there is a shape $\shape$, such that $u(\shape) =u^0 = \bar{u}$, then $\shape$ is a stationary point of the reduced objective functional \eqref{eq_reduced_functional}.
\section{Numerical experiments} \label{sec:shape_numerics}
In this section, we want to put the above derived formula \eqref{eq:shapeder_reduced_final_general} for the shape derivative of the reduced objective functional into numerical practice.  
In the following numerical examples we test one singular symmetric and one nonsymmetric square integrable kernel. Specifically,
\begin{align*}
    \g^{sym}_\shape(\xb,\yb) = \phi^{sym}(\xb,\yb) \chi_{B_{\delta}(\xb)}(\yb),
\end{align*}
where 
$$
\phi^{sym}(\xb,\yb) = 
\begin{cases}
100 d_{\delta}\frac{1}{||\xb - \yb||_2^{2+2s}} & \text{if } (\xb,\yb) \in \Omega_1 \times \Omega_1, \\
1.0 d_{\delta}\frac{1}{||\xb - \yb||_2^{2+2s}} & \text{if } (\xb,\yb) \in \Omega_2 \times \Omega_2, \\
10 d_{\delta}\frac{1}{||\xb - \yb||_2^{2+2s}} & \text{else}, 
\end{cases}
$$
with with scaling constants $d_{\delta}\defas\frac{2-2s}{\pi \delta^{2-2s}}$ and 
\begin{align*}
\g^{nonsym}_\shape(\xb,\yb) = \phi^{nonsym}(\xb,\yb) \chi_{B_{\delta}(\xb)}(\yb),
\end{align*}
where 
\begin{align*}
\phi^{nonsym}(\xb,\yb) &= 
\begin{cases}
5.0 c_\delta & \text{if } \xb \in \Omega_1,\\
3.0 c_\delta & \text{if } \xb \in \Omega_2,
\end{cases}
\end{align*}
with scaling constants $c_\delta \defas \frac{1}{\del^4}$. We  truncate all kernel functions by  $\|\cdot\|_2$-balls of radius $\delta = 0.1$ so that $\OuO \subset [-\del,1+\del]^2$.
As a right-hand side we choose a piecewise constant function 
\begin{equation*}
f_\shape(\xb) = 100\chi_{\Omega_1}(\xb) -10\chi_{(\Omega_2)}(\xb),
\end{equation*}  
i.e., $f_1 = 100$ and $f_2 = -10$. 
We note that the nonsymmetric kernel $\g^{nonsym}$ satisfies the conditions for the class of integrable kernels considered in \cite{Vollmann_diss}, such that the corresponding nonlocal problem is well-posed and also assumption (P1) can easily be verified in this case. The symmetric kernel $\g^{sym}$ is a special case of Example \ref{ex:singular_kernel} and therefore the assumptions (P0) and (P1) are met. The well-posedness of the nonlocal problem regarding the singular kernel is shown in \cite{wellposedness}.
As a perimeter regularization we choose $\nu = 0.002$ and, since we only utilize $\Vb$ with ${\text{supp}(\Vb) \cap \shape_k \neq \emptyset}$, we additionally assume that the nonlocal boundary has no direct influence on the shape derivative of the nonlocal bilinear form $D_\shape \AAA_\shape$, such that for all $\Vb \in \vecfields$ with ${\text{supp}(\Vb) \cap \shape_k \neq \emptyset}$ we have for the square integrable kernel 
\begin{align*}
    D_\shape\AAA_{\shape}(u^0, \advar^0) [\Vb] = \sum_{i,j=1,2}~& \iint\limits_{\Omega_i\times\Omega_j} \left( \advar^0(\xb) - \advar^0(\yb) \right)\left(u^0(\xb)\nabla_x \gij(\xb,\yb) - u^0(\yb)\nabla_y \gji(\yb,\xb) \right)^\top \Vb(\xb)\\ 
    &+ (\advar^0(\xb) - \advar^0(\yb))(u^0(\xb)\gij(\xb,\yb) - u^0(\yb)\gji(\yb,\xb))\di \Vb(\xb) ~d\yb d\xb
\end{align*}
and for the singular symmetric kernel
\begin{align*}
    &D_\shape\AAA_{\shape}(u^0, \advar^0) [\Vb] \\
    &= \sum_{i,j=1,2}~ \iint\limits_{\Omega_i\times\Omega_j}\frac{1}{2}(u^0(\xb) -u^0(\yb))(\advar^0(\xb) - \advar^0(\yb))\left( \nabla_{\xb}\gij(\xb,\yb)\Vb(\xb) + \nabla_{\yb}\gij(\xb,\yb)\Vb(\yb) \right) \\
    &\hspace{22mm} + (u^0(\xb) -u^0(\yb))(\advar^0(\xb) - \advar^0(\yb))\gij(\xb,\yb) \di \Vb(\xb) ~d\yb d\xb.
\end{align*}
In order to solve problem \eqref{def:shape_opt_problem}, we apply a finite element method, where
we employ continuous piecewise linear basis functions on triangular grids for the discretization of the nonlocal constraint equation. In particular we use the free meshing software Gmsh \cite{gmsh} to construct the meshes and the Python package nlfem \cite{nlfem} to assemble the stiffness matrices of the nonlocal state and adjoint equation. Moreover, to compute the load vector of the state and adjoint equation and the shape derivatives $D_\shape J$ and $D_\shape F_\shape$, we employ the open-source finite element software FEniCS \cite{fenics_article,fenics_book}.
For a detailed discussion on the assembly of the nonlocal stiffness matrix we refer to \cite{fe_cookbook,nlfem}. Here we solely emphasize how to implement a subdomain--dependent kernel of type \eqref{def:shape_special_mixed_kernel}. During the mesh generation each triangle is labeled according to its subdomain affiliation. Thus, whenever we integrate over a pair of two triangles, we can read out the labels $(i,j)$ and choose the corresponding kernel $\g_{ij}$.

The data $\bar{u}$ is generated as solution $u(\overline{\shape})$ of the constraint equation associated to a target shape $\overline{\shape}$. Thus the data is represented as a linear combination of basis functions from the finite element basis. For the interpolation task in Line \ref{shape_algo:interpolation} of Algorithm \ref{alg:shape_optimization} we solely need to translate between (non-matching) finite element grids by using the project function of  FEniCS. In all examples below the target shape $\overline{\shape}$ is chosen to be a circle of radius $0.25$ centered at $(0.5,0.5)$.

\begin{figure}[h!] 
	\centering
	~\\~\\
	\textbf{\textsf{Example 1: singular symmetric kernel}}
	\begin{center}
	\begin{small}
	\begin{tabular}{cccc}
		\includegraphics[width = 0.2\textwidth]{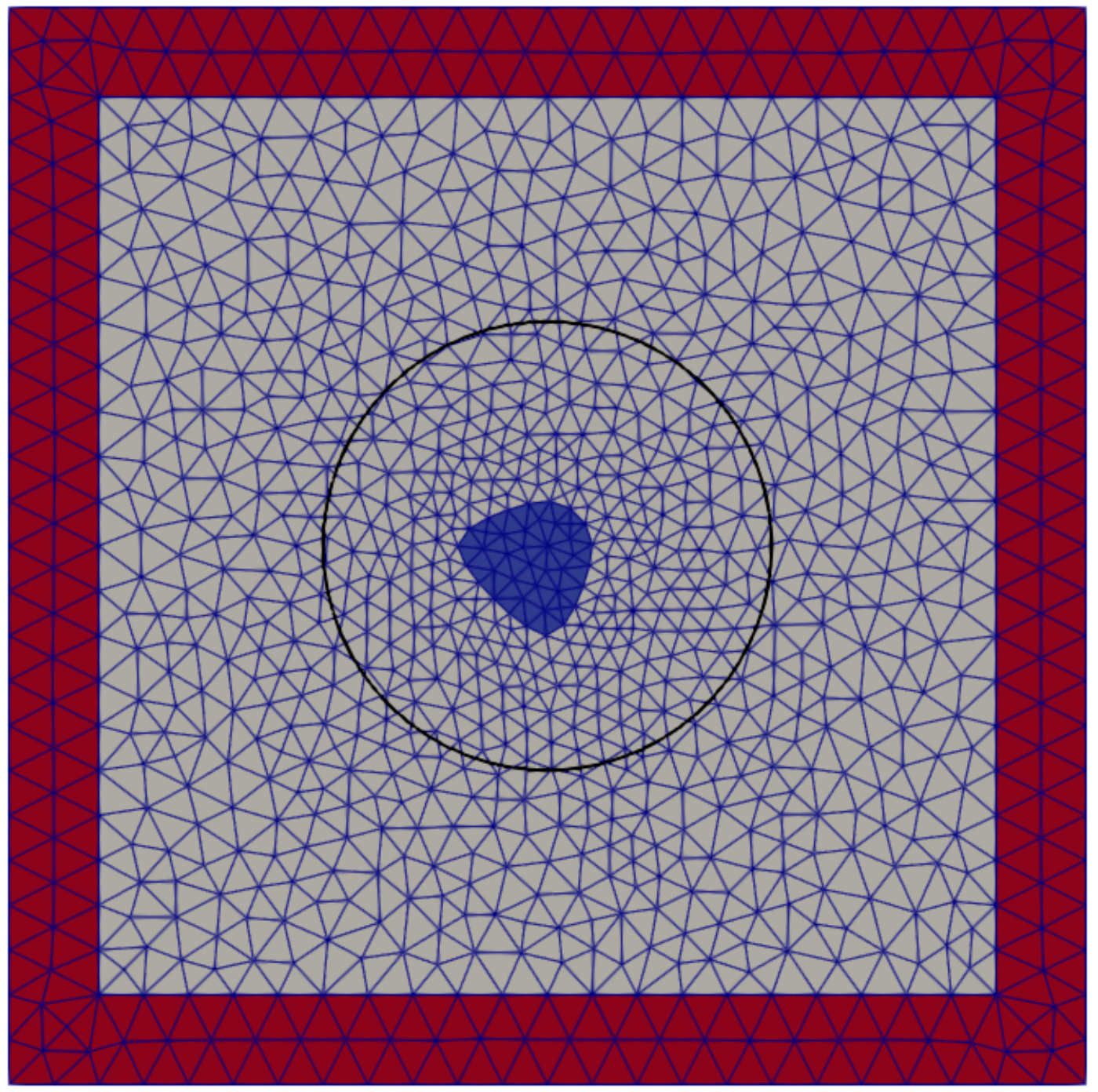}
		&\includegraphics[width = 0.2\textwidth]{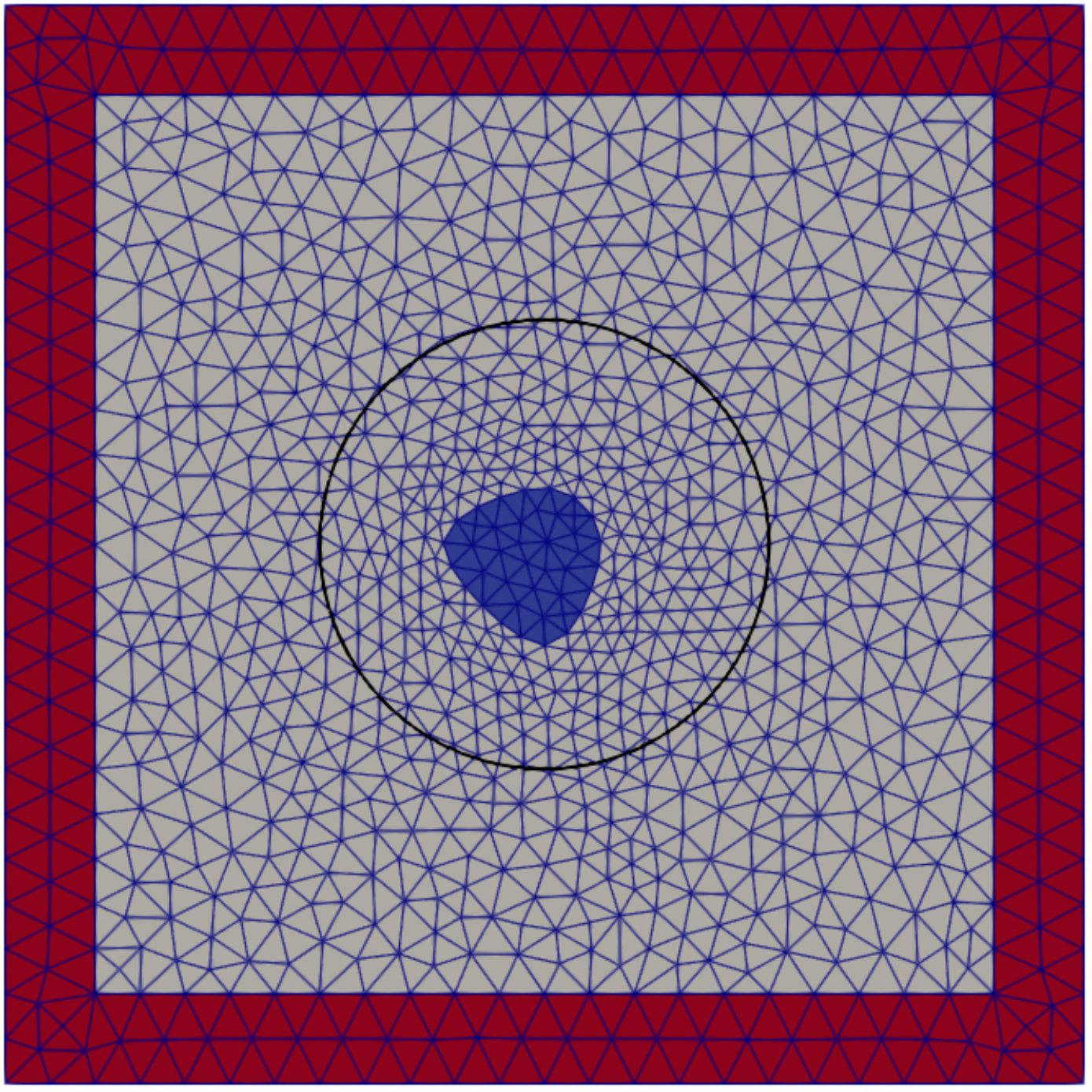}
		&\includegraphics[width = 0.2\textwidth]{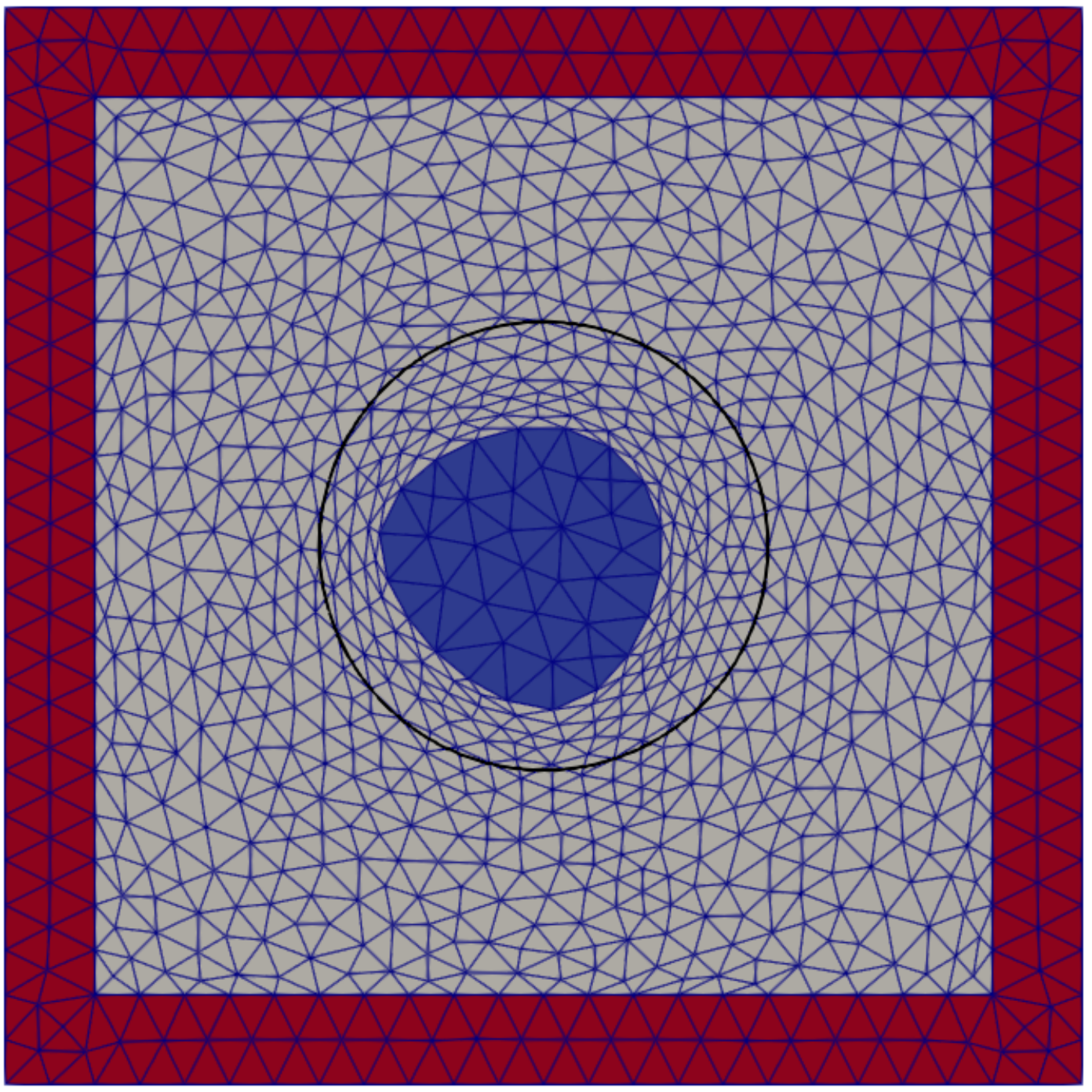}
		&\includegraphics[width = 0.2\textwidth]{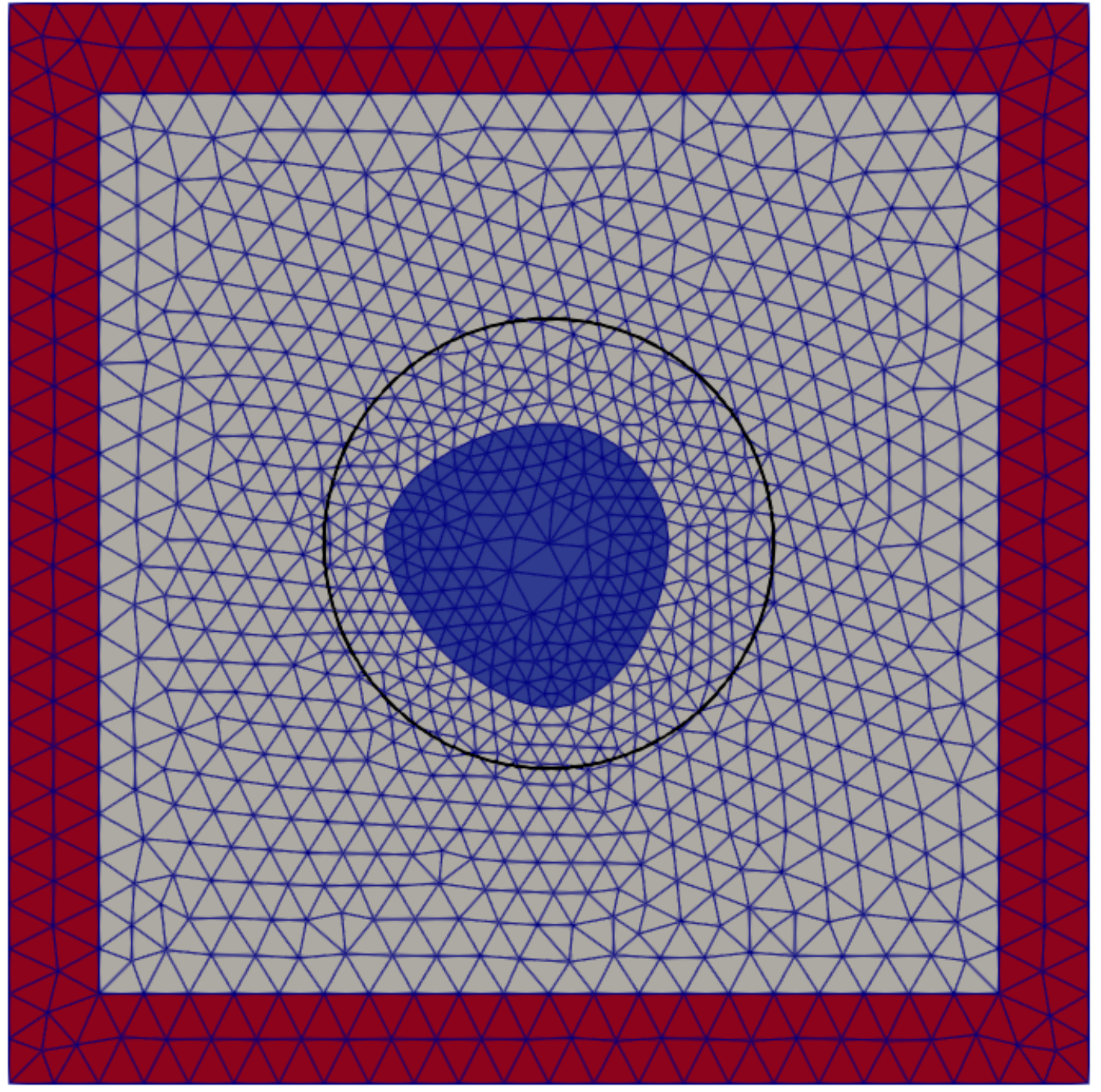}\\
		Start setup & Iteration 1 & Iteration 5 & Iteration 5 remeshed\\
		\includegraphics[width = 0.2\textwidth]{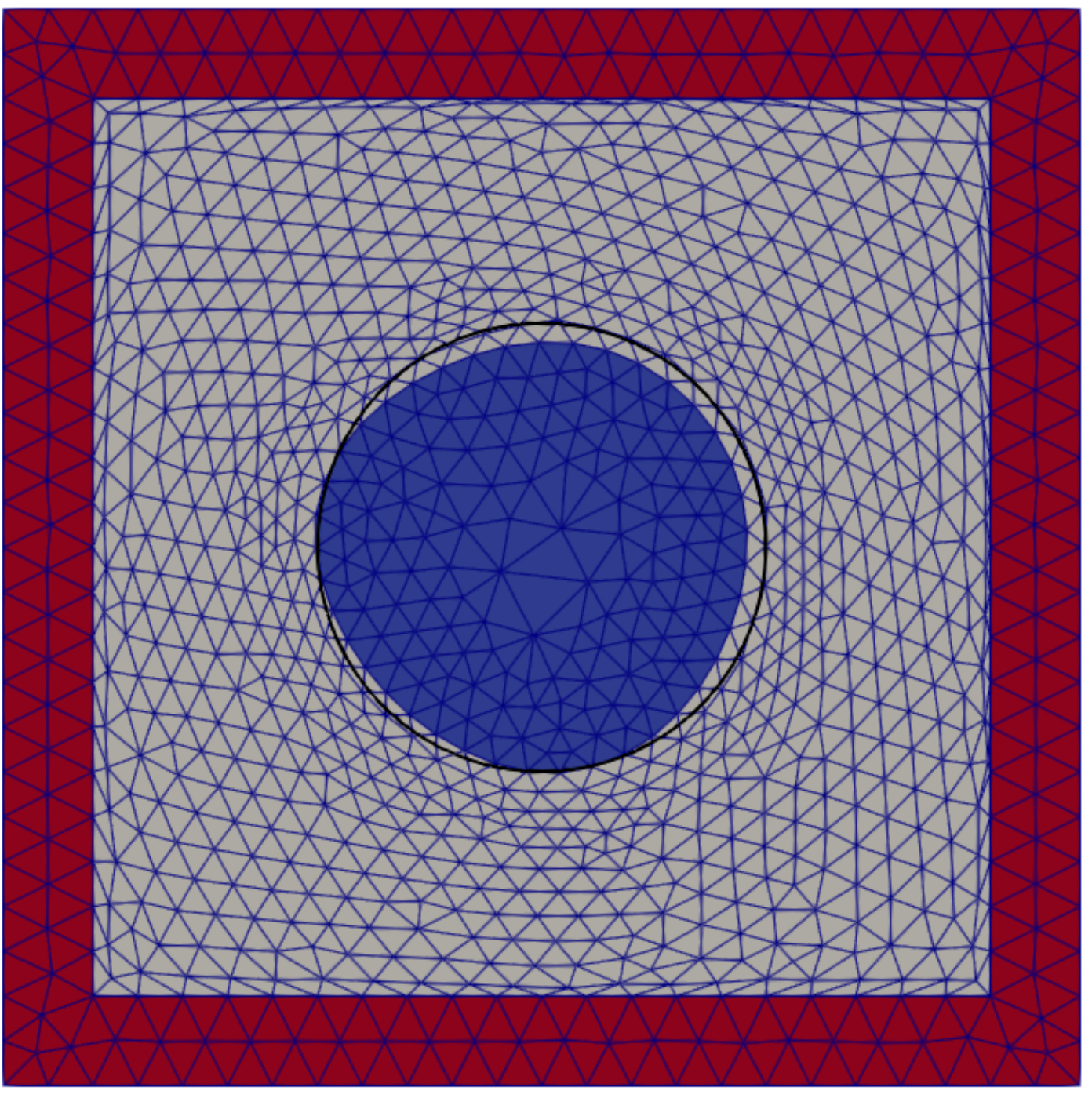}
		&\includegraphics[width = 0.2\textwidth]{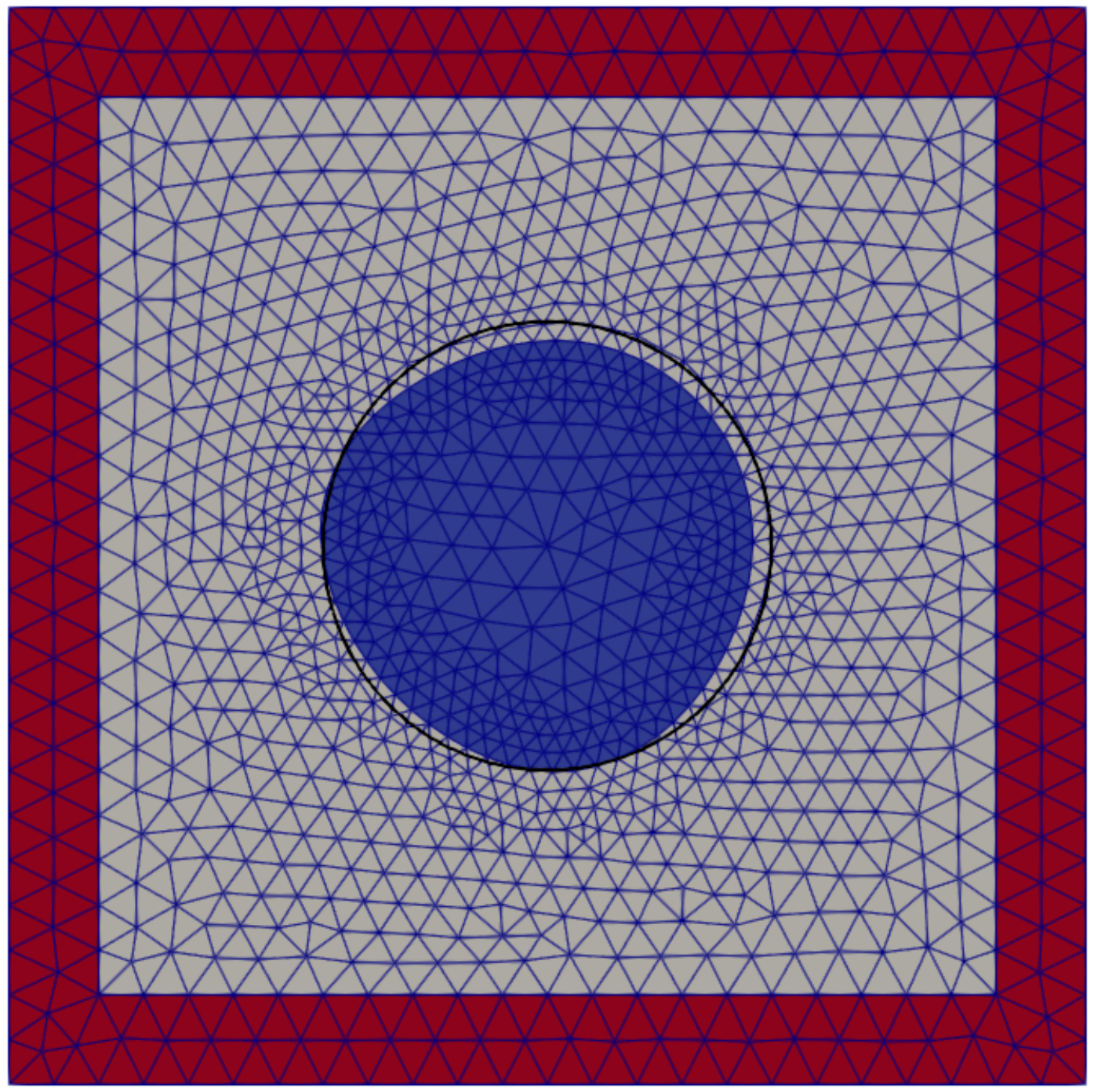}
		&\includegraphics[width = 0.2\textwidth]{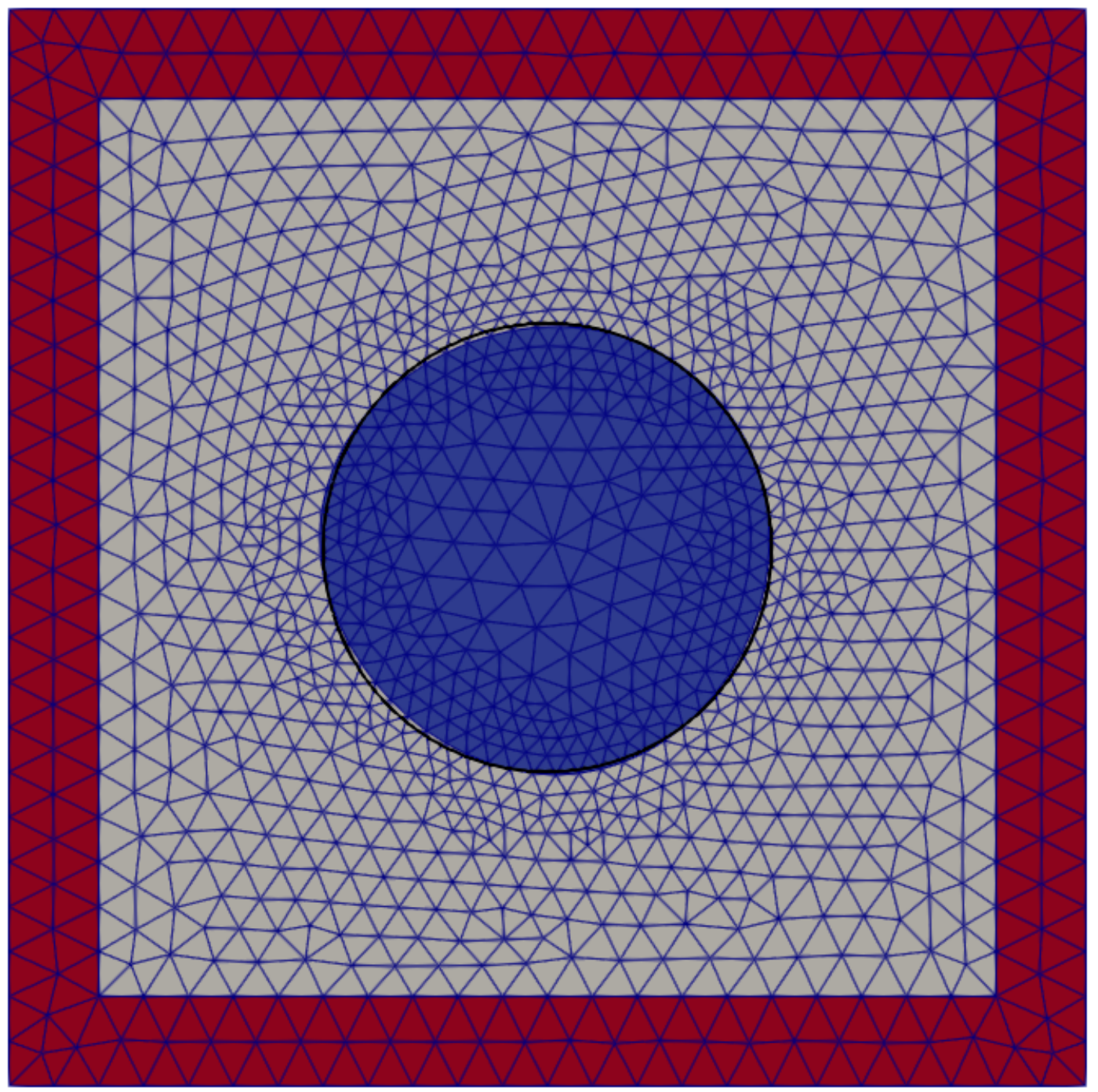}
		&\includegraphics[width = 0.2\textwidth]{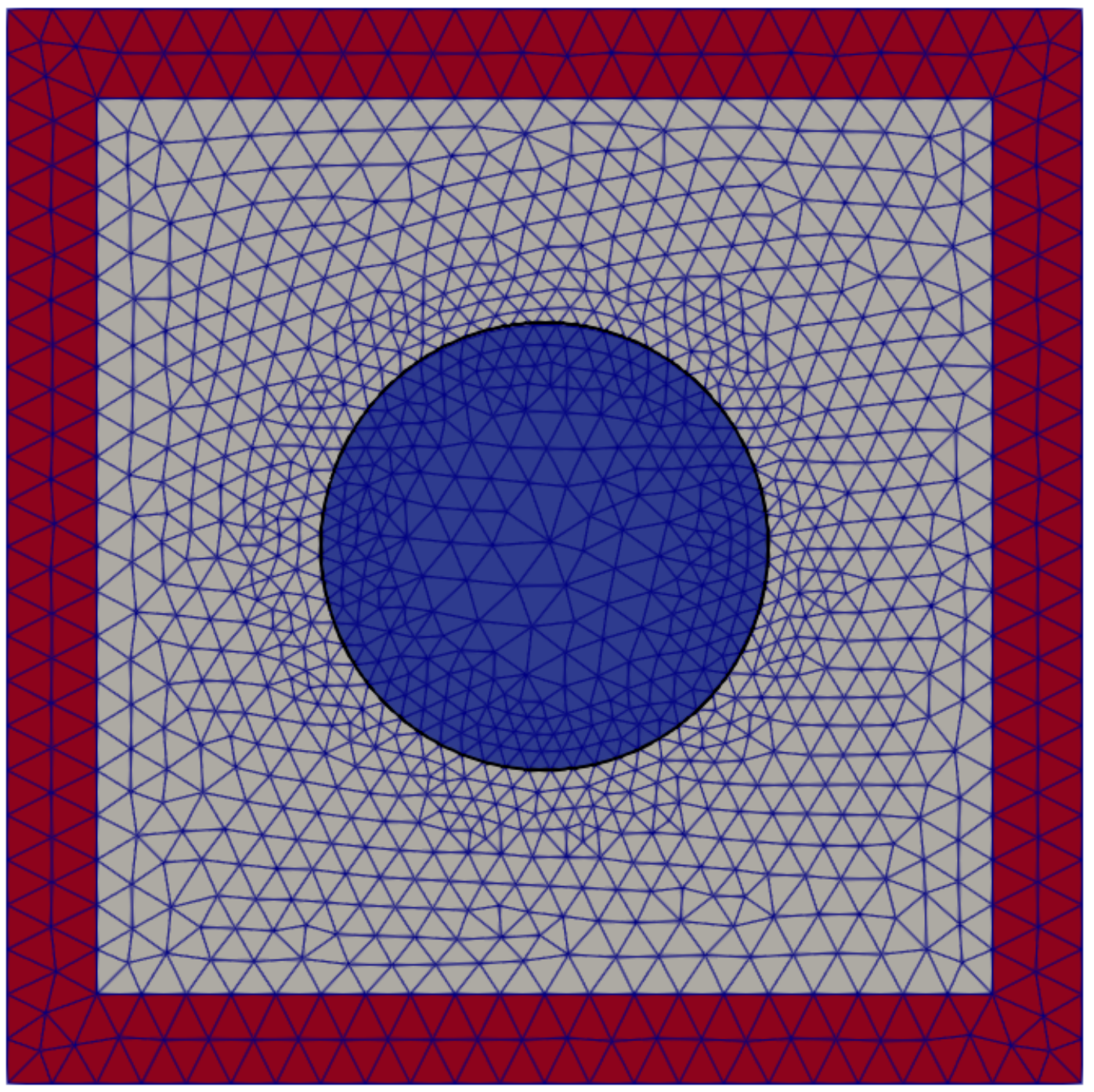}\\
		Iteration 10 & Iteration 10 remeshed & Iteration 15 & Iteration 20
	\end{tabular}
	\end{small}
	\end{center}
	\caption{}
	\label{fig:shape_ex_2_sym}
	~\\
\end{figure}

\begin{figure}[h!] 
	\centering
	~\\~\\
	\textbf{\textsf{Example 1: nonsymmetric integrable kernel}}
	\begin{center}
	\begin{small}
	\begin{tabular}{cccc}
		\includegraphics[width = 0.2\textwidth]{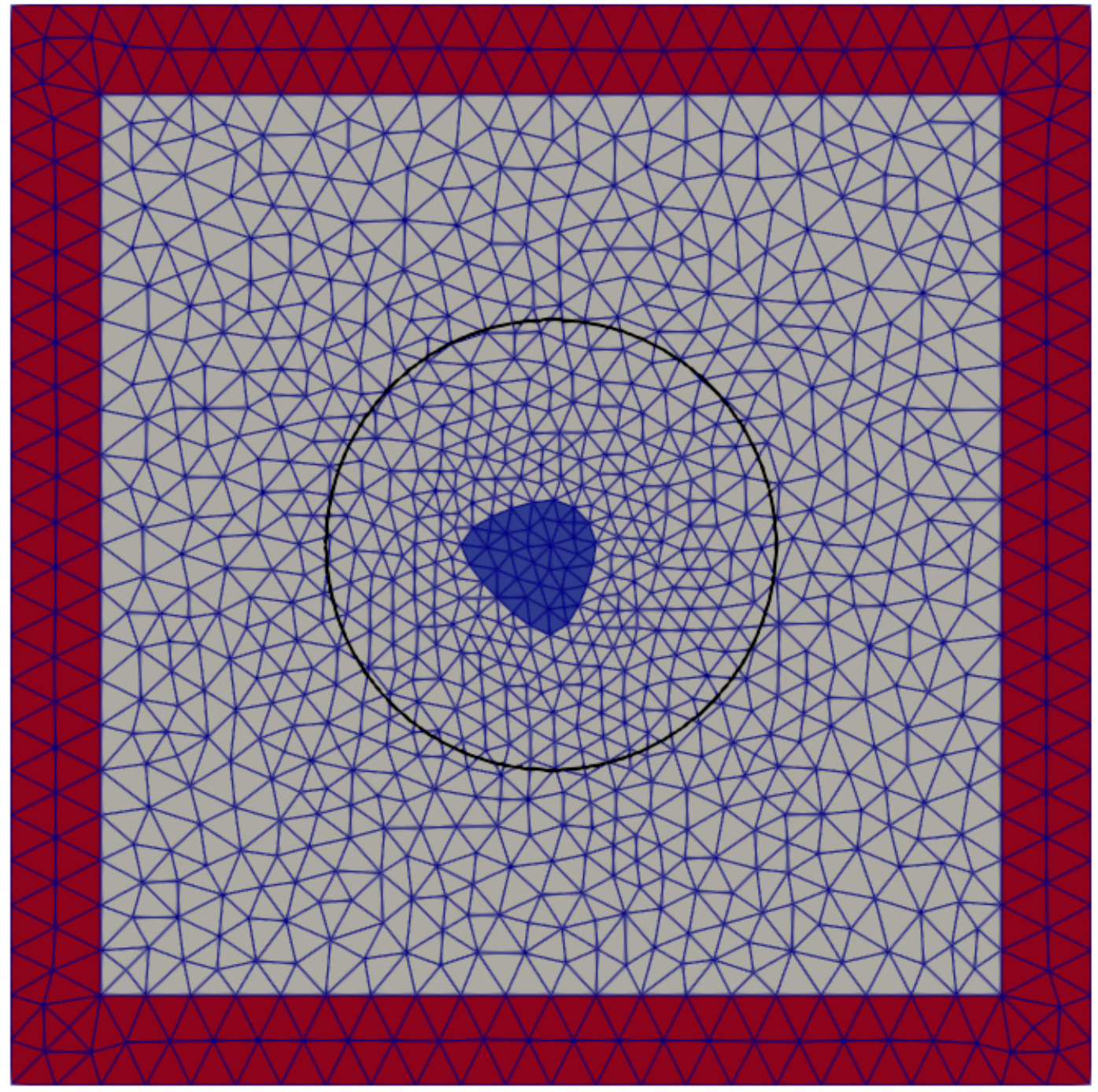}
		&\includegraphics[width = 0.2\textwidth]{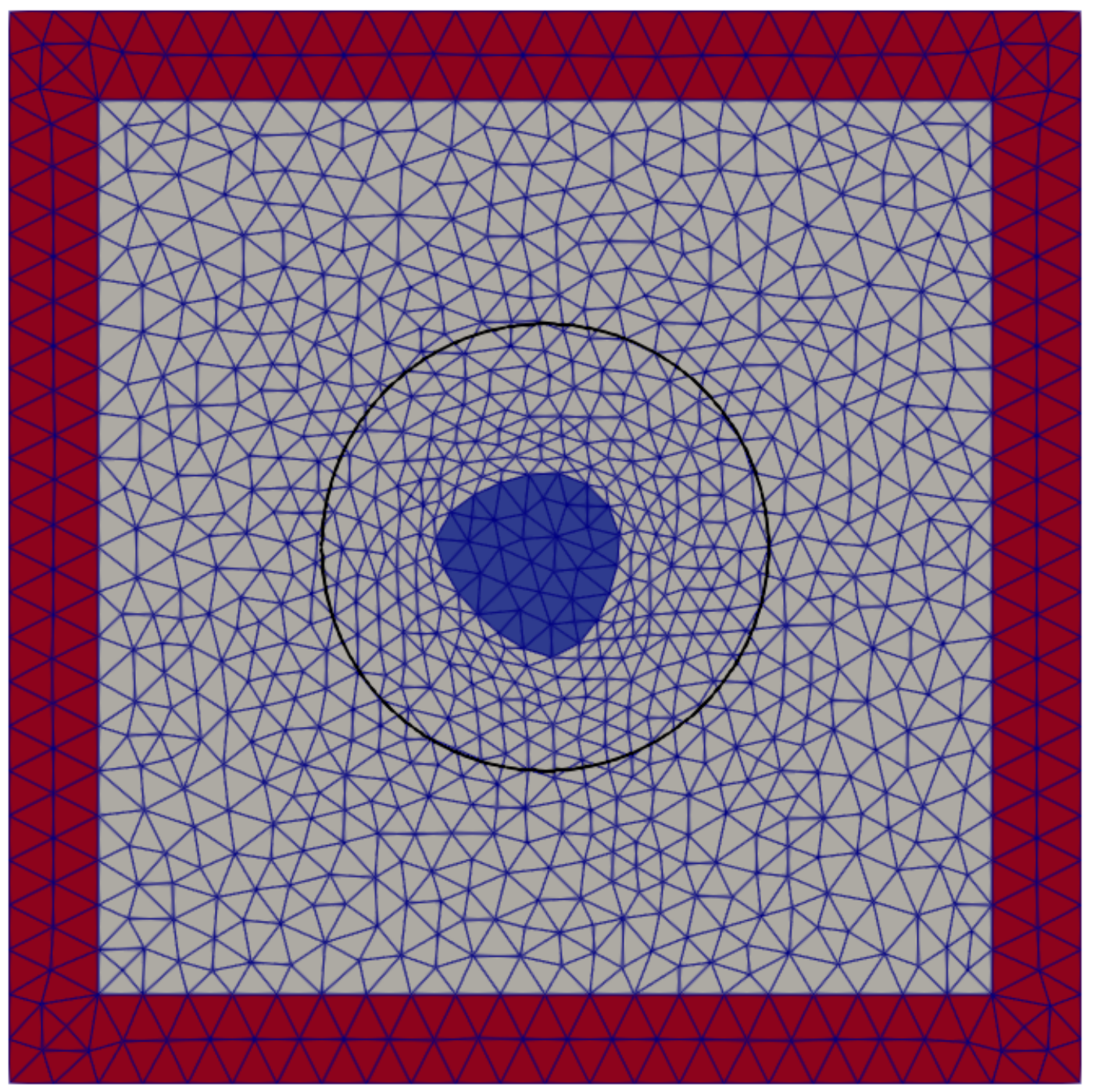}
		&\includegraphics[width = 0.2\textwidth]{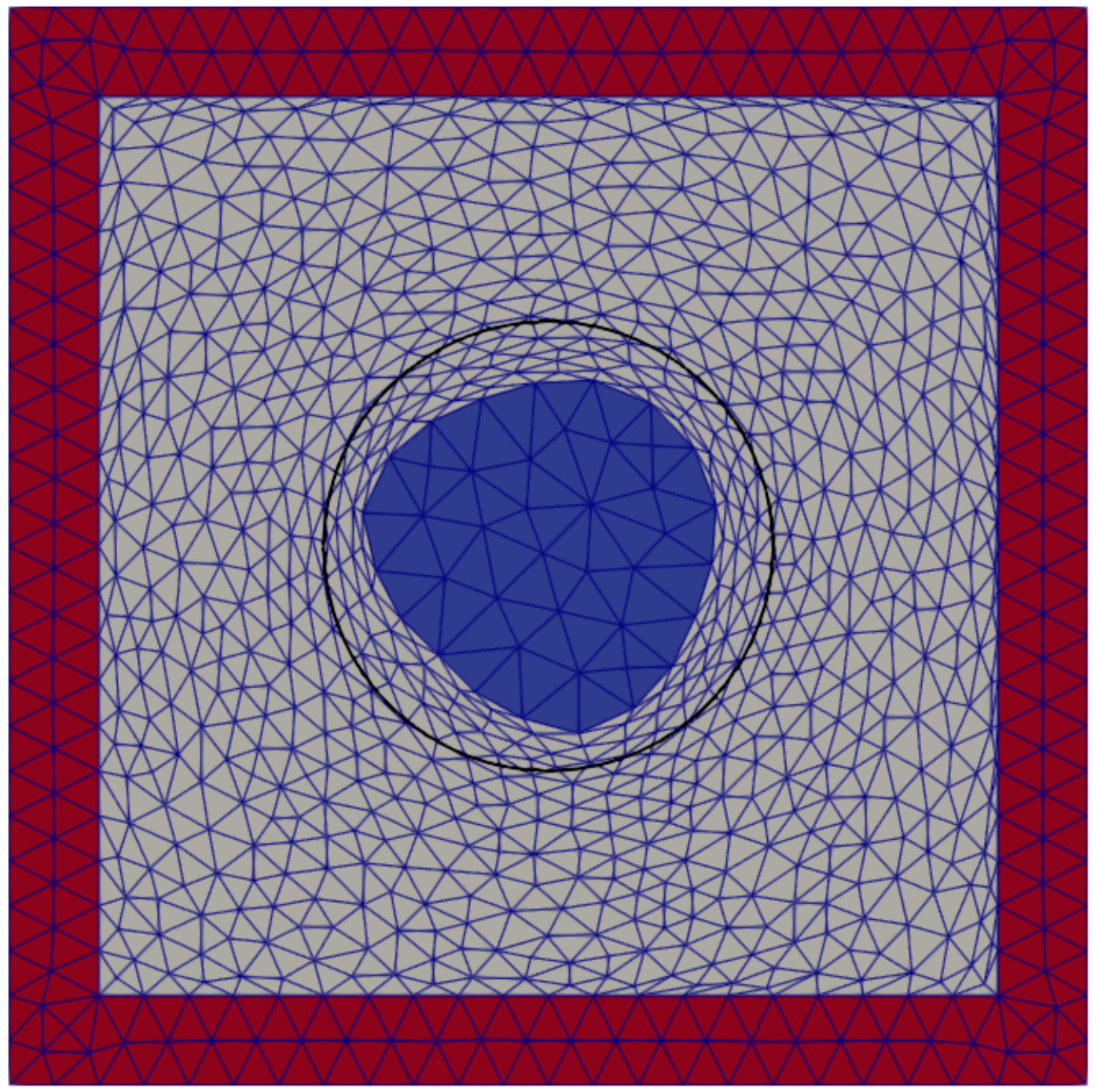}
		&\includegraphics[width = 0.2\textwidth]{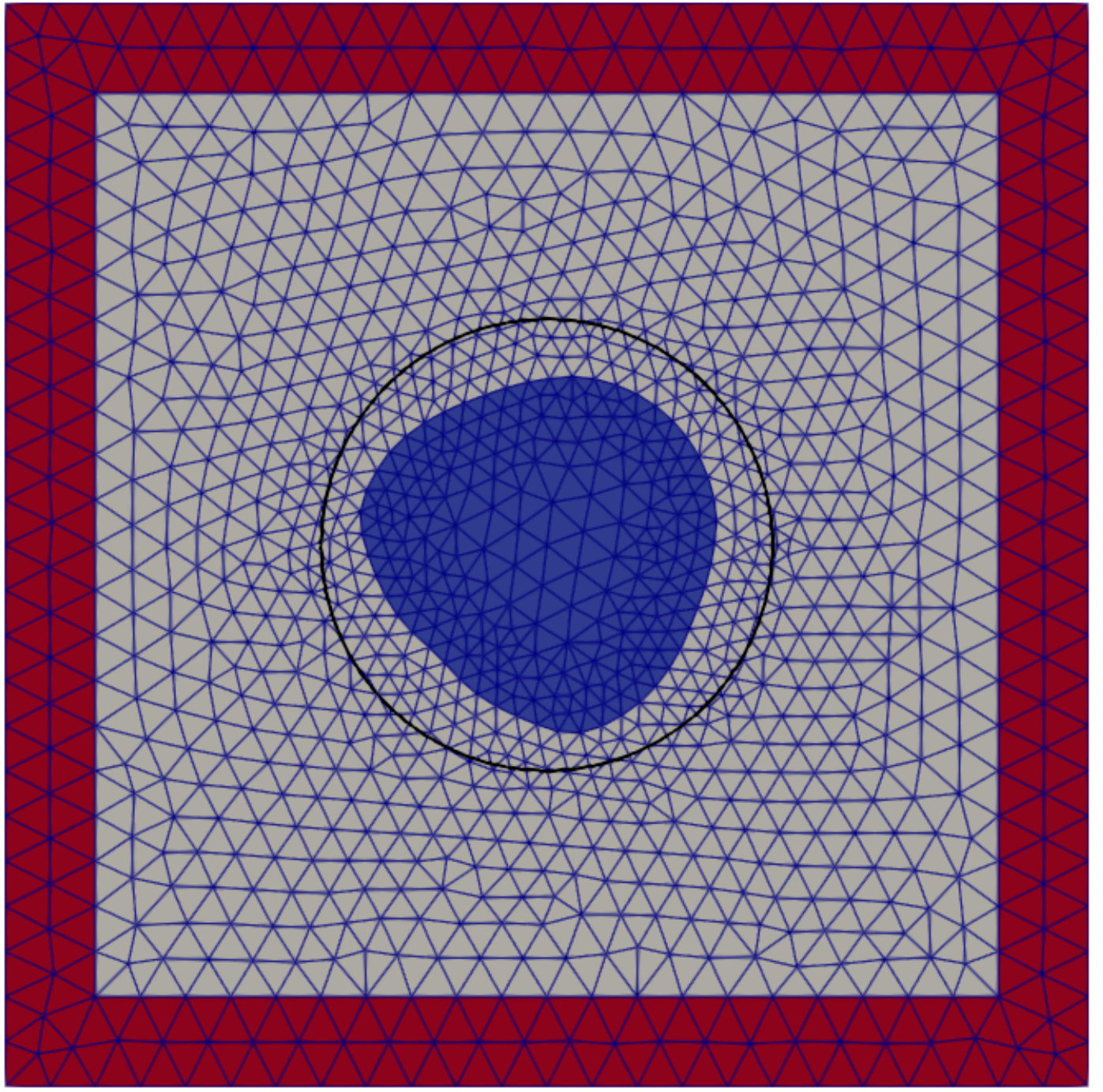}\\
		Start setup & Iteration 1 & Iteration 5 & Iteration 5 remeshed\\
		\includegraphics[width = 0.2\textwidth]{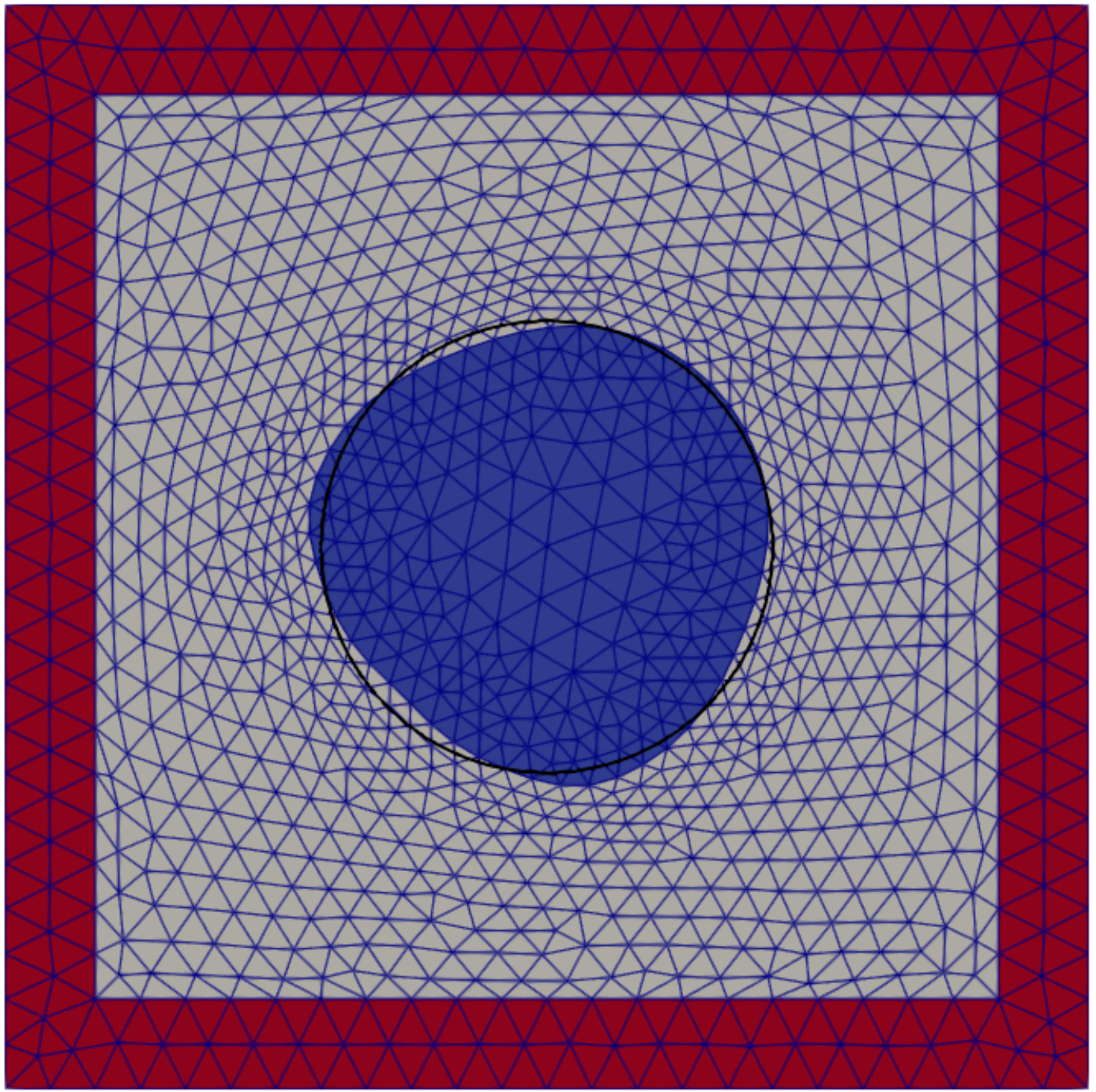}
		&\includegraphics[width = 0.2\textwidth]{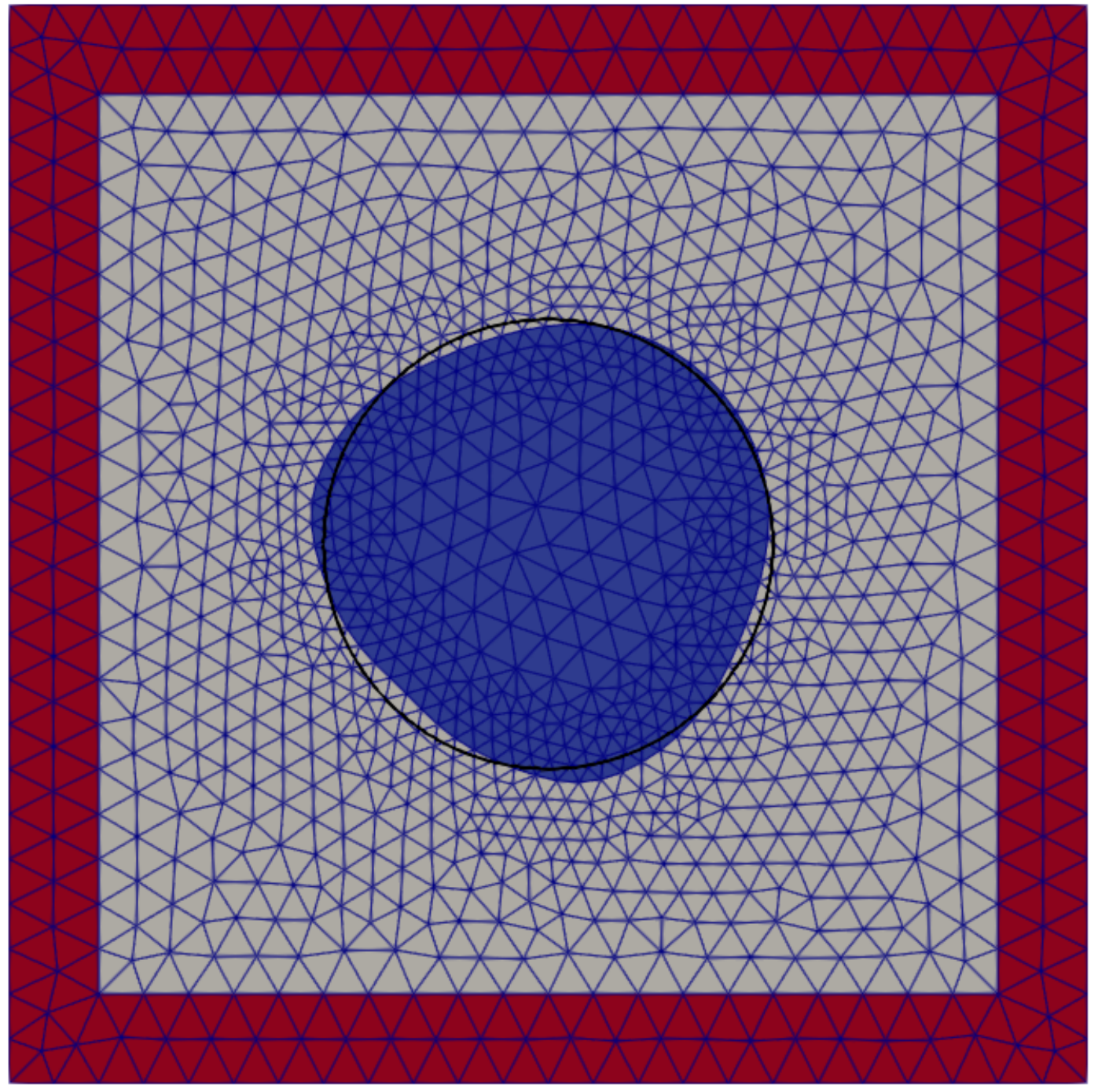}
		&\includegraphics[width = 0.2\textwidth]{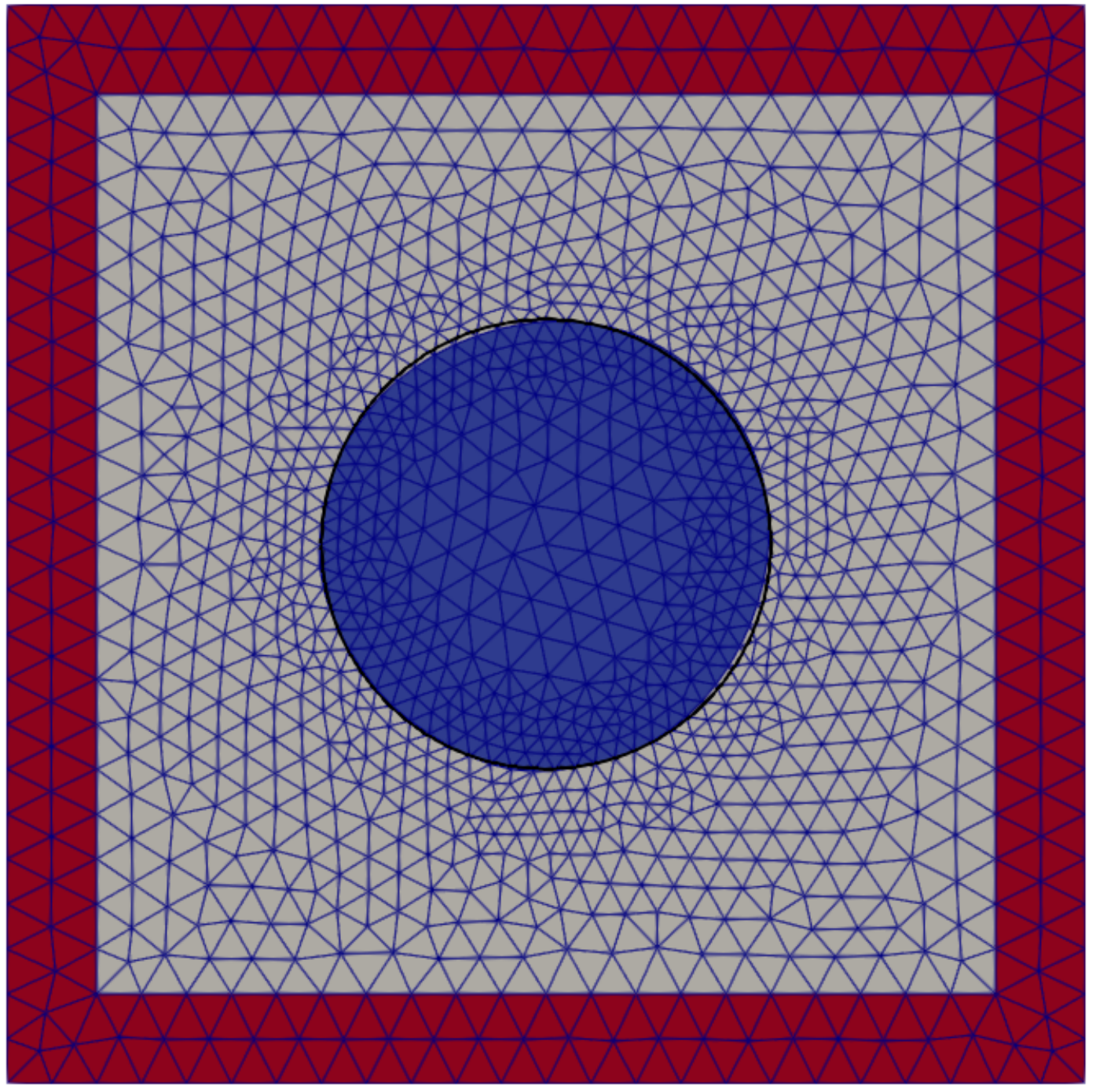}
		&\includegraphics[width = 0.2\textwidth]{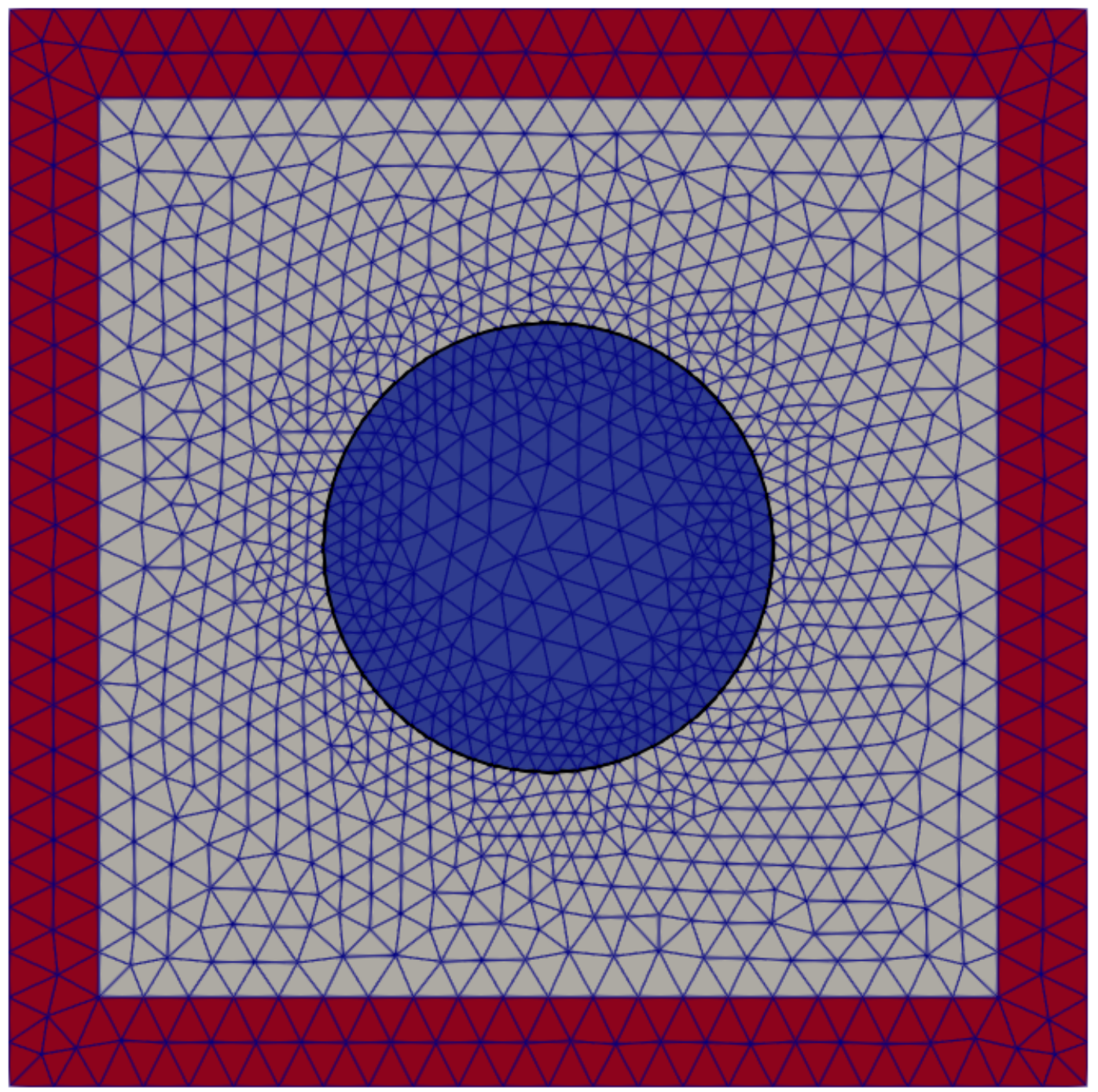}\\
		Iteration 10 & Iteration 10 remeshed & Iteration 15 & Iteration 20
	\end{tabular}
	\end{small}
	\end{center}
	\caption{}
	\label{fig:shape_ex_2_nonsym}
	~\\
\end{figure}

\begin{figure}[h!] 
	\centering
	%
	%
	\textbf{\textsf{Example 2: singular symmetric kernel}}
	\begin{center}
	\begin{small}
	\begin{tabular}{cccc}
		\includegraphics[width = 0.2\textwidth]{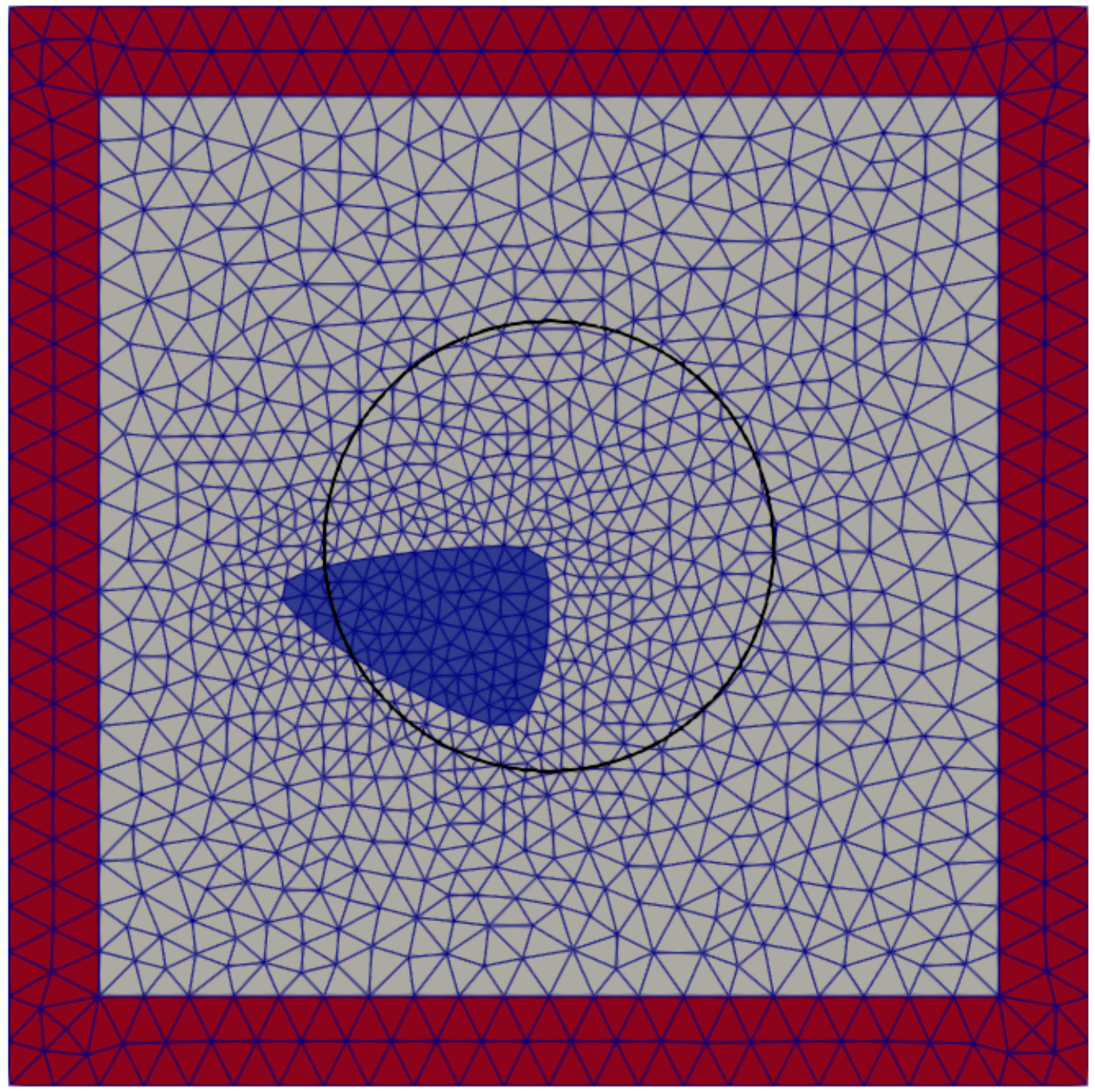}
		&\includegraphics[width = 0.2\textwidth]{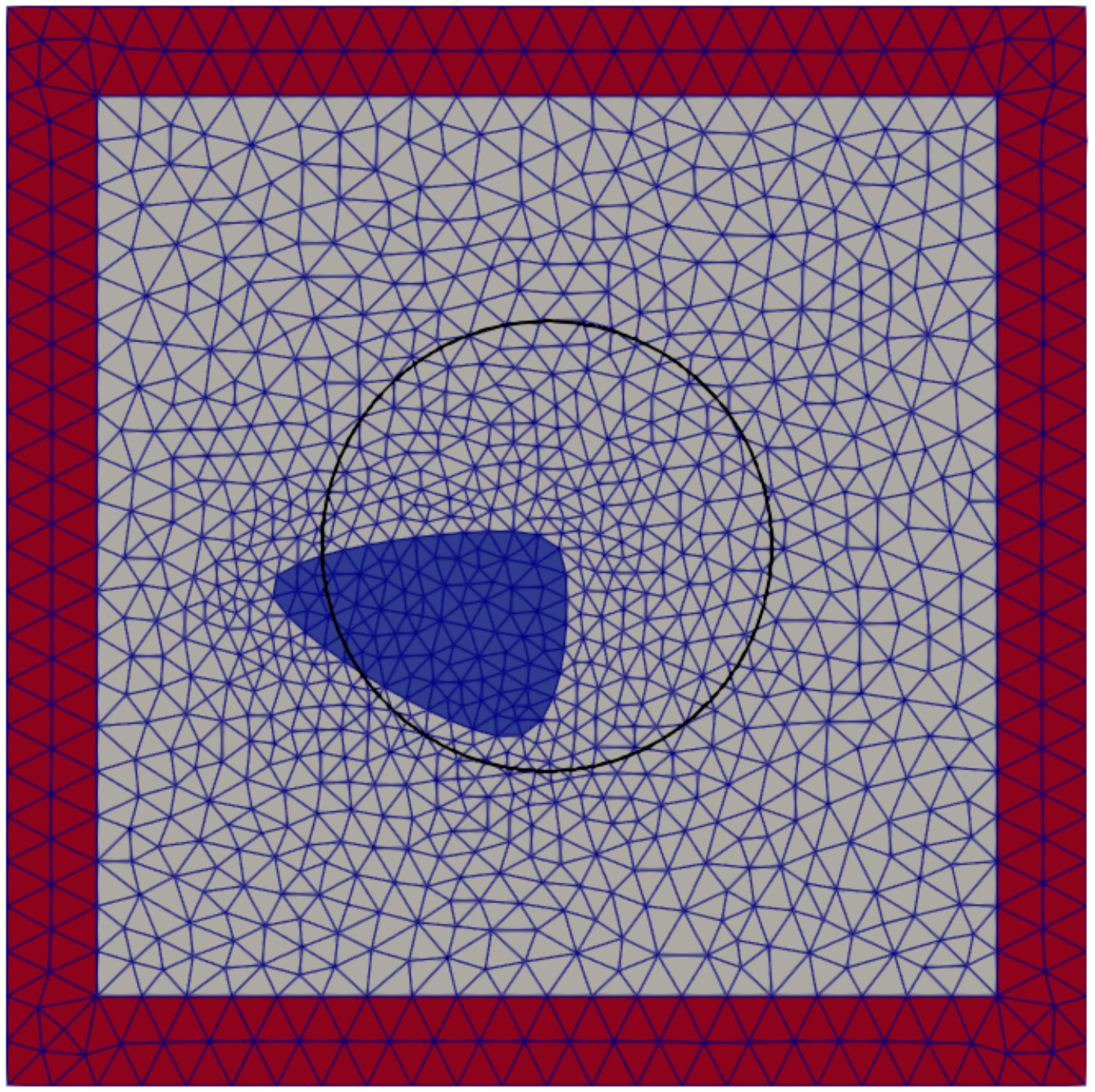}
		&\includegraphics[width = 0.2\textwidth]{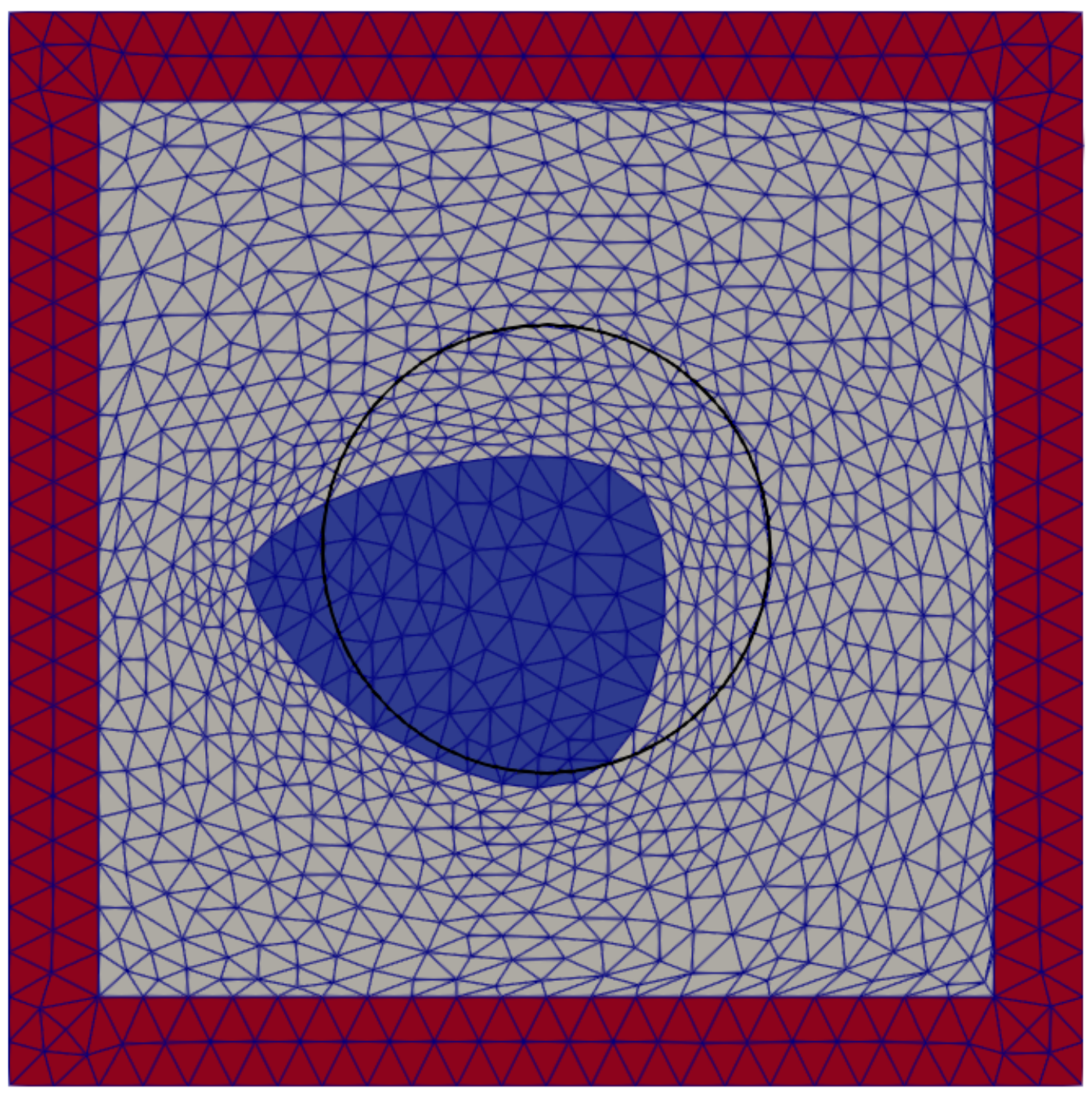}
		&\includegraphics[width = 0.2\textwidth]{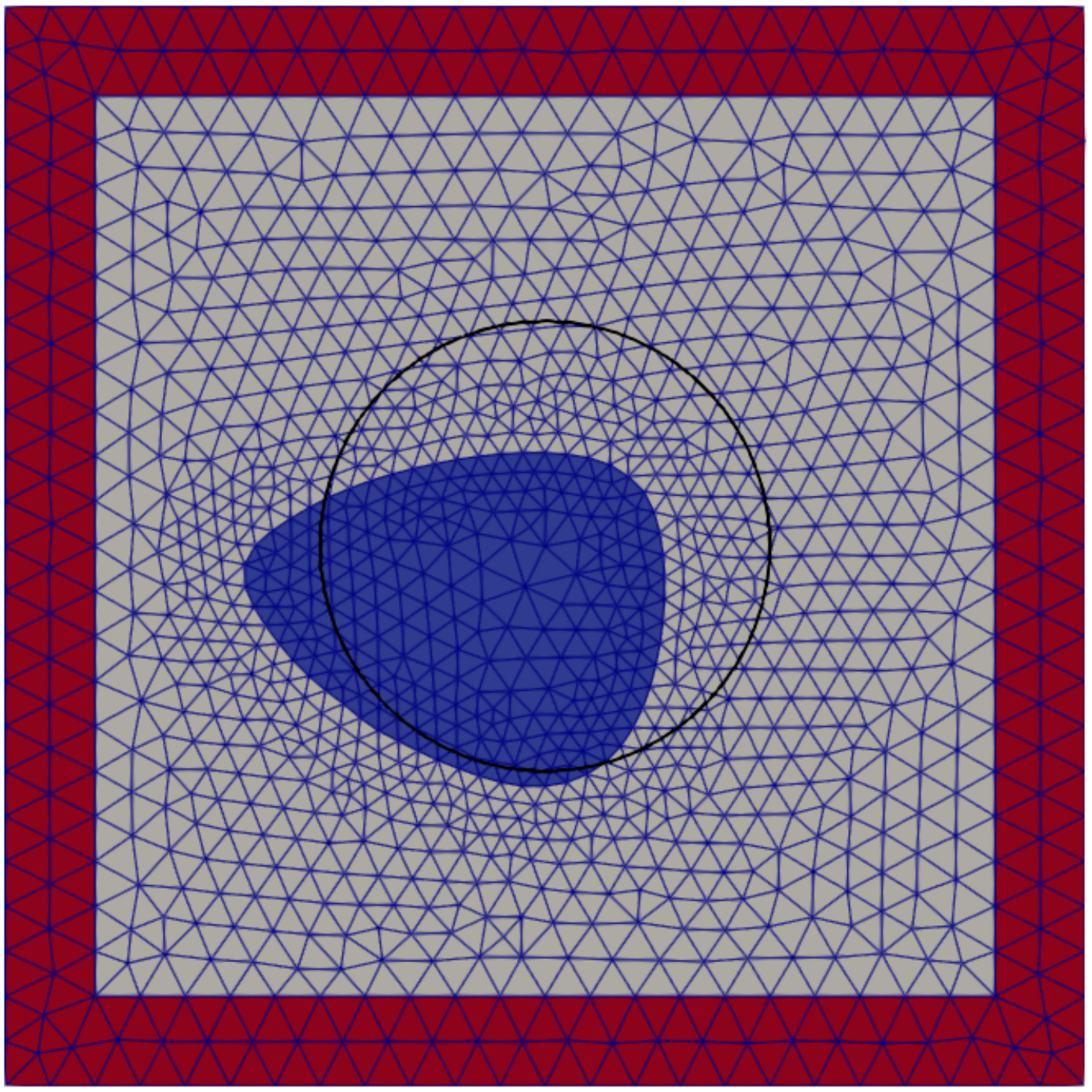}\\
		Start setup & Iteration 1 & Iteration 5 & Iteration 5 remeshed \\
		\includegraphics[width = 0.2\textwidth]{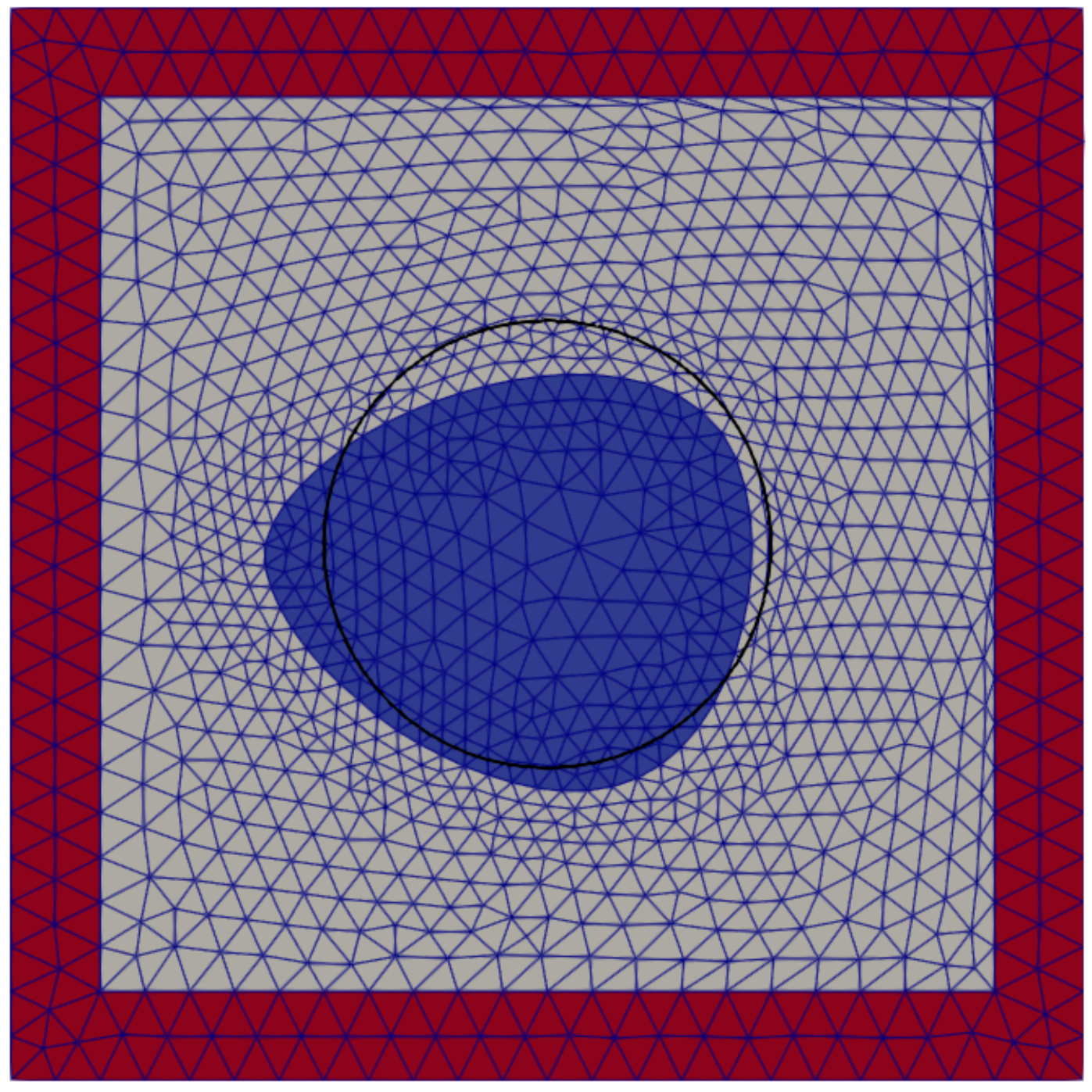}
		&\includegraphics[width = 0.2\textwidth]{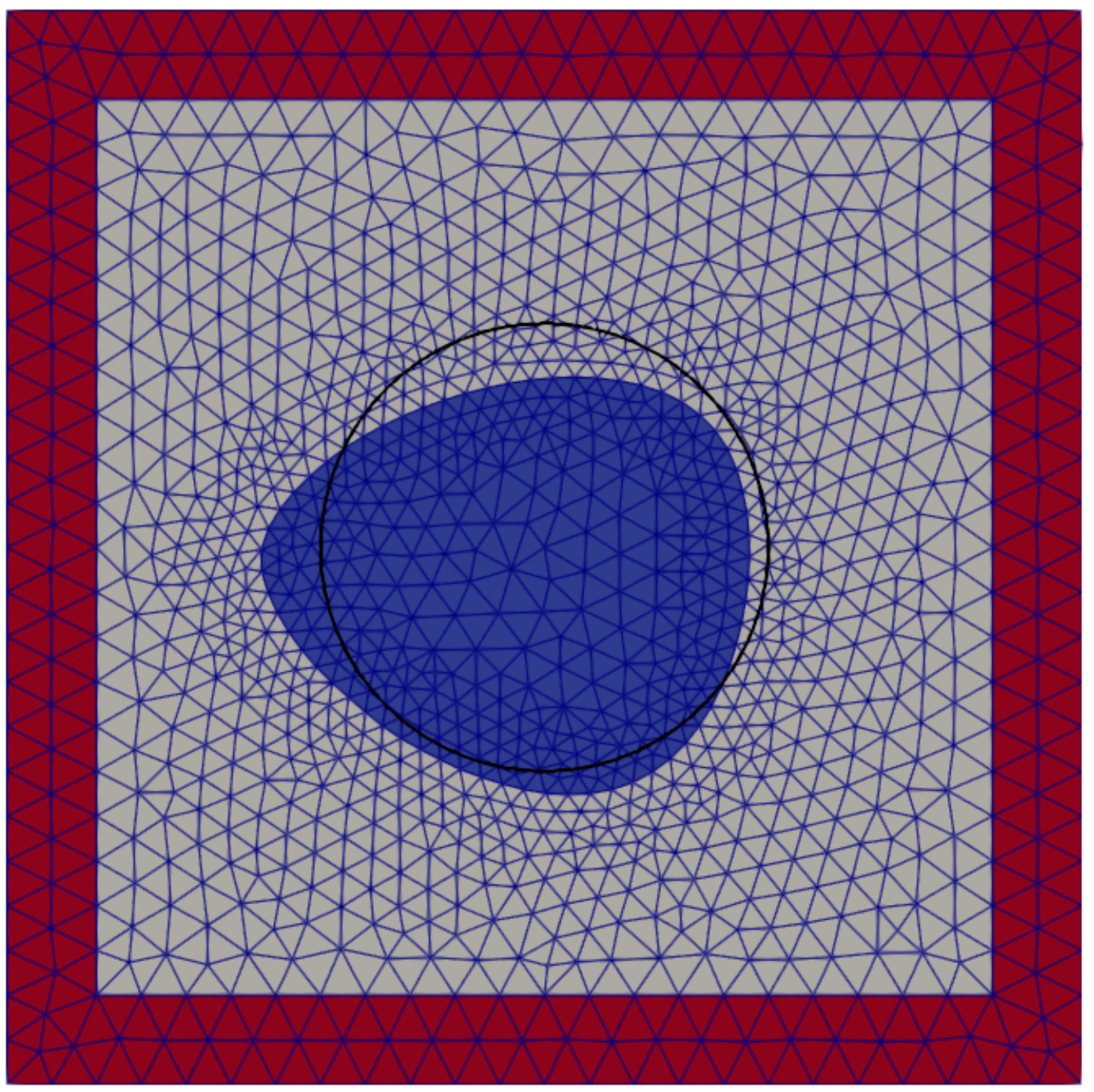}
		&\includegraphics[width = 0.2\textwidth]{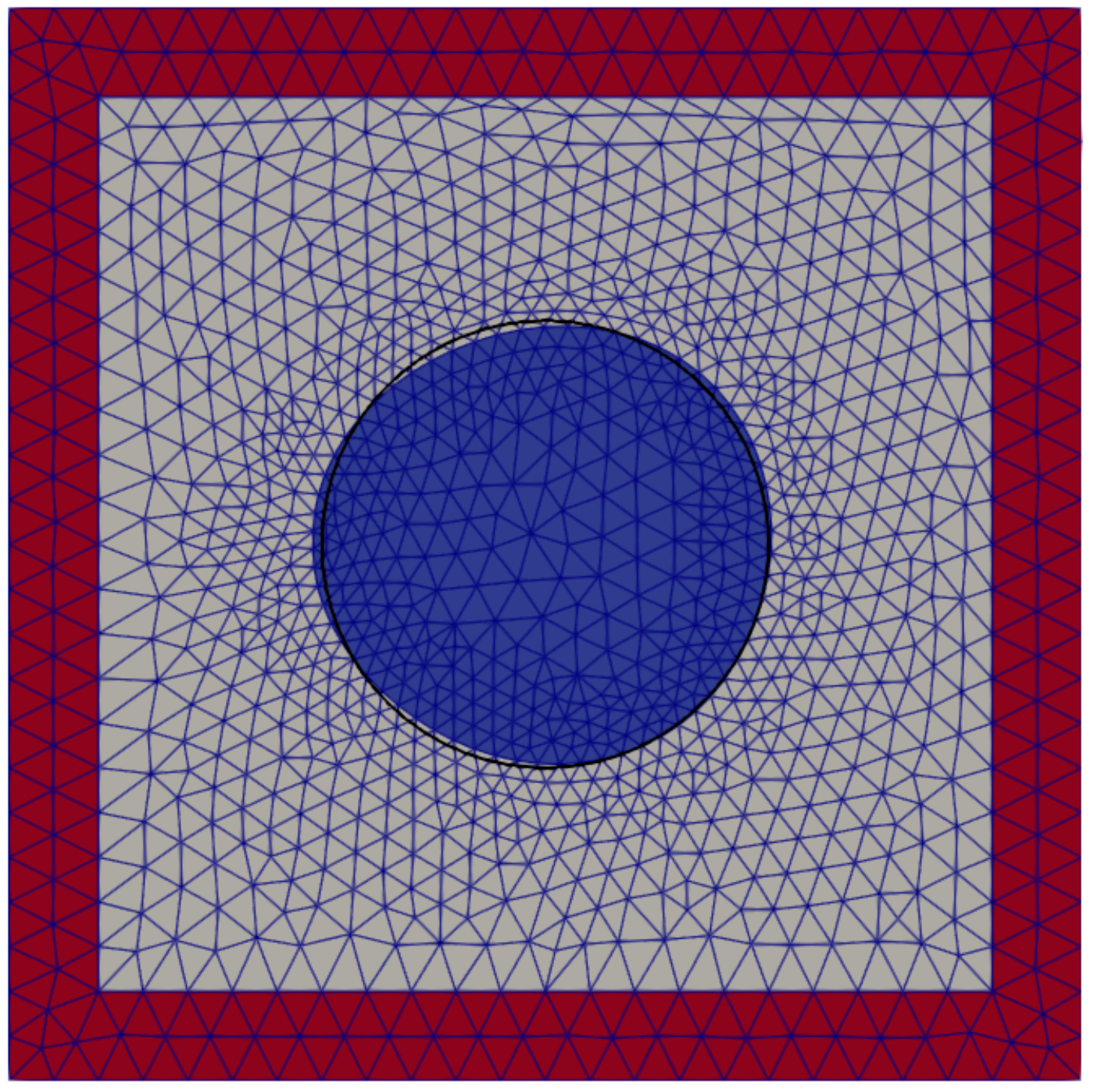}
		&\includegraphics[width = 0.2\textwidth]{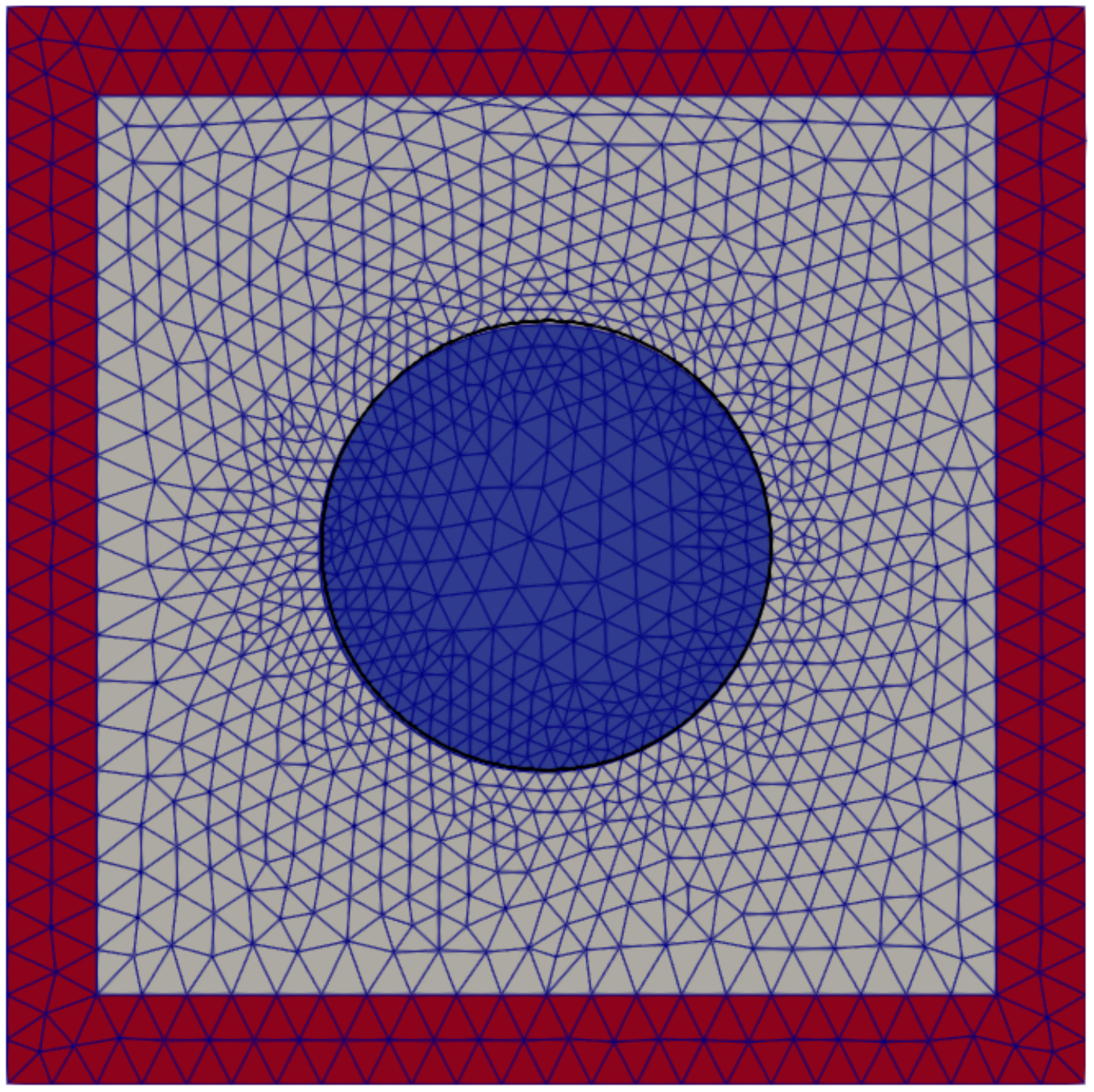}\\
		Iteration 10 & Iteration 10 remeshed & Iteration 15 & Iteration 20
	\end{tabular}
	\end{small}
	\end{center}
	\caption{}
	\label{fig:shape_ex_3_sym}
	%
\end{figure}

\begin{figure}[h!]
	\centering
	%
	%
	\textbf{\textsf{Example 2: nonsymmetric integrable kernel}}
	\begin{center}
	\begin{small}
	\begin{tabular}{cccc}
	     \includegraphics[width = 0.2\textwidth]{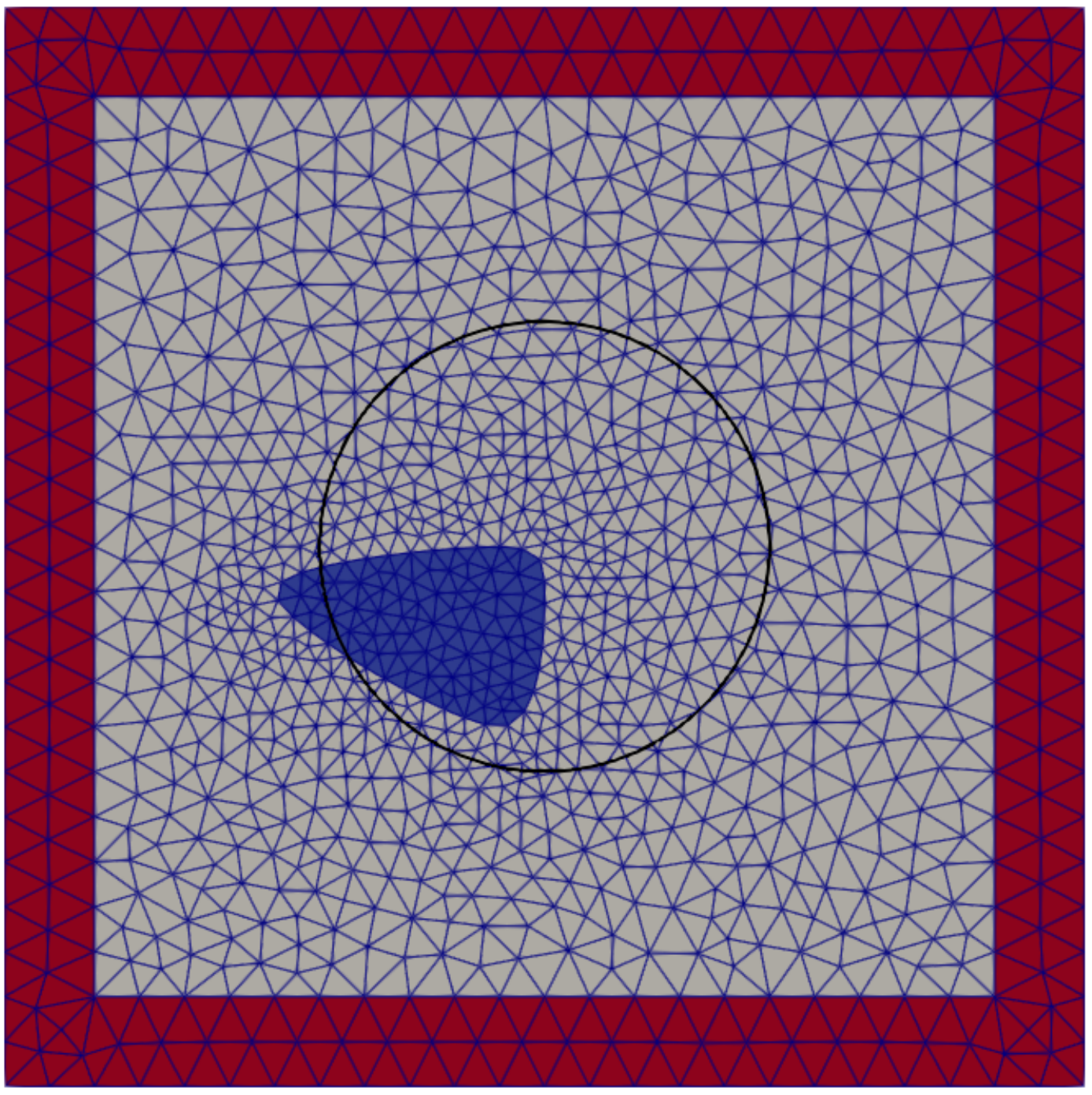}
		&\includegraphics[width = 0.2\textwidth]{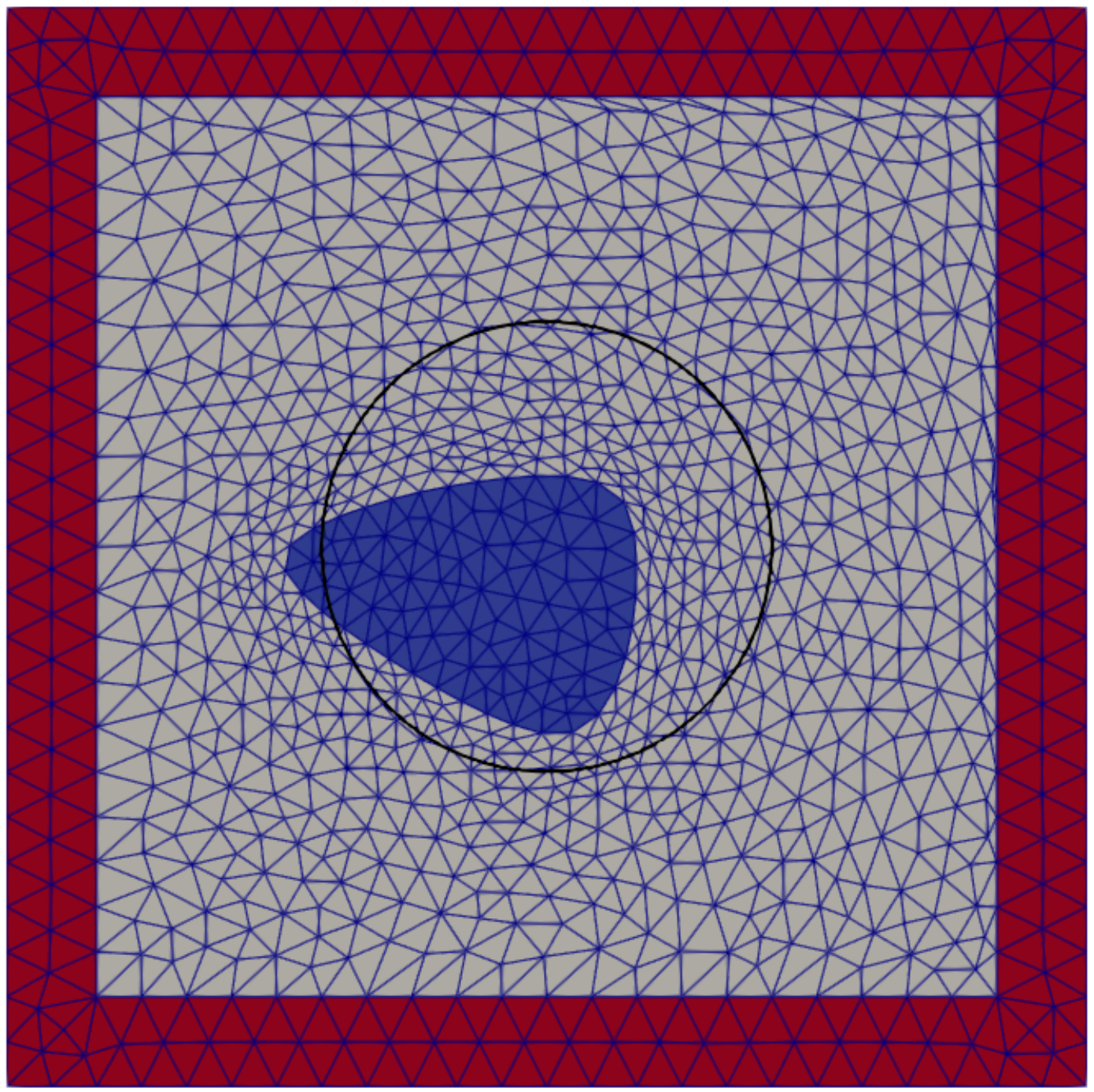}
		&\includegraphics[width = 0.2\textwidth]{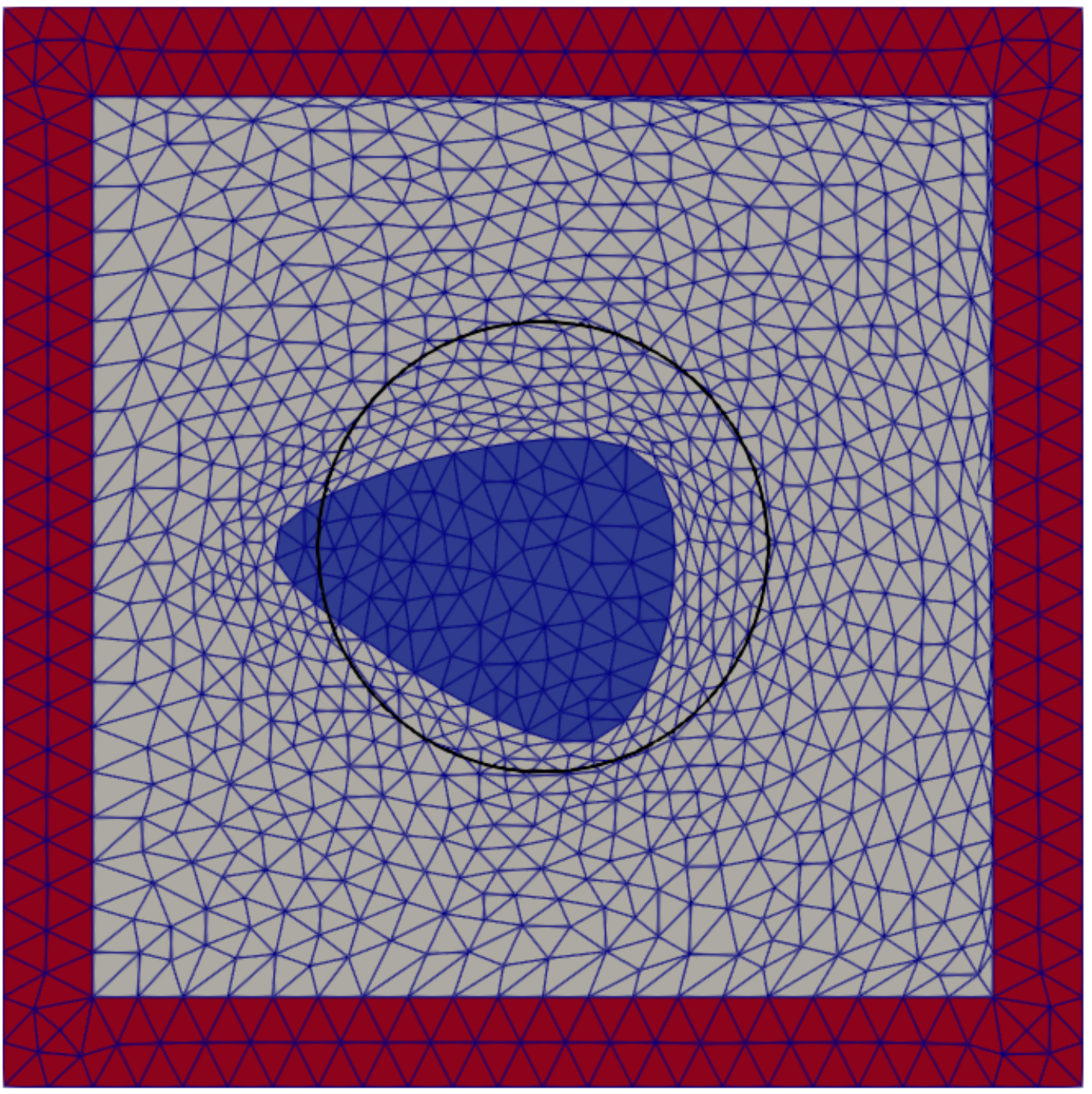}
		&\includegraphics[width = 0.2\textwidth]{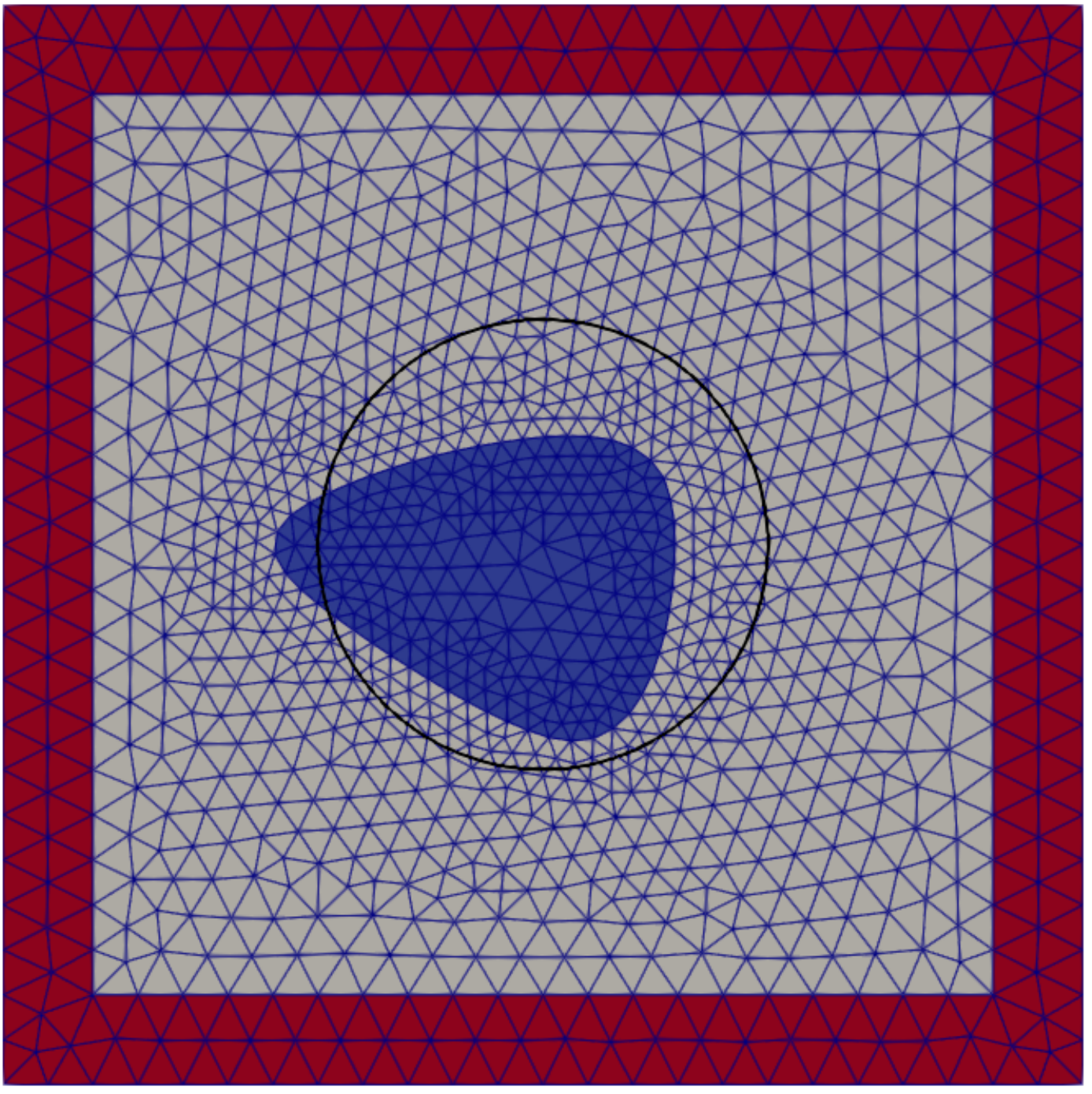}\\
	    Start setup & Iteration 1 & Iteration 5 & Iteration 5 remeshed\\ 
		\includegraphics[width = 0.2\textwidth]{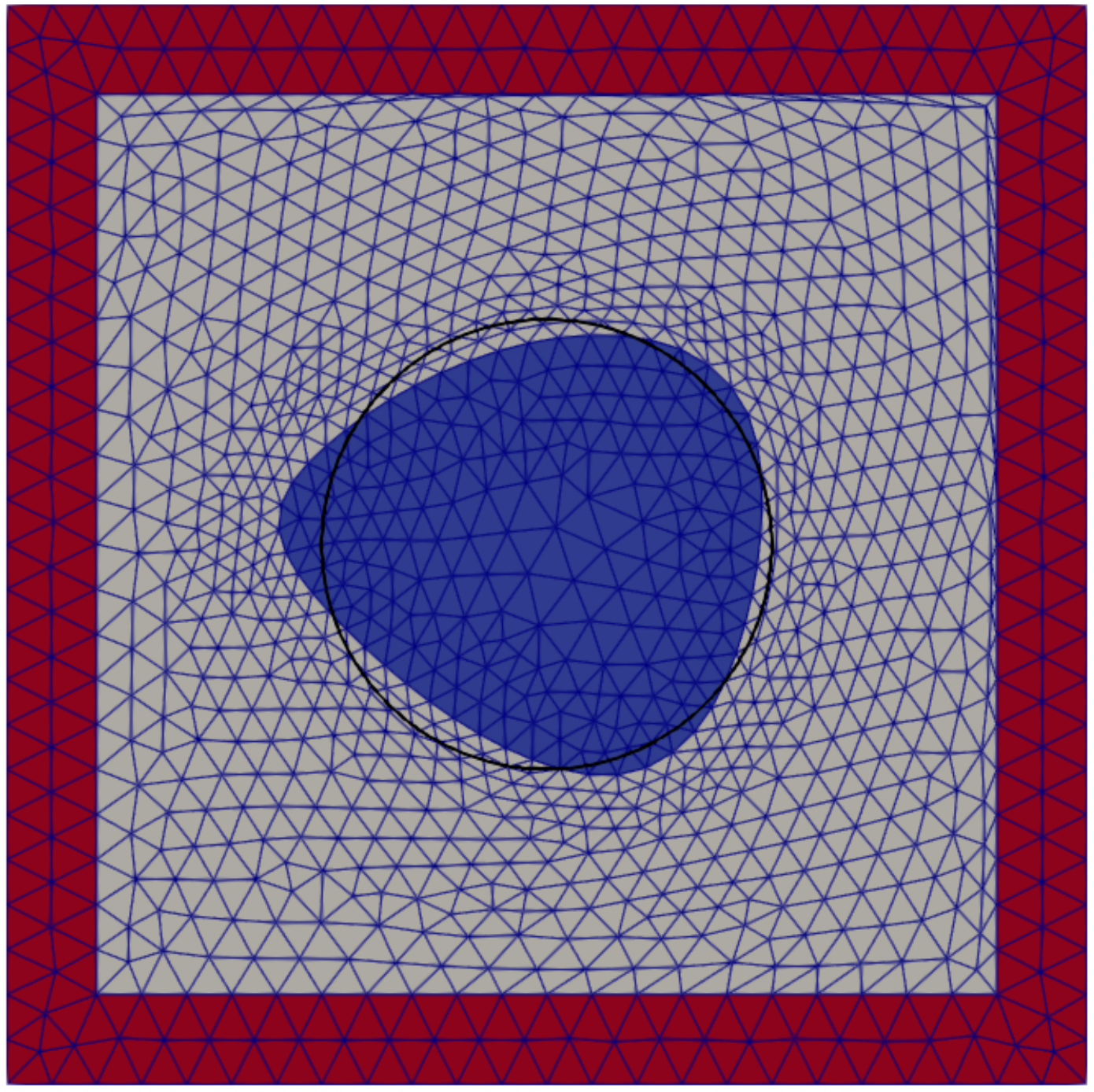}
		&\includegraphics[width = 0.2\textwidth]{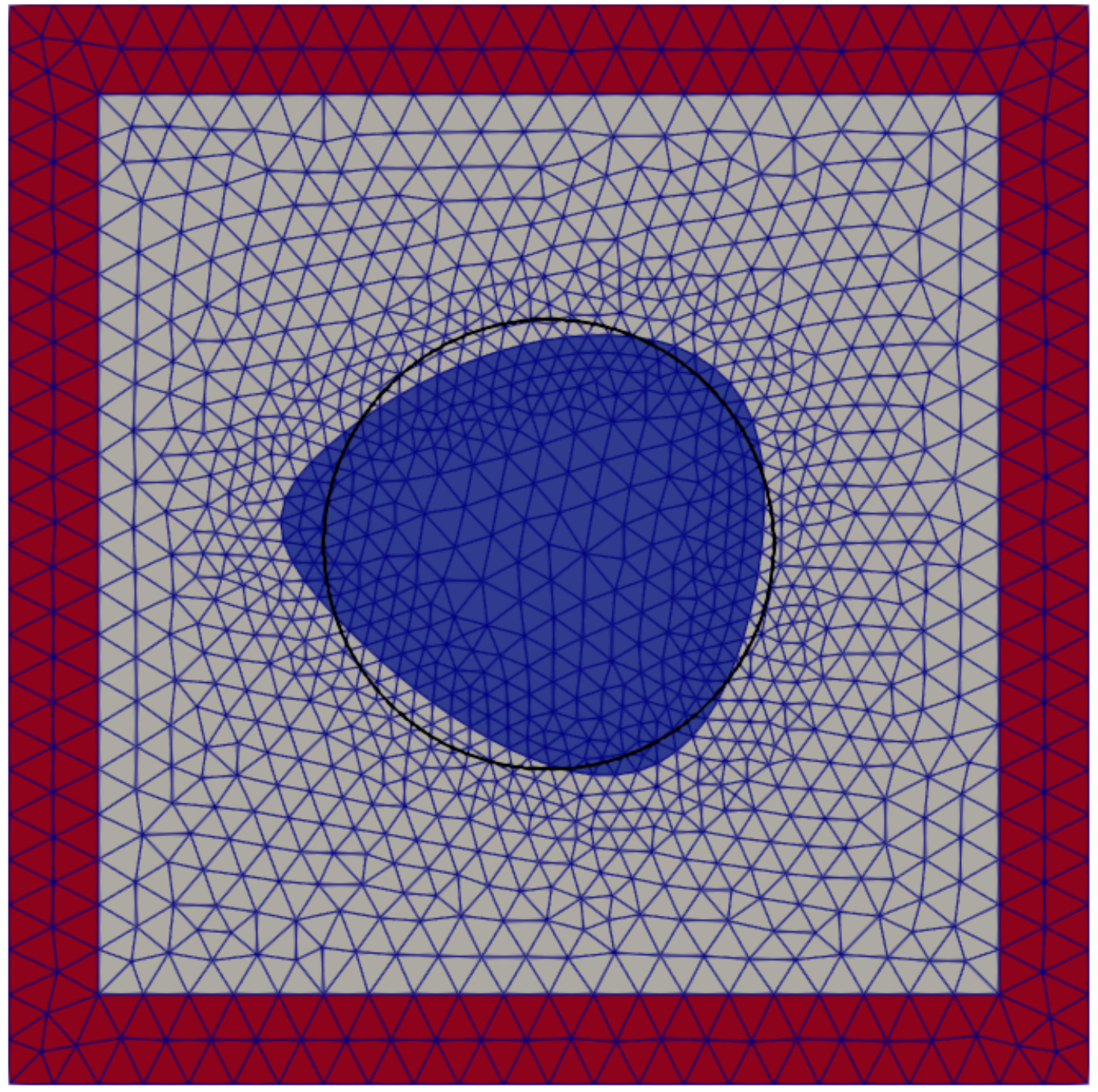}
		&\includegraphics[width = 0.2\textwidth]{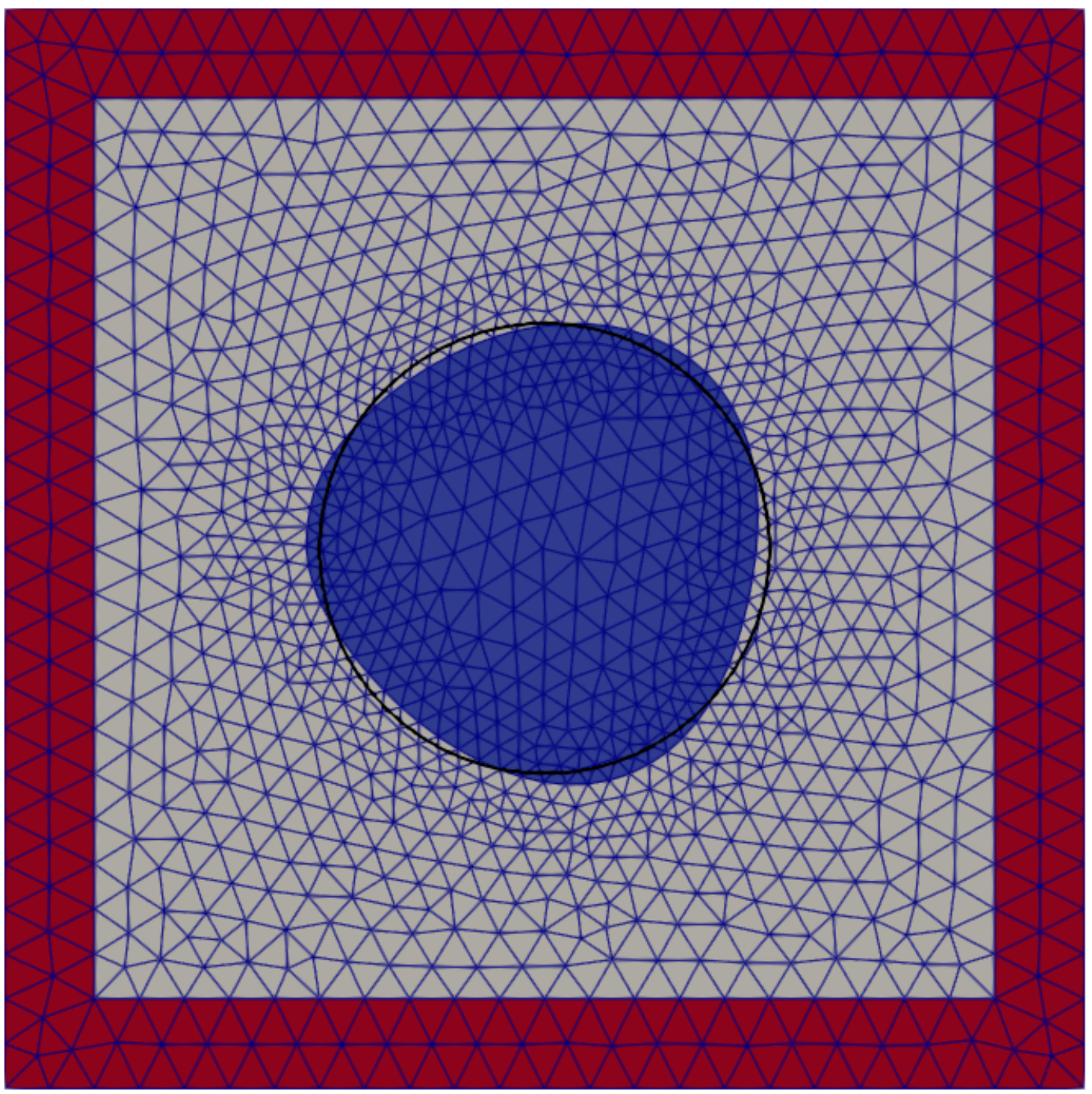}
		&\includegraphics[width = 0.2\textwidth]{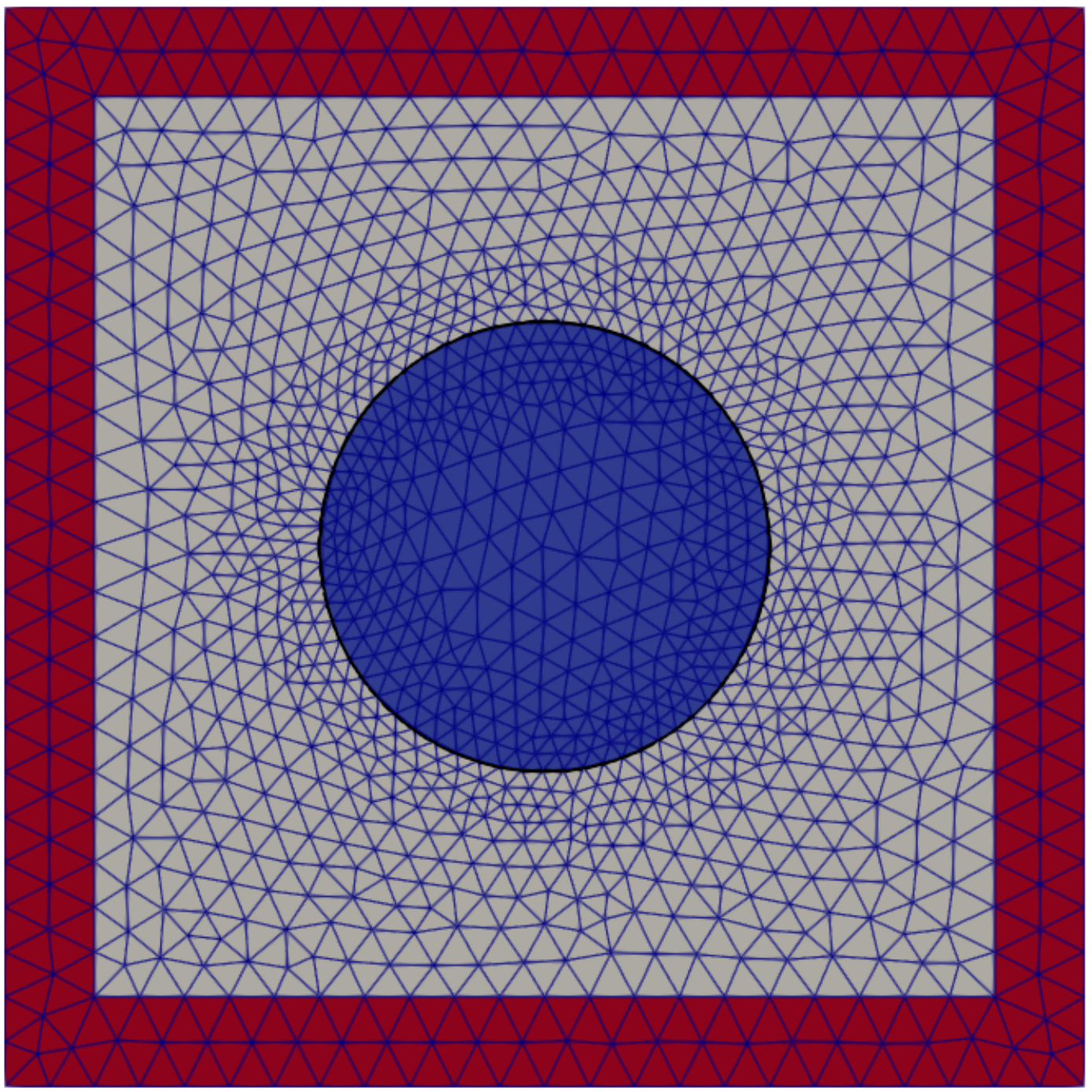}\\
		Iteration 10 & Iteration 10 remeshed & Iteration 15 & Iteration 25 
	\end{tabular}
	\end{small}
	\end{center}
	\caption{}
	\label{fig:shape_ex_3_nonsym}
	%
\end{figure}

We now present two different non-trivial examples which differ in the choice of the initial guess $\shape_0$. They are presented and described in the Figures \ref{fig:shape_ex_2_sym} and \ref{fig:shape_ex_3_sym} for the singular symmetric kernel and in the Figures \ref{fig:shape_ex_2_nonsym} and \ref{fig:shape_ex_3_nonsym} for the nonsymmetric integrable kernel. In each plot of the aforementioned figures the black line represents the target interface $\overline{\shape}$. Moreover the blue area depicts $\Omega_1$, the grey area $\Omega_2$ and the red area the nonlocal interaction domain $\Omega_I$.

Since the start shapes are smaller than the target shape, the shape needs to expand in the first few iterations. Thereby the nodes of the mesh are pushed towards the boundary, so that the mesh quality decreases and the algorithm stagnates, because nodes are prohibited to be pushed outside of $\Omega$. Therefore, we apply a re-meshing technique, where we re-mesh after the fifth and tenth iteration. In order to re-mesh, we save the points of our current shape as a spline in a dummy .geo file, that also contains the information of the nonlocal boundary, and then compute a new mesh with Gmsh. In this new mesh the distance between the nodes and the boundary is sufficiently large enough to attain a better improvement regarding the objective function value by the new mesh deformations.

It is important to mention that computation times and the performance of Algorithm \ref{alg:shape_optimization} in general are very sensitive to the choice of parameters and may strongly vary, which is why reporting exact computation times is not very meaningful at this stage. Particularly delicate choices are those of the system parameters including the kernel (diffusion and convection) and the forcing term, which both determine the identifiability of the model. But also the choice of Lam\'e parameters to control the step size, specifically $\mu_{max}$ (we set $\mu_{min} = 0$ in all experiments, since we want the boundary of $\Omega$ to be fixed).

Moreover, especially in the case of system parameters with high interface-sensitivity in combination with an inconveniently small $\mu_{max}$, mesh deformations may be large in the early phase of the algorithm. Thus, such mesh deformations $\tilde{\Ub}_k$ of high magnitude lead to destroyed meshes so that an evaluation of the reduced objective functional $J^{red}((id + \alpha\tilde{\Ub}_k)(\Omega_k))$, which requires the assembly of the nonlocal stiffness matrix, becomes a pointless computation. In order to avoid such computations we first perform a line search depending on one simple mesh quality criterion. More precisely, we downscale the step size, i.e., $\alpha = \tau \alpha$, until all finite element nodes of the resulting mesh $(id+\alpha\tilde{\Ub}_k)(\Omega_k)$ are a subset of $\Omega$. After that, we continue with the backtracking line search in Line \ref{shape_algo:line search} of Algorithm \ref{alg:shape_optimization}.
\begin{figure}[h!]
    \centering
    \includegraphics[width = 0.8\textwidth]{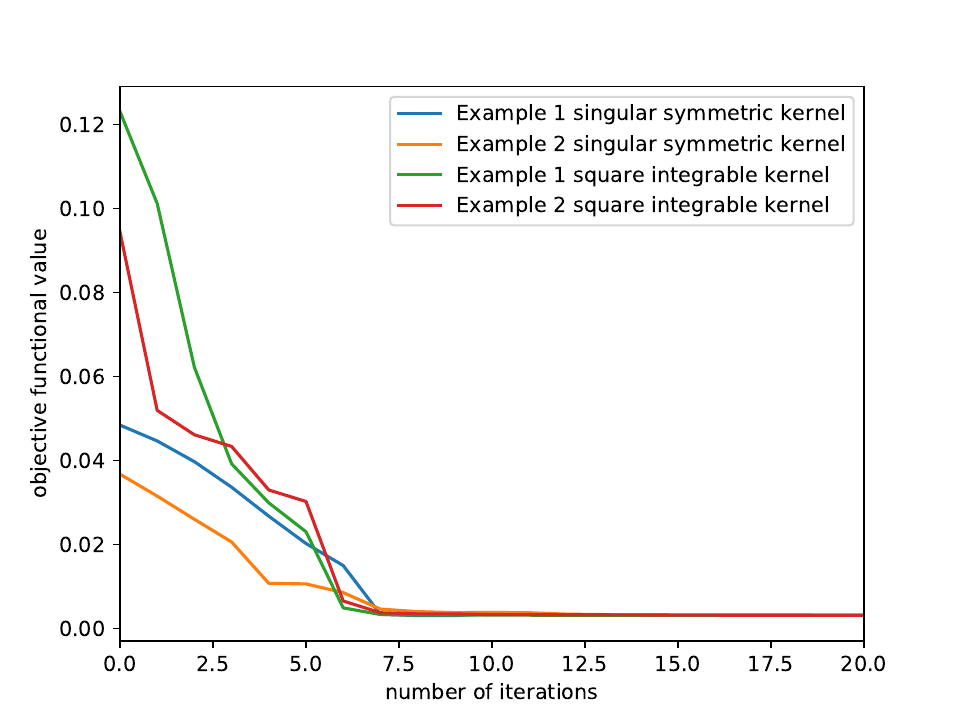}
    \caption{In the first six or seven iterations the improvement regarding the objective function value is quite high. After that the objective function value decreases in a much slower fashion. Due to the regularization term the objective functional will not converge to zero.}
    \label{fig:comparison_objective_funtional}
\end{figure}
\section{Concluding remarks and future work} \label{sec:shape_conclusion_outlook}
We have conducted a numerical investigation of shape optimization problems which are constrained by nonlocal system models. We have proven through numerical experiments the applicability of established shape optimization techniques for which the shape derivative of the nonlocal bilinear form represents the novel ingredient. All in all, this work is only a first step along the exploration of the interesting field of nonlocally constrained shape optimization problems.

\section*{Acknowledgement}
This work has been partly supported by the German Research Foundation
(DFG) within the Research Training Group 2126: ‘Algorithmic Optimization’.

\appendix
\section{Proof of Lemma \ref{lem:AAM_nonlocal_shape_problem}}
\label{chap:Proof_AAM_nl_shape_problem}
\begin{proof}
Define $\xt(\xb) \defas \det DF_t(\xb),\ \gt_\shape(\xb,\yb) \defas \g_\shape(F_t(\xb),F_t(\yb)),\ f_\shape^t(\xb) \defas f_\shape(F_t(\xb))$\\
and ${\nabla f_\shape^t(\xb) \defas \nabla f_\shape(F_t(\xb))}$.\\
\textbf{Assumption(H0):}\\
$G(t,su^t + (1-s)u^0,\advar)$ is absolutely continuous in $s$, if there exists a function $g \in L^1[0,1]$, such that
\begin{align*}
G(t,su^t + (1-s)u^0,\advar) = G(t,u^0,\advar) + \int_0^s g(\tilde{s}) d\tilde{s}.
\end{align*}
We can compute
\begin{align*}
&G(t,su^t + (1-s)u^0,\advar) = \AAA(t,su^t + (1-s)u^0, \advar) - F(t,\advar) + J(t,su^t + (1-s)u^0) \\
&= \frac{1}{2} \int_\Omega \int_\Omega (\advar(\xb) - \advar(\yb))((su^t + (1-s)u^0)(\xb)\gt_\shape(\xb,\yb)\\ 
&~~~~~~~~~~~~~ - (su^t + (1-s)u^0)(\yb)\gt_\shape(\yb,\xb))\xt(\xb)\xt(\yb)~d\yb d\xb + \int_\Omega \int_{\Omega_I} (su^t + (1-s)u^0) \advar \xt ~d\yb d\xb \\
&~~ - \int_\Omega f_\shape^t\advar\xt ~d\xb +\frac{1}{2} \int_\Omega (s(u^t -u^0) + (u^0 - \bar{u}))^2\xt ~d\xb + \int_{\shape^t} 1 ds \\
&= s\AAA(t,u^t-u^0,\advar) + \AAA(t,u^0,\advar) - \int_\Omega f_\shape^t\advar \xt ~d\xb + \frac{1}{2} s^2 \int_\Omega (u^t-u^0)^2 \xt ~d\xb \\ 
&~~~+ s \int_\Omega (u^t - u^0)(u^0 - \bar{u})\xt ~d\xb + \frac{1}{2} \int_\Omega (u^0 - \bar{u})^2 \xt ~d\xb  + \int_{\Gamma_t} 1 ~ds \\
&= \int_0^s \tilde{s} \int_\Omega (u^t-u^0)^2\xt ~d\xb + \AAA(t,u^t - u^0,\advar) + \int_\Omega (u^t-u^0)(u^0 - \bar{u})\xt ~d\xb d\tilde{s} \\
&~~~ + \AAA(t,u^0,\advar) - \int_\Omega f_\shape^t\advar\xt ~d\xb + \frac{1}{2} \int_\Omega (u^0 - \bar{u})^2 \xt ~d\xb  + \int_{\Gamma_t} 1 ~ds = \int_0^s M_1\tilde{s} + M_2 ~d\tilde{s} + G(t,u^0,\advar),
\end{align*}
where
\begin{align*}
M_1 := \int_\Omega (u^t-u^0)^2\xt ~d\xb \text{ and } M_2:=\AAA(t,u^t-u^0,\advar) + \int_\Omega (u^t-u^0)(u^0 - \bar{u})\xt ~d\xb.
\end{align*} 
Furthermore the second criterion of (H0) is also satisfied:
\begin{align*}
&\int_0^1 |d_u G(t,su^t + (1-s)u^0,\advar)[\tilde{u}]| ~ds \\
&= \int_0^1 |\AAA(t,\tilde{u},\advar) + \int_\Omega (su^t + (1-s)u^0 - \bar{u})\tilde{u}\xt ~d\xb|ds \\
&\leq |\AAA(t,\tilde{u},\advar)| + \int_0^1 s ~ds |\int_\Omega (u^t - u^0)\tilde{u}\xt ~d\xb| +  |\int_\Omega (u^0 - \bar{u})\tilde{u}\xt ~d\xb| \\
&= |\AAA(t,\tilde{u},\advar)| + \frac{1}{2} |\int_\Omega (u^t - u^0)\tilde{u}\xt ~d\xb| + |\int_\Omega (u^0 - \bar{u})\tilde{u}\xt ~d\xb|  < \infty.
\end{align*}
As we saw above in the proof of assumption (H0), the left hand side of the average adjoint equation can be written as follows
\begin{align*}
\int_0^1 d_u G(t,su^t + (1-s)u^0,\advar^t)[\tilde{u}] ~ds = \AAA(t,\tilde{u},\advar^t) + \int_\Omega \left(\frac{1}{2}(u^t + u^0) - \bar{u}\right)\tilde{u}\xt ~d\xb.
\end{align*}
Then, (\ref{AAE}) can be reformulated as
\begin{align}
\label{nonlocalAAE}
\AAA(t,\tilde{u},\advar^t) = - \int_\Omega \left(\frac{1}{2}(u^t + u^0) - \bar{u}\right)\tilde{u}\xt ~d\xb \quad \text{for all } \tilde{u} \in V_c(\Omega).
\end{align}
Since the right hand side of (\ref{nonlocalAAE}) is a linear and continuous operator with regards to $\tilde{u}$, equation (\ref{nonlocalAAE}) is a well-defined nonlocal problem, which has a unique solution $\advar^t$ due to the assumptions (P0). 
By further using assumptions (P0), we can conclude for $u^t$, that there exists a $C_2 > 0$, such that
\begin{align*}
    ||u^t||^2_{L^2(\Omega)} &\leq C_1 |A(t,u^t,u^t)| = C_1 |(f^t\xt,u^t)_{L^2(\Omega)}| \leq C_1 ||f^t \xt||_{L^2(\Omega)} ||u^t||_{L^2(\Omega)} \\
    \Rightarrow ||u^t||_{L^2(\Omega)} &\leq C_1 ||f^t \xt||_{L^2(\Omega)} \leq C_2,
\end{align*}
where we used in the last step, that $||f_\shape^t\xt||_{L^2(\Omega)} \rightarrow ||f_\shape||_{L^2(\Omega)}$ (\cite[Lemma 2.16]{Sturm_diss}). Since ${\xt(\xb) = \det (I + tDV(\xb))}$ is continuous on $[0,T]\times \bar{\Omega}$, there exists a $\bar{\xi}$, such that $|\xt(\xb)|\leq \bar{\xi}$ for all $(t,\xb) \in [0,T]\times \bar{\Omega}$. Then we derive for $\advar^t$ that
\begin{align*}
    ||\advar^t||^2_{L^2(\Omega)} &\leq C_1 |A(t,\advar^t,\advar^t)| = C_1 |(\frac{1}{2}(u^t+u^0) - \bar{u})\xt,\advar^t)_{L^2(\Omega)}| \\ 
    &\leq \bar{\xi} C_1 ||\frac{1}{2}(u^t+u^0) - \bar{u}||_{L^2(\Omega)} ||\advar^t||_{L^2(\Omega)} \\
    \Rightarrow ||\advar^t||_{L^2(\Omega)} &\leq \bar{\xi} C_1 ||\frac{1}{2}(u^t+u^0) - \bar{u}||_{L^2(\Omega)} \leq \bar{\xi} C_1 ||\bar{u}||_{L^2(\Omega)} + \frac{1}{2} \bar{\xi} C_1  (||u^t||_{L^2(\Omega)} + ||u^0||_{L^2(\Omega)}) \\
    &\leq \bar{\xi} C_1  ||\bar{u}||_{L^2(\Omega)} + \bar{\xi}C_1C_2.
\end{align*}
\textbf{Assumption (H1):}\\
Since $u^t$ and $\advar^t$ are bounded for $t \in [0,T]$, then for every sequence $\{t_n\}_{n \in \Nb}$ with $t_n \rightarrow 0$ there exist subsequences $\{t_{n_k}\}_{k \in \Nb}$ and $\{t_{n_l}\}_{l \in \Nb}$, such that there exist functions $q_1, q_2 \in L^2(\OuO)$ with $u^{t_{n_l}} \rightharpoonup q_1$ $\advar^{t_{n_k}} \rightharpoonup q_2$ in $L^2(\OuO)$. 
In the following we use that for functions \\${\{g^t\}_{t \in [0,T]},\{h^t\}_{t \in [0,T]} \in L^2(\OuO)^{[0,T]}}$ with $g^t \rightarrow g^0$ and $h^t \rightharpoonup h^0$ in $L^2(\OuO)$, we obtain
\begin{align}
\label{weak_convergence_lemma}
\begin{split}
    \left| \int_\Omega h^t g^t ~d\xb - \int_\Omega h^0 g^0 ~d\xb \right| \leq \left| \int_\Omega h^t(g^t - g^0) ~d\xb \right| + \left| \int_\Omega (h^t - h^0)g^0 ~d\xb \right| \\ 
    \leq ||h^t||_{L^2(\Omega)}||g^t-g^0||_{L^2(\Omega)} + |\int_\Omega (h^t - h^0)g^0 ~d\xb| \rightarrow 0.
\end{split}
\end{align}
\textbf{Case 1: Proof of (H1) for square integrable kernels}\\
Since $\phi_{ij}$ is essentially bounded on $\Omega \times \Omega$, $\phi_{iI}$ is essentially bounded on $\Omega \times \Omega_I$ and $\xt$ is continuous and therefore bounded on $\bar{\Omega}$, we can conclude that $$\psi_t(\yb) \defas (\advar(\xb) - \advar(\yb))\gt_\shape(\xb,\yb)\xt(\xb)\xt(\yb) \in L^2(\OuO)$$ and by using the dominated convergence theorem we get $$\psi_t(\yb) \rightarrow  (\advar(\xb) - \advar(\yb))\g_\shape(\xb,\yb) \text{ in } L^2(\OuO) \text{ for } t \searrow 0.$$
Thus, by (\ref{weak_convergence_lemma}) we derive
\begin{align*}
&\int_\Omega(\advar(\xb) - \advar(\yb))(u^{t_{n_l}}(\xb)\g_\shape^{t_{n_l}}(\xb,\yb) - u^{t_{n_l}}(\yb)\g_\shape^{t_{n_l}}(\yb,\xb))\xi^{t_{n_l}}(\xb)\xi^{t_{n_l}}(\yb)~d\yb \\ 
&\rightarrow \int_\Omega(\advar(\xb) - \advar(\yb))(u^0(\xb)\g_\shape(\xb,\yb) - u^0(\yb)\g_\shape(\yb,\xb))~d\yb \text{ for } \xb \in \Omega, \advar \in L^2_c(\OuO).
\end{align*}
Due to the continuity of parameter integrals, we have
\begin{align*}
    &\int_\Omega \int_\Omega (\advar(\xb) - \advar(\yb))(u^{t_{n_l}}(\xb)\g_\shape^{t_{n_l}}(\xb,\yb) - u^{t_{n_l}}(\yb)\g_\shape^{t_{n_l}}(\yb,\xb))\xi^{t_{n_l}}(\xb)\xi^{t_{n_l}}(\yb)~d\yb d\xb\\ 
    &\rightarrow \int_\Omega \int_\Omega (\advar(\xb) - \advar(\yb))(u^0(\xb)\g_\shape(\xb,\yb) - u^0(\yb)\g_\shape(\yb,\xb))~d\yb d\xb.
\end{align*}
Analogously, we can show
\begin{align*}
    \int_\Omega u^{t_{n_l}}(\xb)\advar(\xb) \int_{\Omega_I} \g_\shape^{t_{n_l}}(\xb,\yb)\xi^{t_{n_l}}(\xb) ~d\yb d\xb \rightarrow \int_\Omega u^0(\xb)\advar(\xb) \int_{\Omega_I} \g_\shape(\xb,\yb) ~d\yb d\xb.
\end{align*}
So we can conclude $\lim_{l \rightarrow \infty} \AAA(t_{n_l},u^{t_{n_l}},\advar) = A(0,q_1,\advar)$.
Because $f_\shape^t \xt \rightarrow f_\shape$ in $L^2(\Omega)$ according to \cite[Lemma 2.16]{Sturm_diss}, we can compute for all $\advar \in L_c^2(\Omega)$ 
\begin{align*}
\AAA(0,q_1,\advar) = \lim_{l \rightarrow \infty} \AAA(t_{n_l},u^{t_{n_l}},\advar) = \lim_{l \rightarrow \infty} \int_\Omega f_\shape^{t_{n_l}} \advar \xi^{t_{n_l}} ~d\xb = \int_\Omega f_\shape \advar ~d\xb.
\end{align*}
Since the solution is unique we derive $q_1 = u^0$ and $u^t \rightharpoonup u^0$. Similarly, we have for $q_2$ and for all $\tilde{u} \in L^2_c(\Omega)$
\begin{align*}
\AAA(0, \tilde{u}, q_2) = \lim_{k \rightarrow \infty} \AAA(t_{n_k}, \tilde{u}, \advar^{t_{n_k}}) = - \lim_{k \rightarrow \infty} \int_\Omega (\frac{1}{2}(u^{t_{n_k}} + u^0) - \bar{u})\tilde{u} \xi^{t_{n_k}} ~d\xb = - \int_\Omega (u^0 - \bar{u})\tilde{u} ~d\xb.
\end{align*}
So we conclude $q_2=\advar^0$ and $\advar^t \rightharpoonup \advar^0(t \searrow 0)$. 
By using the mean value theorem, there exist $s_t \in (0,t)$, s.t. 
${s_t \rightarrow 0 (t \searrow 0)}$ and
\begin{align*} 
    \frac{G(t,u^0,\advar^t) - G(0,u^0,\advar^t)}{t} = \partial_t G(s_t,u^0,\advar^t).
\end{align*}
Therefore we now prove assumption (H1) by showing
\begin{align*}
    \lim_{s,t \searrow 0} \partial_t G(t,u^0,\advar^s) = \partial_t G(0,u^0,\advar^0).
\end{align*}
Computing the derivative regarding $t$ yields
\begin{align*}
\partial_t G(t,u^0,\advar^s) &= \partial_t \AAA(t,u^0,\advar^s) - \partial_t F(t,\advar^s) + \partial_t J(t, u^0).
\end{align*}
First we can show
\begin{align}
\begin{split}
\label{eq:AAM_H1_ell}
 \partial_t F(t,\advar^s) = \int_\Omega (\nabla f_\shape^t)^\top \Vb \advar^s \xt ~d\xb + \int_\Omega f_\shape^t  \advar^s \left. \frac{d}{dr}\right|_{r=t^+}\xi^r ~d\xb \rightarrow& \int_\Omega (\nabla f_\shape)^\top \Vb \advar^0 ~d\xb\\ &+ \int_\Omega f_\shape  \advar^0 \di \Vb ~d\xb
 = \partial_t F(t,\advar^0).
 \end{split}
\end{align}
By applying \cite[Lemma 2.16]{Sturm_diss}, we obtain $\nabla f_\shape^t \xt \rightarrow \nabla f_\shape$ in $L^2(\Omega,\mbR^2)$ and $\ f_\shape^t \rightarrow f_\shape$ in $L^2(\Omega)$. Since every ${\Vb \in \vecfields}$ is bounded, we can conclude  $(\nabla f_\shape^t)^\top \Vb \xt \rightarrow \nabla f_\shape^\top \Vb$ in $L^2(\Omega)$. Moreover for every $t \in [0,T)$ the derivative ${\left. \frac{d}{dr}\right|_{r=t^+}\xi^r = \left. \frac{d}{dr}\right|_{r=t^+} \det (I + r D\Vb)}$ is continuous in $r$ and ${\left. \frac{d}{dr}\right|_{r=0^+}\xi^r = \di \Vb}$(see e.g. \cite{Schmidt_diss}), so we derive $f_\shape^t \left. \frac{d}{dr}\right|_{r=t^+}\xi^r \rightarrow f_\shape \di \Vb$ in $L^2(\Omega)$. Again by using \eqref{weak_convergence_lemma}, we obtain the convergence in \eqref{eq:AAM_H1_ell}.\\
Furthermore, we now employ representation \eqref{bilinear_form_representation_1} of the nonlocal bilinear form $\AAA$ to compute the partial derivative of $\AAA$ regarding $t$
\begin{align*}
    &\partial_t \AAA(t,u^0,\advar^s) = \partial_t \AAA(t,u^0,\advar^s) \\
    &= \int_\Omega \int_{\OuO} \advar^s(\xb) \underbrace{\left( u^0(\xb)\nabla_{\xb} \gt_\shape(\xb,\yb) - u^0(\yb) \nabla_{\yb} \gt_\shape(\yb,\xb)\right)^\top \Vb(\xb) \xt(\xb)\xt(\yb)}_{=: A_1(t,u^0)(\xb,\yb)} ~d\yb d\xb \\
    &+ \int_\Omega \int_{\OuO} \advar^s(\xb) \underbrace{\left( u^0(\xb)\nabla_{\yb} \gt_\shape(\xb,\yb) - u^0(\yb) \nabla_{\xb} \gt_\shape(\yb,\xb)\right)^\top \Vb(\yb) \xt(\xb)\xt(\yb)}_{A_2(t,u^0)(\xb,\yb)} ~d\yb d\xb \\
    &+ \int_\Omega \int_{\OuO} \advar^s(\xb) \underbrace{\left( u^0(\xb)\gt_\shape(\xb,\yb) - u^0(\yb) \gt_\shape(\yb,\xb) \right) \left. \frac{d}{dr}\right|_{r=t^+} \left(\xi^r(\xb)\xi^r(\yb)\right)}_{A_3(t,u^0)(\xb,\yb)} ~d\yb d\xb.
\end{align*}
Since $\left. \frac{d}{dr}\right|_{r=t^+} \xi^r(\xb)$ and $\xt(\xb)$ are continuous in ${\xb \in \bar{\Omega}}$, ${\phi_{ij},\nabla \phi_{ij}}$ are essentially bounded for ${(\xb,\yb) \in \bar{\Omega} \times \bar{\Omega}}$ and $\phi_{iI},\nabla \phi_{iI}$ are essentially bounded for ${(\xb,\yb) \in \bar{\Omega} \times \bar{\Omega}_I}$, we can conclude in the same manner as above that $A_i(t,u^0)(\xb,\cdot) \in L^2(\OuO)$ for all $\xb \in \Omega \setminus \Gamma$ and therefore 
\begin{align*}
 \int_\Omega \int_{\OuO} \advar^s(\xb) A_i(t,u^0)(\xb,\yb) ~d\yb d\xb \rightarrow \int_\Omega \int_{\OuO} \advar^0(\xb) A_i(0,u^0)(\xb,\yb) ~d\yb d\xb\quad (i=1,2,3).   
\end{align*}

As a consequence, we derive $\lim_{s,t \searrow 0} \partial_t \AAA(t,u^0,\advar^s) = \partial_t \AAA(0,u^0,\advar^0)$.\\

All in all, we obtain
\begin{align*}
    \lim_{s,t \searrow 0} \partial_t G(t,u^0,\advar^s) &= \lim_{s,t \searrow 0} \partial_t \AAA(t,u^0,\advar^s) - \lim_{s,t \searrow 0} \partial_t F(t,\advar^s) + \lim_{t \searrow 0} \partial_t J(t, u^0)\\
    &=\partial_t \AAA(0,u^0,\advar^0) - \partial_t F(0,\advar^0) + \partial_t J(0, u^0) = \partial_t G(0,u^0,\advar^0).
\end{align*}
\textbf{Case 2: Proof of (H1) for singular kernels}\\
Define $D_n^t := \{(x,y) \in (\OuO)^2: ||F_t(\xb) - F_t(\yb)||_2 > \frac{1}{n} \}$ for $t \in [0,T]$ and $n \in \Nb$. Since $\gt(\xb,\yb) \leq \frac{\gamma^*}{||F_t(\xb) - F_t(\yb)||_2^{2+2s}} < n^{2+2s} \gamma^*$ for all $t \in [0,T],(x,y) \in D_n^t$ and $\xt$ is continuous on $\bar{\Omega} \cup \bar{\Omega}_I$, we can conclude that
\begin{align*}
\iint\limits_{(\OuO)^2} \left( (\advar(\xb) - \advar(\yb))\gamma^{t_l}(\xb,\yb)\xi^{t_l}(\xb)\xi^{t_l}(\yb) \chi_{D_n^{t_l}}(\xb,\yb)\right)^2 ~d\yb d\xb < \infty
\end{align*}
and by using \cite[Lemma 2.16]{Sturm_diss} that
\begin{align*}
&\lim_{l \rightarrow \infty} \iint\limits_{(\OuO)^2} (\advar(\xb) - \advar(\yb))\gamma^{t_l}(\xb,\yb)\xi^{t_l}(\xb)\xi^{t_l}(\yb) \chi_{D_n^{t_l}}(\xb,\yb) ~d\yb d\xb \\ 
&= \iint\limits_{(\OuO)^2} (\advar(\xb) - \advar(\yb))\gamma(\xb,\yb) \chi_{D_n^0}(\xb,\yb) ~d\yb d\xb  .
\end{align*}
With this convergence and \eqref{weak_convergence_lemma}, we derive the second step and with the dominated convergence theorem we get the first and third step of the following computation
\begin{align*}
&\frac{1}{2}\iint\limits_{(\OuO)^2} (\advar(\xb) - \advar(\yb))(q_1(\xb) - q_1(\yb))\g(\xb,\yb) ~d\yb d\xb \\ &= \lim_{n \rightarrow \infty} \frac{1}{2} \iint\limits_{D_n^0} (\advar(\xb) - \advar(\yb))(q_1(\xb) - q_1(\yb))\g(\xb,\yb) ~d\yb d\xb \\ 
&= \lim_{n \rightarrow \infty} \lim_{l \rightarrow \infty} \frac{1}{2} \iint\limits_{D_n^{t_l}} (\advar(\xb) - \advar(\yb))(u^{t_l}(\xb) - u^{t_l}(\yb))\gamma^{t_l}(\xb,\yb) \xi^{t_l}(\xb) \xi^{t_l}(\yb) ~d\yb d\xb \\
&=\lim_{l \rightarrow \infty} \frac{1}{2} \iint\limits_{(\OuO)^2} (\advar(\xb) - \advar(\yb)) (u^{t_l}(\xb) - u^{t_l}(\yb)) \gamma^{t_l}(\xb,\yb) \xi^{t_l}(\xb) \xi^{t_l}(\yb) ~d\yb d\xb \\
&=\lim_{l \rightarrow \infty} \int_{\Omega} f^{t_l} \advar \xi^{t_l} ~d\xb = \int_\Omega f \advar ~d\xb.
\end{align*}
So we can conclude, that $q_1 = u^0$ and $u^t \rightharpoonup u^0$.
Analogously, we can show
\begin{align*}
& \frac{1}{2} \iint\limits_{(\OuO)^2} (u(\xb) - u(\yb))(q_2(\xb) - q_2(\yb))\g(\xb,\yb) ~d\yb d\xb \\
&= \lim_{n \rightarrow \infty} \frac{1}{2} \iint\limits_{D_n^0} (u(\xb) - u(\yb))(q_2(\xb) - q_2(\yb))\g(\xb,\yb) ~d\yb d\xb \\
&= \lim_{k \rightarrow \infty} \lim_{n \rightarrow \infty} \frac{1}{2} \iint\limits_{D_n^{t_k}} (u(\xb) - u(\yb))(\advar^{t_k}(\xb) - \advar^{t_k}(\yb))\gamma^{t_k}(\xb,\yb)\xi^{t_k}(\xb)\xi^{t_k}(\yb) ~d\yb d\xb \\
&= \lim_{k \rightarrow \infty} \frac{1}{2} \iint\limits_{(\OuO)^2} (u(\xb) - u(\yb))(\advar^{t_k}(\xb) - \advar^{t_k}(\yb))\g^{t_k}(\xb,\yb)\xi^{t_k}(\xb)\xi^{t_k}(\yb) ~d\yb d\xb \\
&= - \lim_{k \rightarrow \infty} \int_\Omega \left(\frac{1}{2}(u^{t_k} + u^0) - \bar{u} \right) u \xi^{t_k} ~d\yb d\xb = - \int_\Omega (u^0 - \bar{u}) u ~d\xb 
\end{align*}
and therefore derive $q_2 = \advar^0$ and $\advar^t \rightharpoonup \advar^0$. As in case 1, the next step is to proof
\begin{align*}
\lim_{s,t \searrow 0} \partial_t G(t, u^0,\advar^s) = \partial_t G(0,u^0,\advar^0). 
\end{align*}
By again applying \eqref{weak_convergence_lemma} and the dominated convergence theorem we conclude
\begin{align*}
&\lim_{s,t \searrow 0} \partial_t A(t, u^0,\advar^s) \\ 
&= \lim_{s,t \searrow 0} \frac{1}{2} \iint\limits_{(\OuO)^2} (\advar^s(\xb) - \advar^s(\yb))(u^0(\xb) - u^0(\yb))\gt(\xb,\yb) \left. \frac{d}{dr}\right|_{r=t^+}(\xi^r(\xb) \xi^r(\yb)) ~d\yb d\xb \\
&+\lim_{s,t \searrow 0} \frac{1}{2} \iint\limits_{(\OuO)^2} (\advar^s(\xb) - \advar^s(\yb))(u^0(\xb) - u^0(\yb))(\nabla_x \gt(\xb,\yb) \Vb(\xb) \\
&\hspace{71mm} + \nabla_y \gt(\xb,\yb) \Vb(\yb))\xt(\xb)\xt(\yb) ~d\yb d\xb \\
&= \lim_{n \rightarrow \infty} \lim_{s,t \searrow 0} \frac{1}{2} \iint\limits_{D_n^t}(\advar^s(\xb) - \advar^s(\yb))(u^0(\xb) - u^0(\yb))\gt(\xb,\yb) \left. \frac{d}{dr}\right|_{r=t^+}(\xi^r(\xb) \xi^r(\yb)) ~d\yb d\xb \\
&+ \lim_{n \rightarrow \infty} \lim_{s,t \searrow 0} \frac{1}{2} \iint\limits_{D_n^t} (\advar^s(\xb) - \advar^s(\yb))(u^0(\xb) - u^0(\yb))(\nabla_x \gt(\xb,\yb) \Vb(\xb) \\
&\hspace{75mm} + \nabla_y \gt(\xb,\yb) \Vb(\yb))\xt(\xb)\xt(\yb) ~d\yb d\xb \\
&= \lim_{n \rightarrow \infty} \frac{1}{2} \iint\limits_{D_n^0} (\advar^0(\xb) - \advar^0(\yb))(u^0(\xb) - u^0(\yb))\g(\xb,\yb)(\di \Vb(\xb) + \di \Vb(\yb)) ~d\yb d\xb \\
&+ \lim_{n \rightarrow \infty} \frac{1}{2} \iint\limits_{D_n^0} (\advar^0(\xb) - \advar^0(\yb))(u^0(\xb) - u^0(\yb))(\nabla_x \gamma(\xb,\yb) \Vb(\xb) + \nabla_y \gamma(\xb,\yb) \Vb(\yb)) ~d\yb d\xb \\
&= \frac{1}{2} \iint\limits_{(\OuO)^2} (\advar^0(\xb) - \advar^0(\yb))(u^0(\xb) - u^0(\yb))\g(\xb,\yb)(\di \Vb(\xb) + \di \Vb(\yb)) ~d\yb d\xb \\
&+ \frac{1}{2} \iint\limits_{(\OuO)^2} (\advar^0(\xb) - \advar^0(\yb))(u^0(\xb) - u^0(\yb))(\nabla_x \gamma(\xb,\yb) \Vb(\xb) + \nabla_y \gamma(\xb,\yb) \Vb(\yb)) ~d\yb d\xb\\ 
&= \partial_t A(0,u^0,\advar^0)
\end{align*}
Analogously to case 1, we obtain
\begin{align*}
    \lim_{s,t \searrow 0} \partial_t G(t,u^0,\advar^s) &= \lim_{s,t \searrow 0} \partial_t \AAA(t,u^0,\advar^s) - \lim_{s,t \searrow 0} \partial_t F(t,\advar^s) + \lim_{t \searrow 0} \partial_t J(t, u^0)\\
    &=\partial_t \AAA(0,u^0,\advar^0) - \partial_t F(0,\advar^0) + \partial_t J(0, u^0) = \partial_t G(0,u^0,\advar^0).
\end{align*}
\end{proof}

\pagestyle{plain}
\bibliographystyle{plain}
\bibliography{literature.bib}
\end{document}

%% file: preamble.tex
\usepackage{amsfonts,amsmath,amssymb} 
\usepackage[pdftex]{graphicx} 
\usepackage[utf8]{inputenc}
\usepackage{color}   
\usepackage[dvipsnames]{xcolor}
\usepackage{hyperref}
\usepackage{mathtools}
\usepackage{mathrsfs} 
\usepackage[T1]{fontenc}
\usepackage{lmodern}
\usepackage[hmarginratio=3:2]{geometry} 
\usepackage[ruled, lined, longend, linesnumbered]{algorithm2e}
\SetKwProg{testwhile}{while}{ do}{~~~~~~end while}
\SetKwProg{testif}{if}{ then}{~~~~~end if}
\SetKwProg{testift}{if}{ then}{ }
\usepackage{amsthm}
\usepackage[nottoc,notlot,notlof]{tocbibind}
\usepackage{longtable} 
\usepackage[margin=1cm]{caption} 

\usepackage[rightcaption,raggedright]{sidecap}
\usepackage{framed} 
\usepackage{soul} 
\usepackage{floatpag}
\usepackage{float} 
\usepackage[section]{placeins}
\usepackage{comment}
\usepackage{booktabs,array} 

\usepackage{chngcntr}
\counterwithin{figure}{section}

\allowdisplaybreaks

\usepackage{tikz}
\usetikzlibrary{calc,trees,positioning,arrows,chains,shapes.geometric,%
	decorations.pathreplacing,decorations.pathmorphing,shapes,%
	matrix,shapes.symbols,backgrounds,shadows}

\usepackage{tikz-cd} 
\hypersetup{
	colorlinks=true, 
	linktoc=all,     
	linkcolor=blue,  
	citecolor=green,
	urlcolor=blue
}

\usepackage{sectsty}

\chapternumberfont{\LARGE} 
\chaptertitlefont{\Huge}

\usepackage{fancyhdr}
\fancypagestyle{plain}{ %
	\fancyhf{} 
	
	\fancyfoot[LE,RO]{\thepage}
}

\usepackage{enumitem}
 
\let\oldproof\proof
\renewcommand{\proof}{\color{black}\oldproof}

\newcommand{\shape}{\Gamma}
\newcommand{\shapet}{\Gamma^t}
\newcommand{\gij}{\gamma_{ij}}
\newcommand{\gji}{\gamma_{ji}}
\newcommand{\giI}{\gamma_{iI}}
\newcommand{\ballxy}{\chi_{B_\varepsilon(\xb)}(\yb)}

\DeclareMathOperator{\di}{div}

\DeclareMathOperator{\Emb}{Emb}
\DeclareMathOperator{\Diff}{Diff}
\DeclareMathOperator{\Id}{{\bf Id}}

\newcommand{\R}{\mathbb{R}}

\newcommand{\lk}{\left\{}
\newcommand{\rk}{\right\}}
\newcommand{\defas}{:=}
\newcommand{\g}{\gamma}

\newcommand{\Lu}{\mathcal{L}}
\newcommand{\OuO}{\Omega \cup \Omega_I}
\newcommand{\OtuO}{\Omega^t \cup \Omega_I}

\newcommand{\mcL}{\mathcal{L}}

\newcommand{\mbR}{\mathbb{R}}

\newcommand{\Nb}{\mathbb{N}}

\newcommand{\ee}{e_{1}}

\newcommand{\del}{\delta}

\newcommand{\vecfields}{C_0^k(\Omega,\R^d)} 
\newcommand{\advar}{v} 
\newcommand{\gt}{\gamma^t}
\newcommand{\gijt}{\gamma_{ij}^t}
\newcommand{\gjit}{\gamma_{ji}^t}
\newcommand{\xt}{\xi^t}

\def\ggg{{\phi}}
\def\SSS{{S}}

\def \xb{\mathbf{x}}

\def \nb{\mathbf{n}}
\def \yb{\mathbf{y}}

\def \Vb{\mathbf{V}}
\def \Ub{\mathbf{U}}

\def\AAA{{A}}
\def\A{\mathcal{A}}

\newtheorem{theo}{Theorem}[section]
\newtheorem{lem}[theo]{Lemma}

\newtheorem{ex}[theo]{Example}

\def \be{\begin{equation}}
\def \ee{\end{equation}}
\def \ba{\begin{array}\ds}
	\def \ea{\end{array}}
\def \ds{\displaystyle}

\def\omgs{{\Omega}}
\def\omgc{{\Omega_{I}}}

\def\bfff{{f}}



%% file: interface_nonlocal_2.pdf_tex
\begingroup%
  \makeatletter%
  \providecommand\color[2][]{%
    \errmessage{(Inkscape) Color is used for the text in Inkscape, but the package 'color.sty' is not loaded}%
    \renewcommand\color[2][]{}%
  }%
  \providecommand\transparent[1]{%
    \errmessage{(Inkscape) Transparency is used (non-zero) for the text in Inkscape, but the package 'transparent.sty' is not loaded}%
    \renewcommand\transparent[1]{}%
  }%
  \providecommand\rotatebox[2]{#2}%
  \newcommand*\fsize{\dimexpr\f@size pt\relax}%
  \newcommand*\lineheight[1]{\fontsize{\fsize}{#1\fsize}\selectfont}%
  \ifx\svgwidth\undefined%
    \setlength{\unitlength}{240.75549749bp}%
    \ifx\svgscale\undefined%
      \relax%
    \else%
      \setlength{\unitlength}{\unitlength * \real{\svgscale}}%
    \fi%
  \else%
    \setlength{\unitlength}{\svgwidth}%
  \fi%
  \global\let\svgwidth\undefined%
  \global\let\svgscale\undefined%
  \makeatother%
  \begin{picture}(1,1)%
    \lineheight{1}%
    \setlength\tabcolsep{0pt}%
    \put(0,0){\includegraphics[width=\unitlength,page=1]{interface_nonlocal_2.pdf}}%
    \put(0.46499008,0.49252966){\color[rgb]{0,0,0}\makebox(0,0)[lt]{\lineheight{0}\smash{\begin{tabular}[t]{l}$\Omega_1$\end{tabular}}}}%
    \put(0.51363543,0.08973653){\color[rgb]{0,0,0}\makebox(0,0)[lt]{\lineheight{0}\smash{\begin{tabular}[t]{l}${\color{red}\Gamma}$\end{tabular}}}}%
    \put(0.47918062,0.70223471){\color[rgb]{0,0,0}\makebox(0,0)[lt]{\lineheight{0}\smash{\begin{tabular}[t]{l}$\g_{11}$\end{tabular}}}}%
    \put(0,0){\includegraphics[width=\unitlength,page=2]{interface_nonlocal_2.pdf}}%
    \put(0.09004584,0.82452655){\color[rgb]{0,0,0}\makebox(0,0)[lt]{\lineheight{0}\smash{\begin{tabular}[t]{l}$\Omega_2$\end{tabular}}}}%
    \put(0.53215712,0.93982352){\color[rgb]{0,0,0}\makebox(0,0)[lt]{\lineheight{0}\smash{\begin{tabular}[t]{l}${\tiny \Omega_I}$\end{tabular}}}}%
    \put(0,0){\includegraphics[width=\unitlength,page=3]{interface_nonlocal_2.pdf}}%
    \put(0.1091484,0.37746432){\color[rgb]{0,0,0}\makebox(0,0)[lt]{\lineheight{0}\smash{\begin{tabular}[t]{l}$\g_{22}$\end{tabular}}}}%
  \end{picture}%
\endgroup%